%% file: ay_revised.tex
\documentclass[reqno]{jgokova}
\usepackage{amsmath,amsthm,amssymb}
\usepackage{epsfig,graphicx}

\usepackage{here}
\usepackage{morefloats}
\usepackage{afterpage}
\newtheorem{theorem}{Theorem}[section]
\newtheorem{proposition}[theorem]{Proposition}
\newtheorem{lemma}[theorem]{Lemma}
\newtheorem{corollary}[theorem]{Corollary}
\theoremstyle{definition}
\newtheorem{definition}[theorem]{Definition}
\newtheorem{remark}[theorem]{Remark}

\newtheorem{example}[theorem]{Example}
\newtheorem{question}[theorem]{Question}

\volume{2}

\begin{document}
\input goksty.tex
\title[Corks, Plugs and exotic structures]{Corks, Plugs and exotic structures}
\author[AKBULUT and YASUI]{Selman Akbulut and Kouichi Yasui}
\thanks{The first author is partially supported by NSF, and the second author is partially supported by JSPS Research Fellowships for Young Scientists.}
\address{Department~of~Mathematics, Michigan State University, E.~Lansing, MI, 48824, USA}
\email{akbulut@math.msu.edu}

\address{Department~of~Mathematics, Graduate~School~of~Science, Osaka~University, Toyonaka, Osaka 560-0043, Japan}
\email{kyasui@cr.math.sci.osaka-u.ac.jp}


\begin{abstract}
We discuss corks, and introduce new objects which we call plugs. Though plugs are fundamentally different objects, they also detect exotic smooth structures in 4-manifolds like corks. We discuss relation between corks, plugs and rational blow-downs. 
We show how to detect corks and plugs inside of some exotic manifolds. Furthermore, we construct knotted corks and plugs. 
\end{abstract}
\maketitle

\section{Introduction}
\begin{figure}[ht!]
\begin{center}
\includegraphics[width=2.9in]{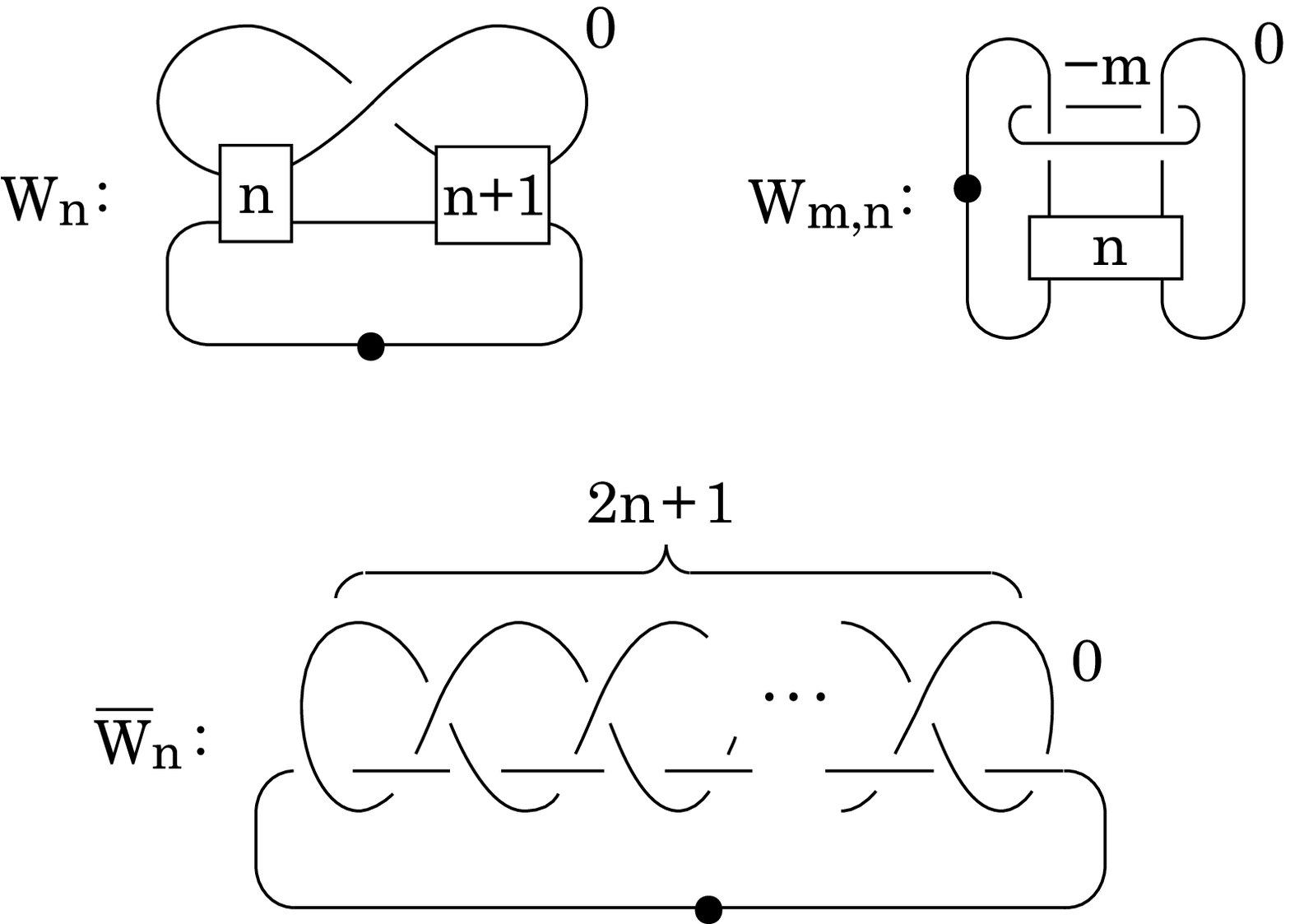}
\caption{}
\label{fig1}
\end{center}
\end{figure}

Let $W_n$, $\overline{W}_n$ and $W_{m,n}$ be the smooth $4$-manifolds given by Figure~\ref{fig1}. Notice that $W_1$ is a version of the manifold defined in  ~\cite{Ma} by Mazur, and  $\overline{W}_n$ is the ``positron'' introduced by the first author and Matveyev~\cite{AM2}. The first author~\cite{A1} proved that $E(2)\#\overline{\mathbf{CP}}^2$ changes its diffeomorphism type by regluing an imbedded copy of $W_1$ inside via a natural  involution on the boundary $\partial W_1$. This was later generalized to  $E(n)\#\overline{\mathbf{CP}}^2$ $(n\geq 2)$  by
 Bi\v zaca-Gompf~\cite{BG}.  The following general theorem was first proved independently by Matveyev~\cite{M}, Curtis-Freedman-Hsiang-Stong~\cite{C}, and later on strengthened by the first author and Matveyev~\cite{AM2}: 
\begin{theorem}[\cite{M}, \cite{C},  \cite{AM2}]\label{th:1.1}
For every homeomorphic but non-diffeomorphic pair of simply connected closed $4$-manifolds, one is obtained from the other by removing a
contractible $4$-manifold and gluing it via an involution on the boundary. Such a contractible $4$-manifold has since been called a {\it Cork}. Furthermore, corks and their complements can always be made compact Stein $4$-manifolds.
\end{theorem}
Also in  \cite{K} Kirby gave a description of the $5$-dimensional $h$-cobordisms induced by corks. So, clearly corks are very important in $4$-manifold topology, but unfortunately not much is known about them. Even though corks determine exotic copies of any manifold, we only know a few concrete examples. Sometime ago the first author posed the question: ``Is $W_1$ a universal cork?", i.e.\ whether $W_1$ is sufficient to detect every exotic structure? This paper came out of our searches for corks in some concrete examples of exotic manifolds. 
 
 \vspace{.05in}

In this work we realized that besides a {\it Cork} there is another basic fundamental object in $4$-manifolds which detects exoticity; we named it a {\it Plug}. Plugs naturally appear when we do surgeries, for example, rational blow-downs. So searching plug structures of $4$-manifolds may be easier than investigating cork structures. Even though plugs are defined in a  similar way to corks, they are different objects and they can't be explained by corks. However, in some cases, plug operations has the same effect as cork operations (Remark~\ref{remark:E_3}). Thus plugs might be helpful to study corks. Similarly to corks, enlarging plugs easily provide us exotic pairs of $4$-manifolds. It turns out that just as corks generalize the Mazur manifold, plugs generalize the ``Gluck construction''.  The lack of understanding how different corks are related to each other have prevented us from defining $4$-manifold invariants from the Stein decomposition theorem of $4$-manifolds \cite{AM2}. We hope that plugs will shed light on understanding corks (e.g.\ they might be deformations of corks). 
   \vspace{.05in}

This paper is organized as follows. In Section $2$, we recall the definition of \textit{Corks}, and then define new objects which we call \textit{Plugs}. Plugs have a property similar to corks and they naturally appear in $0$-log transform and rational blow-down operations. We prove for $n\geq 2$ and $m\geq 1$ that $W_n$ and $\overline{W}_n$ are corks and that $W_{m,n}$ is a plug. 
It is an interesting question whether the sequences of corks $W_n$, $\overline{W}_n$ and the plugs $W_{m,n}$  are sufficient to detect exoticity of all $4$-manifolds. 

In Section $3$, we give examples of $4$-manifolds which contain $W_n$ (resp.\ $W_{m,n}$) as a cork (resp.\ as a plug). For example, $W_1$ is a cork of  elliptic surfaces $E(n)_{p,q}\# \overline{\mathbf{CP}}^2$ $(n\geq 2)$ and knot surgered elliptic surfaces $E(n)_K\# \overline{\mathbf{CP}}^2$ $(n\geq 2)$. 
 In Section $4$ we review the operations of rational blow-down and logorithmic transform, and later in Section 5 we relate them to corks and plugs. This shows that the cork and plug operations can be naturally identified under some conditions.
 As an interesting test case,  in Secion 6,  we draw a handlebody picture of the exotic $\mathbf{CP}^2\# 9\overline{\mathbf{CP}}^2$ without $1$- and $3$-handles, which was recently constructed (and named $E'_3$) by the second author ~\cite{Y0}, \cite{Y1}, and demonstrate how to locate corks and plugs inside. Moreover, by constructions similar to \cite{Y1}, we show that corks and plugs can be knotted, namely, two different embeddings of the same cork (resp.\ plug) can produce two different smooth structures by cork (resp.\ plug) operations. 
  
   \vspace{.05in}

  In a subsequent paper \cite{AY},  we will show that by enlarging corks and plugs in two different ways we can construct pairs of small Stein manifolds that are exotic copies of each other (existence of exotic Stein manifold pairs has  been established recently in \cite{AEMS}).\medskip \\
{\bf Acknowledgements:} The second author would like to thank his adviser Hisaaki Endo and Motoo Tange for drawing a whole handle diagram of the 4-manifold $E_3$ constructed in \cite{Y1} and valuable discussion on 3-manifold topology, respectively. This work was mainly done during his stay at Michigan State University. He is greatful for their hospitality. 
\section{Corks and Plugs}
In this section, we define corks and plugs, and give examples of exotic smooth structures on Stein manifolds.

\begin{definition}
Let $C$ be a compact Stein $4$-manifold with boundary  and $\tau: \partial C\to \partial C$ an involution on the boundary. 
We call $(C, \tau)$ a \textit{Cork} if $\tau$ extends to a self-homeomorphism of $C$, but cannot extend to any self-diffeomorphism of $C$. 
A cork $(C, \tau)$ is called a cork of a smooth $4$-manifold $X$, if  $C\subset X$ and $X$ changes its diffeomorphism type when removing $C$ and re-gluing it via $\tau$. Note that this operation does not change the homeomorphism type of $X$.
\end{definition}
\begin{definition}
Let $P$ be a compact Stein $4$-manifold with boundary and $\tau: \partial P\to \partial P$ an involution on the boundary, which cannot extend to any self-homeomorphism of $P$. We call $(P, \tau)$ a \textit{Plug} of $X$,  if  $P\subset X$ and $X$ keeps its homeomorphism type and changes its diffeomorphism type when removing $P$ and gluing it via $\tau$. 
We call $(P, \tau)$ a \textit{Plug} if there exists a smooth $4$-manifold $X$ such that $(P, \tau)$ is a plug of $X$. 
\end{definition}

\begin{definition}Let $W_n$, $\overline{W}_n$ and $W_{m,n}$ be smooth 4-manifolds in Figure~$\ref{fig1}$. 
Let \linebreak
$f_n:\partial W_n\to \partial W_n$, $\bar{f}_n: \partial \overline{W}_n\to \partial \overline{W}_n$ and $f_{m,n}:\partial W_{m,n}\to \partial W_{m,n}$ be the obvious involutions obtained by first surgering $S^1\times B^3$ to $B^2\times S^2$ in the interiors of $W_n$, $\overline{W}_n$ and $W_{m,n}$, then surgering the other imbedded $B^2\times S^2$ back to $S^1\times B^3$ (i.e. replacing the dots in Figure ~$\ref{fig1}$), respectively. Notice the diagrams of $W_n$, $\overline{W}_n$, and $W_{m,n}$ are symmetric links.
\end{definition}

\begin{remark} Note the following useful facts:\\
(1) $W_1$ is a Mazur manifold, and $(W_1,f_1)$ is equal to $(W,f)$ in \cite{A1}, \cite{AD}. \\
(2) $\overline{W}_n$ is the positron introduced by the first author and Matveyev in ~\cite{AM2}. \\
(3) $W_n$ and $\overline{W}_{n}$ are contractible, whereas $W_{m,n}$ is homotopy equivalent to $S^2$.  \\
(4) $W_n$ $(n\geq 0)$, $\overline{W}_n$ $(n\geq 0)$ and $W_{m,n}$ $(m, n\geq 1)$ are simply connected compact Stein $4$-manifolds. We check this by changing the $1$-handle notations of $W_n$, $\overline{W}_n$ and $W_{m,n}$, and putting the 2-handles in Legendrian positions as in Figure~\ref{fig2}, and checking  the Eliashberg criterion: the framings on the $2$-handles are less than Thurston-Bennequin number. \\
(5)  If $W_{1,0}\subset X$, removing and regluing $W_{1,0}$ to 
via $f_{1,0}$ has an affect of introducing a Gluck twist to $X$ (i.e. cutting out an imbedded copy of $S^2\times D^2$ from $X$, and then regluing by the nontrivial diffeomorphism of $S^2\times S^1$). So far, we only know one example of smooth manifold (which is non-orientable) that becomes exotic by this Gluck operation,  \cite{A3}.
\begin{figure}[h!]
\begin{center}
\includegraphics[width=3.8in]{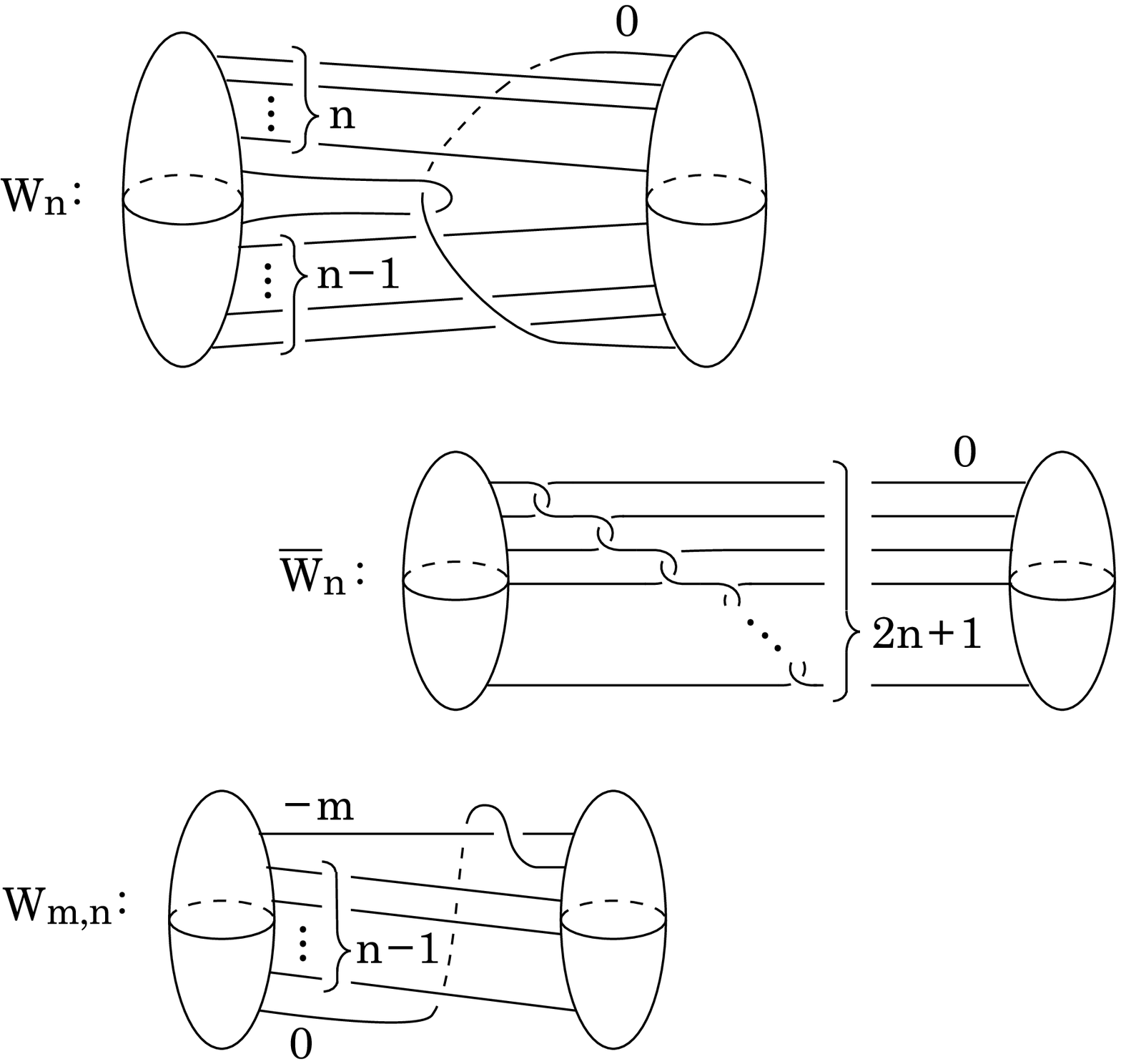}
\caption{}
\label{fig2}
\end{center}
\end{figure}
\end{remark}
We quickly prove the following theorem. 
\begin{theorem}\label{th:cork}
$(1)$ For $n\geq 1$, the pairs $(W_n, f_n)$ and $(\overline{W}_n, \bar{f}_n)$ are corks. \\
$(2)$ For $m\geq 1$ and $n\geq 2$, the pair $(W_{m,n}, f_{m,n})$ is a plug.
\end{theorem}
\begin{remark}
(1) In the case of  $(W_1,f_1)$ this theorem was proved by the first author \cite{A1}. \\
(2) $(W_0, f_0)$ and $(\overline{W}_0, \bar{f}_0)$ are not a corks, because each of the involutions $f_0$  and $\bar{f}_0$ comes from an involution on the boundary of $4$-ball $B^4$, and every diffeomorphism of $S^3$ extends to a self-diffeomorphism of $B^4$. The pair $(W_{m,1}, f_{m,1})$ ($m$: arbitrary integer) is not a plug by the same reason. \\
(3) $(W_{m,0}, f_{m,0})$ ($m$: arbitrary integer) is not a plug because $W_{m,0}=S^2\times B^2$ is not a Stein manifold.
\end{remark}
This theorem clearly follows from the following Lemmas~\ref{lem:topological} and \ref{lem:smooth}. 

\begin{lemma}\label{lem:topological}
$(1)$ For $n\geq 1$, the involution $f_{n}:\partial W_{n}\to \partial W_{n}$ extends to a self homeomorphism of $W_{n}$. \\
$(2)$ For $n\geq 1$, the involution $\bar{f}_n: \partial \overline{W}_n\to \partial \overline{W}_n$ extends to a self homeomorphism of $\overline{W}_n$. \\
$(3)$ For $m\geq 1$ and $n\geq 2$, the involution $f_{m,n}:\partial W_{m,n}\to \partial W_{m,n}$ cannot extend to any self homeomorphism of $W_{m,n}$. 
\end{lemma}

\begin{proof}
$(1)$ (resp.\ $(2)$) is immediate since the boundary of $W_n$ (resp.\ $\overline{W}_n$) is a homology sphere and $W_n$ (resp.\ $\overline{W}_n$) is contractible (see ~\cite{B} for a general discussion). To see $(3)$,  suppose that $f_{m,n}:\partial W_{m,n}\to \partial W_{m,n}$ extends to a self homeomorphism of $W_{m,n}$ and $m\geq 1$ and $n\geq 2$. 
Then the two smooth $4$-manifolds in the Figure~\ref{fig3} have the same intersection form, because one is obtained from the other by removing $W_{m,n}$ and regluing it via $f_{m,n}$. The intersection forms of the left and right $4$-manifolds in Figure~\ref{fig3} are 
\begin{equation*}
\left(
\begin{array}{cc}
-2n-mn^2 &1  \\
1 &-1 
\end{array}
\right) \text{and} 
\left(
\begin{array}{cc}
-2n-mn^2 &-1-mn  \\
-1-mn &-1-m 
\end{array}
\right), 
\end{equation*}
respectively. 
However, we can easily prove that they are not isomorphic. (The right intersection form does not contain any element with self intersection number $-1$. ) This is a contradiction. 
\begin{figure}[ht!]
\begin{center}
\includegraphics[width=2.6in]{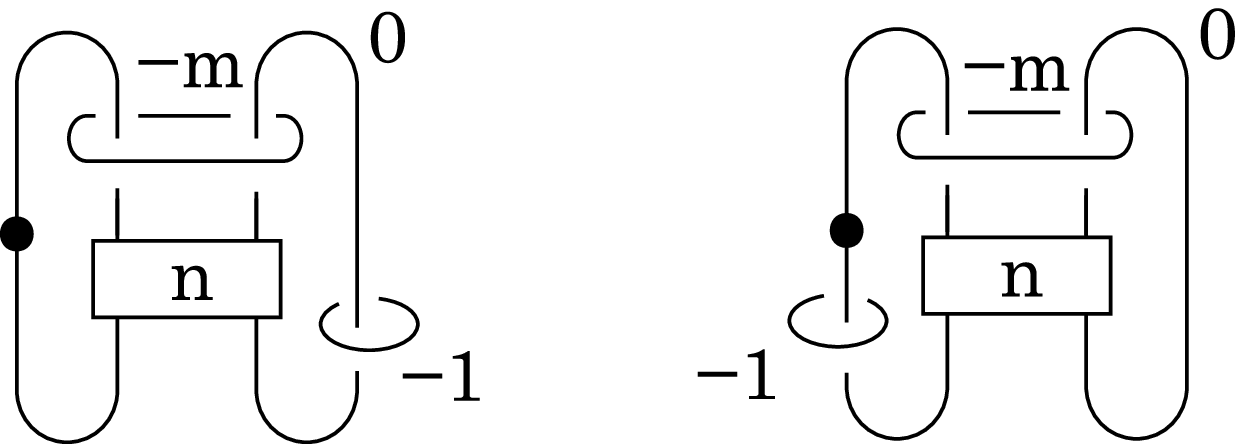}
\caption{}
\label{fig3}
\end{center}
\end{figure}
\end{proof}

The following is a generalization of \cite{A2} and \cite{AM1}, it basically says that by enlarging $W_n$, ${\overline{W}}_n$ and $W_{m,n}$ we can obtain exotic manifolds (hence these manifolds are corks and plugs of their enlargements). 

\vspace{.05in}

\begin{lemma}\label{lem:smooth} 
$(1)$ Let $W^1_n$ and $W^2_n$ be simply connected compact smooth $4$-manifolds defined in Figure~$\ref{fig4}$. Then $W^1_n$ and $W^2_n$ $(n\geq 1)$ are homeomorphic but not diffeomorphic.  In particular, the involution $f_{n}:\partial W_{n}\to \partial W_{n}$ cannot extend to any self diffeomorphism of $W_n$. Note that $W^2_n$ is obtained from $W^1_n$ by regluing $W_n$ via $f_n$.

\begin{figure}[ht!]
\begin{center}
\includegraphics[width=3.9in]{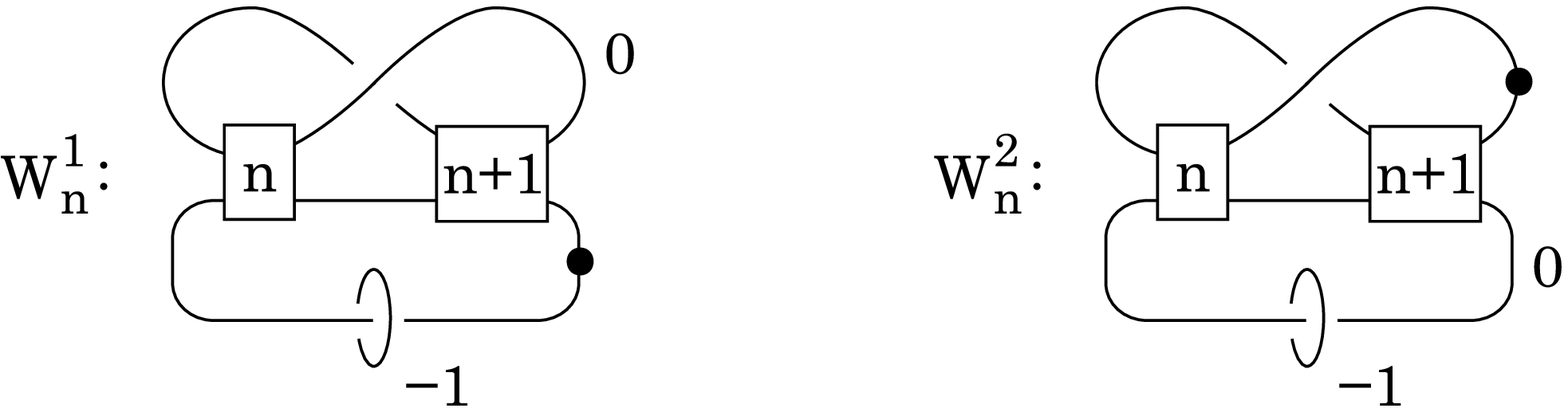}
\caption{}
\label{fig4}
\end{center}
\end{figure}$\:$
\\
\\
$(2)$ Let ${\overline{W}}^1_n$ and ${\overline{W}}^2_n$ be simply connected compact smooth $4$-manifolds defined in Figure~$\ref{fig5}$. Then ${\overline{W}}^1_n$ and ${\overline{W}}^2_n$ $(n\geq 1)$ are homeomorphic but not diffeomorphic. In particular, the involution ${\overline{f}}_{n}:\partial {\overline{W}}_{n}\to \partial \overline{W}_{n}$ cannot extend to any self diffeomorphism of $\overline{W}_n$. Note that $\overline{W}^2_n$ is obtained from $\overline{W}^1_n$ by regluing $\overline{W}_n$ via $\overline{f}_n$. 

\begin{figure}[ht!]
\begin{center}
\includegraphics[width=2.3in]{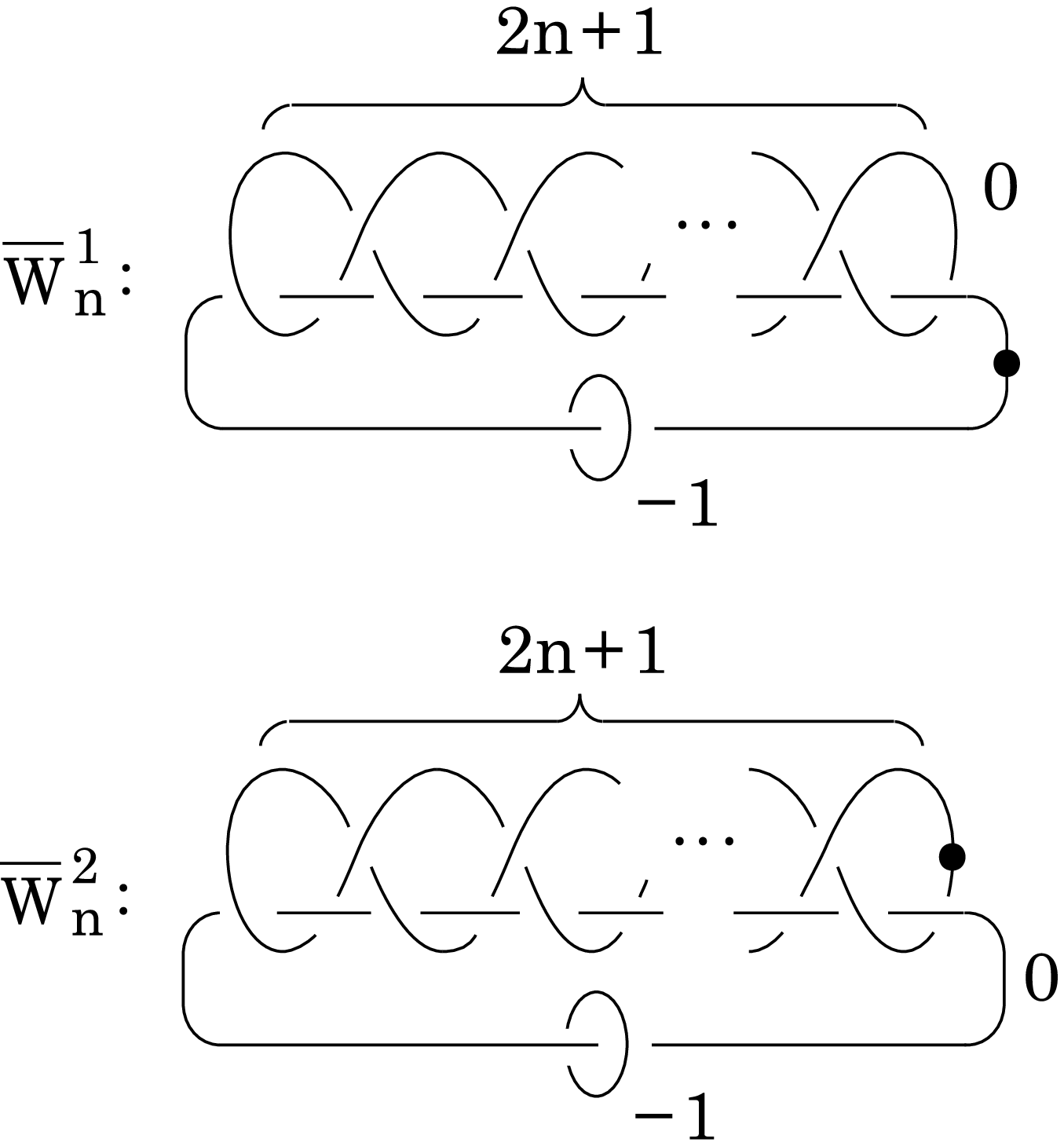}
\caption{}
\label{fig5}
\end{center}
\end{figure}$\:$
\\

\vspace{.05in}

$(3)$ Let $W^1_{m,n}$ and $W^2_{m,n}$ be simply connected compact smooth $4$-manifolds defined in Figure~$\ref{fig6}$. Then $W^1_{m,n}$ and $W^2_{m,n}$ $(m\geq 1, n\geq 2)$ are homeomorphic but not diffeomorphic. Note that $W^2_{m,n}$ is obtained from $W^1_{m,n}$ by regluing $W_{m,n}$ via $f_{m,n}$.
\\

\begin{figure}[ht!]
\begin{center}
\includegraphics[width=4.2in]{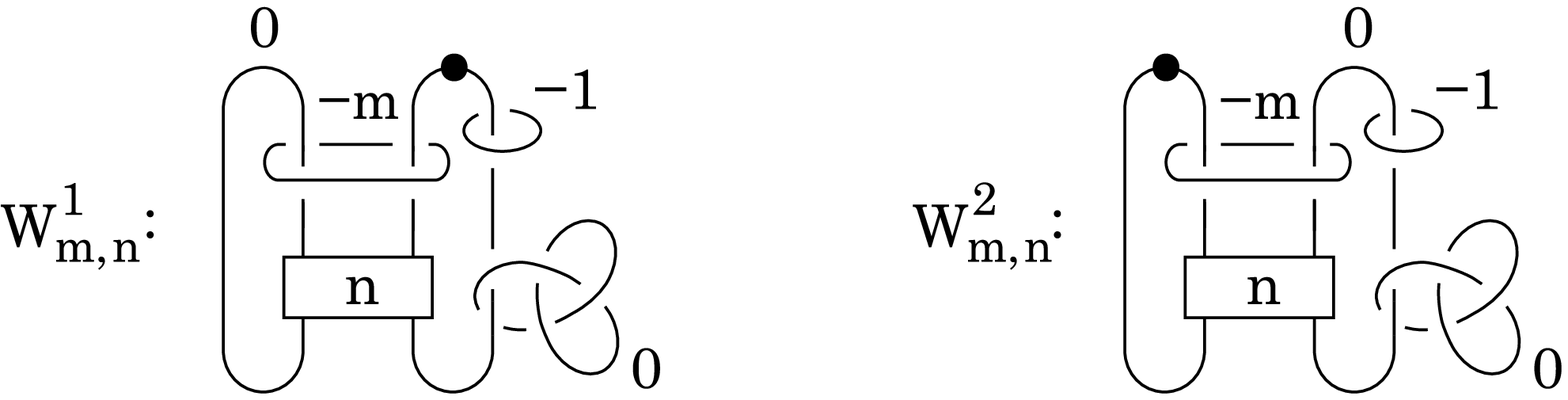}
\caption{}
\label{fig6}
\end{center}
\end{figure}
\end{lemma}
\begin{proof}

(1) Lemma~\ref{lem:topological}.(1) implies that $W^1_n$ and $W^2_n$ are homeomorphic. We can easily check that $W^1_n$ is a compact Stein manifold by an isotopy in Figure~\ref{fig4} (it follows from the Eliashberg criterion: The framing on the $2$-handles are less or equal to the Thurston-Benequin framing). Since every compact Stein manifold can be embedded into a minimal symplectic 4-manifold (e.g. \cite{LM}), $W^1_n$ has no 2-sphere with self intersection number $-1$. But clearly $W^2_n$ contains a 2-sphere with self intersection number $-1$. In particular, $W^2_n$ is not a compact Stein manifold. Hence $W^1_n$ and $W^2_n$ $(n\geq 1)$ are not diffeomorphic. 

(2) Similar to (1), we can show that ${\overline{W}}^1_n$ is a compact Stein manifold and that ${\overline{W}}^2_n$ is not a compact Stein manifold. Therefore ${\overline{W}}^1_n$ and ${\overline{W}}^2_n$ are not diffeomorphic. 

(3) As before, similar to (1), we can show that $W^1_{m,n}$ is a compact Stein manifold  and that $W^2_{m,n}$ is not a compact Stein manifold. Therefore, $W^1_{m,n}$ and $W^2_{m,n}$ $(n\geq 1)$ are not diffeomorphic. We check they are homeomorphic. The intersection forms of $W^1_{m,n}$ and $W^2_{m,n}$ are 
\begin{equation*}
\left(
\begin{array}{ccc}
-1 &0 &0  \\
0 &-m &1 \\
0 &1 &0 
\end{array}
\right) \text{and} 
\left(
\begin{array}{ccc}
-2n-n^2m &1 &1  \\
1 &0 &0 \\
1 &0 &-1 
\end{array}
\right),
\end{equation*}
respectively. We can easily show that they are isomorphic to $\langle 1 \rangle \oplus \langle -1 \rangle \oplus \langle -1 \rangle $. It thus follows from Boyer's theorem~\cite{B} that $W^1_{m,n}$ and $W^2_{m,n}$ are homeomorphic. Note that their boundaries are diffeomorphic. 
\end{proof}
\begin{remark}
Here we showed $W_n$, $\overline{W}_n$ and $W_{m,n}$ are corks and plugs by using a quick argument of the first author and Matveyev~\cite{AM2} based on property of Stein manifolds. Notice that this technique gives $W_{n}^{i}$ are $W_{m,n}^{i}$ are Stein only when $i=1$, whereas when  $i=2$ they are not Stein (they have an imbedded $-1$ sphere). In  a subsequent paper \cite{AY}, by a more sophisticated argument, by enlarging corks and plugs we construct  simply connected exotic compact Stein manifold pairs  realizing any  Betti number $b_{2}\geq 1$.
\end{remark}

\vspace{.05in}
The technical proposition below makes arguments of stabilizations easy. We use this proposition in Section~\ref{section:relation}.

\vspace{.2in}
 
\begin{proposition}\label{prop:stabilization}
$(1)$ Let $X$ be a simply connected compact smooth $4$-manifold which contains $W_{n}$ $($resp.\ $\overline{W}_n$$)$. Let $Y$ be the simply connected compact smooth $4$-manifold obtained from $X$ by removing $W_{n}$ $($resp.\ $\overline{W}_n$$)$ and regluing it via $f_{n}$ $($resp.\ $\overline{f}_n$$)$. Then $X\# S^2\times S^2$ is diffeomorphic to $Y\# S^2\times S^2$. \\

$(2)$ Let $X$ be a compact smooth $4$-manifold which contains $W_{m,n}$. Let $Y$ be the compact smooth $4$-manifold obtained from $X$ by removing $W_{m,n}$ and regluing it via $f_{m,n}$. If $m$ is even $(\text{resp.\ odd})$, then $X\# S^2\times S^2$ $(\text{resp.}$ $X\# \mathbf{C}\mathbf{P}^2\# \overline{\mathbf{CP}}^2)$ is diffeomorphic to $Y\# S^2\times S^2$ $(\text{resp.}$ $Y\# \mathbf{C}\mathbf{P}^2\# \overline{\mathbf{CP}}^2)$. 
\end{proposition}
\begin{proof}
(1) $X\# S^2\times S^2$ is obtained from $X$ by surgering a homotopically trivial loop (with the correct framing). We choose to surger $X$ along a meridian of the unique dotted circle of $W_n$ (resp.\ $\overline{W}_n$). This corresponds to turning the unique $1$-handle of $W_n$ (resp.\ $\overline{W}_n$) into a $0$-framed $2$-handle (i.e. replace the dot with 0-framing). Similarly, $Y\# S^2\times S^2$ is obtained from $Y$ by replacing the dot of the dotted circle of $W_n$ (resp.\ $\overline{W}_n$) with $0$. It thus follows from the definition of $Y$ that $X\# S^2\times S^2$ and $Y\# S^2\times S^2$ have the same handlebody diagrams. \\
(2) Since the dotted circle of $W_{m,n}$ has $-m$-framed meridian, changing the dotted circle into $0$-framed circle corresponds to a connected sum with $S^2\times S^2$. This fact implies the required claim. 
\end{proof}

\section{Examples}

\vspace{.05in}

Here, by improving an argument in Gompf-Stipsicz~\cite[Section $9.3$]{GS}, we give examples of closed manifolds containing $W_n$ and $\overline{W}_n$ (resp.\ $W_{m,n}$) as corks (resp.\  plugs).

\vspace{.05in}

Let $E(n)$ be the relatively minimal simply connected elliptic surface with Euler characteristic $12n$ and with no multiple fibers, and $E(n)_{p_1,\dots,p_k}$ the elliptic surface obtained from $E(n)$ by performing logarithmic transformations of multiplicities $p_1,\dots,p_k$. Now recall.

\vspace{.05in}

\begin{theorem}[Gompf-Stipsicz~\cite{GS}]\label{th:GS}
For $n\geq 1$, the elliptic surface $E(n)$ has the handle decomposition as in Figure~$\ref{fig7}$. The obvious cusp neighborhood $($i.e. the dotted circle, $-1$-framed meridian of the dotted circle, and the left most $0$-framed unknot$)$ is isotopic to the regular neighborhood of a cusp fiber of $E(n)$. 
\begin{figure}[ht!]
\begin{center}
\includegraphics[width=3.5in]{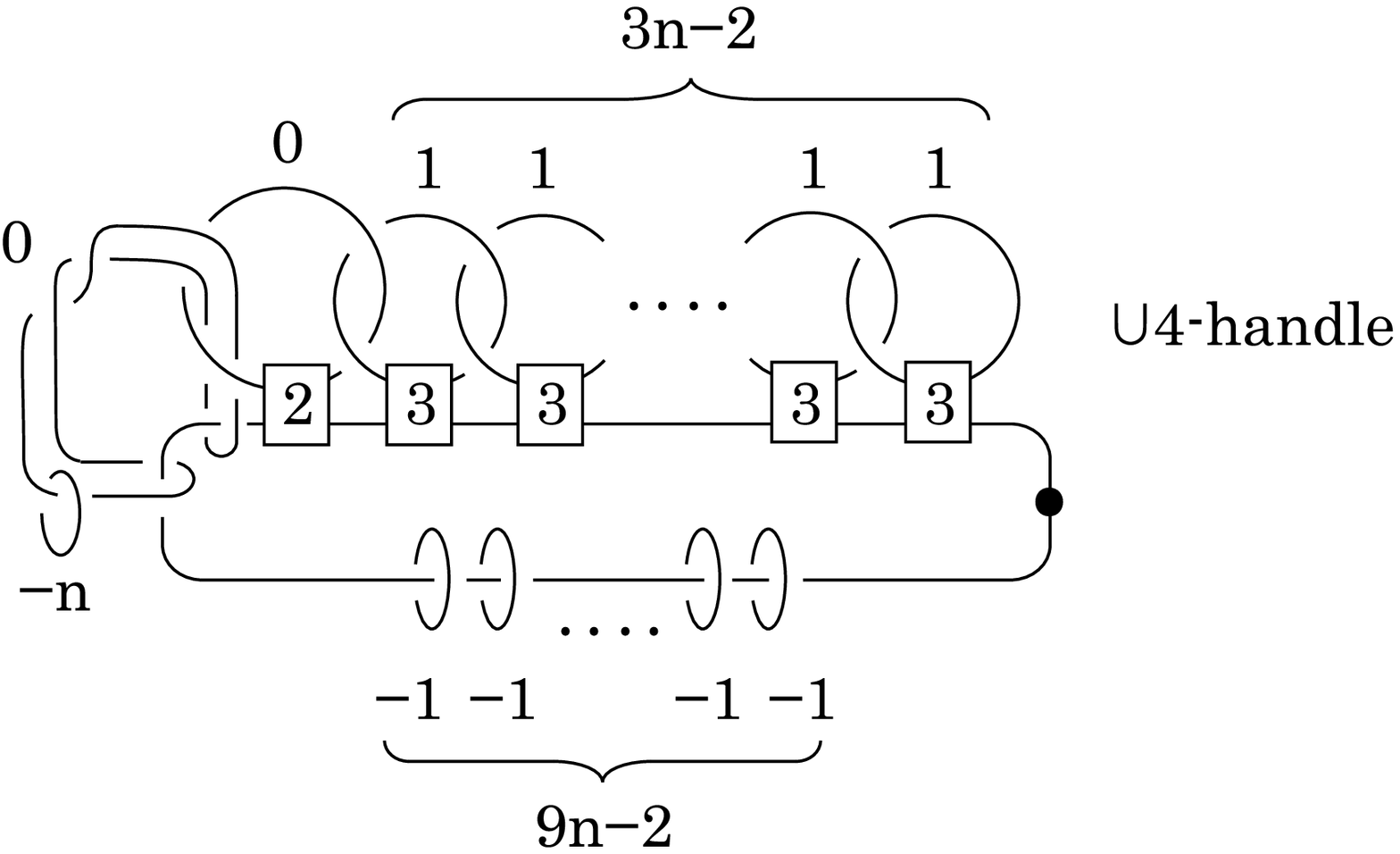}
\caption{$E(n)$}
\label{fig7}
\end{center}
\end{figure}
\end{theorem}

\vspace{.05in}

\begin{corollary}\label{cor:GS}
For $n\geq 1$, the elliptic surface $E(n)$ has the handle decomposition as in Figure~$\ref{fig8}$. The obvious cusp neighborhood $($i.e. $0$-framed trefoil knot$)$ is isotopic to the regular neighborhood of a cusp fiber of $E(n)$. 
\begin{figure}[ht!]
\begin{center}
\includegraphics[width=3.8in]{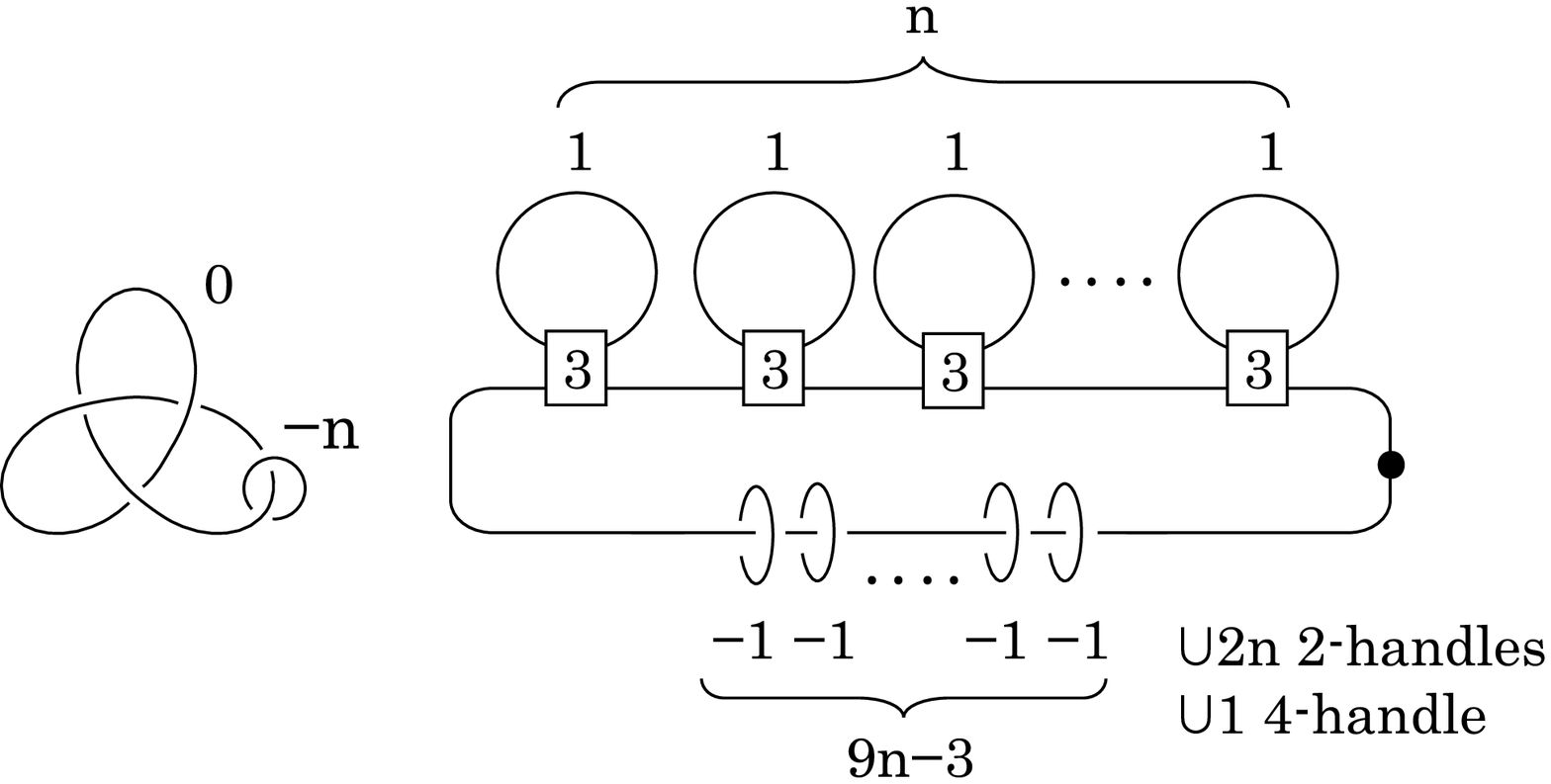}
\caption{$E(n)$}
\label{fig8}
\end{center}
\end{figure}
\end{corollary}
\begin{proof}
In Figure~$\ref{fig7}$, pull off the leftmost 0-framed unknot from the dotted circle by sliding over $-1$-framed knot. 
\end{proof}

\vspace{.1in}

\begin{proposition}\label{prop:cork}
For $n\geq 1$ and $m\geq 1$ we have\\

$(1)$ $E(2n)\# \overline{\mathbf{CP}}^2$ has the handle decompositions as in Figure~$\ref{fig9}$. \\

$(2)$ $E(2n)\# m\overline{\mathbf{CP}}^2$ has the handle decompositions as in Figure~$\ref{fig10}$.\\

$(3)$ $E(2n+1)\# \overline{\mathbf{CP}}^2$ has the handle decomposition as in Figure~$\ref{fig11}$.
\\

\begin{figure}[ht!]
\begin{center}
\includegraphics[width=3.9in]{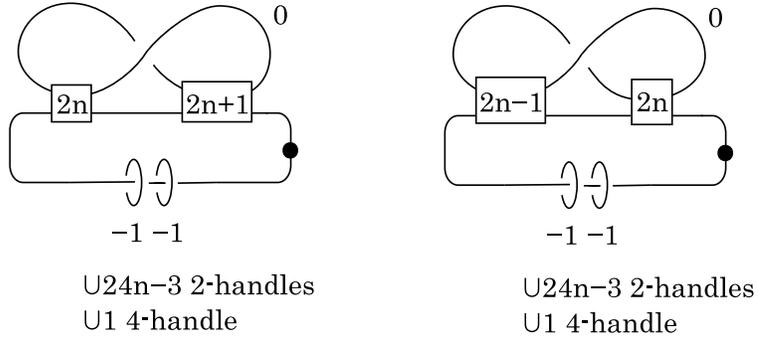}
\caption{Two decompositions of $E(2n)\# \overline{\mathbf{CP}}^2$}
\label{fig9}
\end{center}
\end{figure}

\begin{figure}[ht!]
\begin{center}
\includegraphics[width=3.8in]{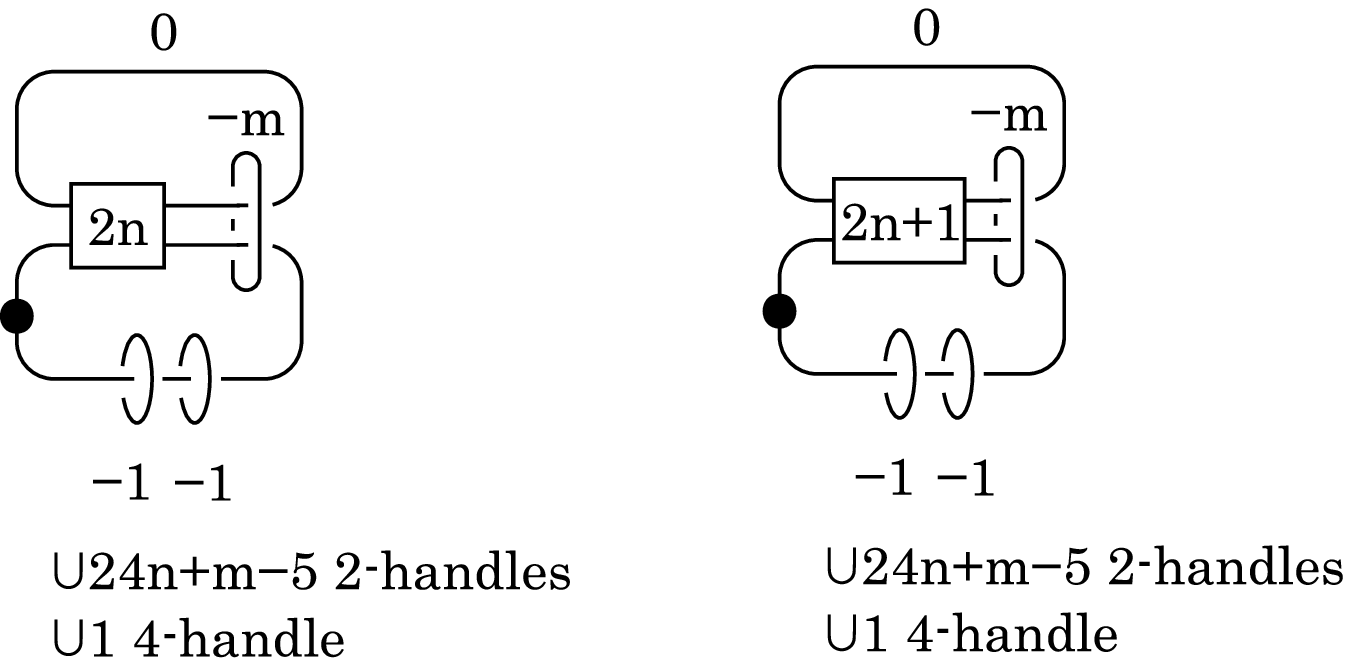}
\caption{Two decompositions of $E(2n)\# m\overline{\mathbf{CP}}^2\:(m\geq 1)$}
\label{fig10}
\end{center}
\end{figure}

\begin{figure}[ht!]
\begin{center}
\includegraphics[width=2.7in]{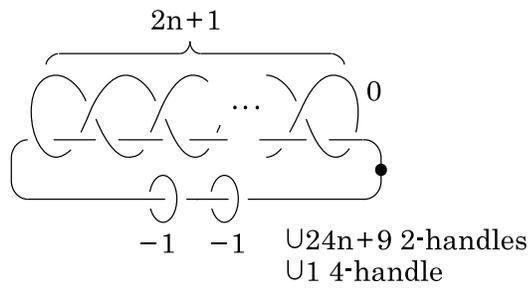}
\caption{$E(2n+1)\# \overline{\mathbf{CP}}^2$ ($n\geq 1$)}
\label{fig11}
\end{center}
\end{figure}
\end{proposition}

\vspace{.1in}

\begin{proof}
$(1)$ Starting with Figure~\ref{fig8}, the handle slides as in Figure~$\ref{fig26}$ give Figure~$\ref{fig27}$. We  easily get Figure~\ref{fig28} by handle slides and isotopies to the $E(2n)$ case of Figure~$\ref{fig27}$. By blowing up, we get Figure~\ref{fig29}. We have the left diagram in Figure~$\ref{fig9}$ by a slide. We also get Figure~\ref{fig30} by a slide in Figure~\ref{fig29}. Notice that the Figure~\ref{fig30} is isotopic to the right diagram in Figure~$\ref{fig9}$. 
\vspace{.05in}

$(2)$ In Figure~\ref{fig29}, blowing up the left (resp.\ right) middle $-1$-framed circle $m-1$ times gives the left (resp.\ right) diagram in Figure~\ref{fig10}. 

\vspace{.05in}

$(3)$ It follows from a diagram of $E(n)$ in Figure~\ref{fig27} that $E(2n+1)$ admits a diagram as in Figure~\ref{fig31}. We obtain a diagram of $E(2n+1)$ in Figure~\ref{fig36} by handle slides and isotopies as in Figure~$\ref{fig31}\sim \ref{fig36}$. By blowing up, we get Figure~\ref{fig11}.  
\end{proof}

This proposition implies the following theorem which says that many exotic 4-manifolds admit  $W_n$ and ${\overline{W}}_n$ (resp.\ $W_{m,n}$) as corks (resp.\ plugs). 

\vspace{.05in}

\begin{theorem}
\label{ex:corks and plugs}
For each $n\geq 1$ and $m\geq 1$ we have\\

$(1)$ $(W_{2n-1}, f_{2n-1})$ and $(W_{2n},f_{2n})$ are corks of $E(2n)\# \overline{\mathbf{CP}}^2$.\\

$(2)$  $(W_{m,2n}, f_{m,2n})$ and $(W_{m,2n+1}, f_{m,2n+1})$ are plugs of $E(2n)\# m\overline{\mathbf{CP}}^2$. \\

$(3)$ $(\overline{W}_{n},\overline{f}_n)$ is a cork of $E(2n+1)\# \overline{\mathbf{CP}}^2$.
\end{theorem}

\begin{proof}
(1) Figure~$\ref{fig9}$ implies that $E(2n)\# \overline{\mathbf{CP}}^2$ splits off $2\overline{\mathbf{CP}}^2$ by regluing $W_{2n}$ (resp.\ $W_{2n-1}$) via $f_{2n}$ (resp.\ $W_{2n-1}$). The Seiberg-Witten invariant and the blow-up formula imply that $E(2n)\# \overline{\mathbf{CP}}^2$ cannot split off $2\overline{\mathbf{CP}}^2$,  therefore $(W_{2n-1}, f_{2n-1})$ and $(W_{2n},f_{2n})$ are corks of $E(2n)\# \overline{\mathbf{CP}}^2$. 
(2) and (3) are similar to (1).
\end{proof}

Let $E(n)_{p,q}$ denote the $4$-manifold obtained from $E(n)$ by the logarithmic trasform on two disjoint torus fibers. Also, for a knot $K$ in $S^3$ with a nontrivial Alexander polynomial, let $E(n)_K$ be the $4$-manifold obtained from $E(n)$ by Fintushel-Stern's knot surgery 
(\cite{FS2}) with $K$ in the regular neighborhood of a cusp fiber. Then we have:
 
 \vspace{.05in}
 
\begin{theorem}
For each $p,q\geq 1$ and  $\gcd(p,q)=1$ we have\\

$(1)$ $W_1$ $(resp.\ \overline{W}_1)$ is a cork of $E(n)_{p,q}\# \overline{\mathbf{CP}}^2$ $(n\geq 2)$ $[resp.\ n\geq 3]$.\\

$(2)$ $W_{1,3}$ is a plug of $E(n)_{p,q}\# \overline{\mathbf{CP}}^2$ $(n\geq 2)$ $(\text{also valid  for $n=1$ and $p,q\geq 2$})$. \\

$(3)$ $W_1$ $(resp.\ \overline{W}_1)$ is a cork of $E(n)_K\# \overline{\mathbf{CP}}^2$ $(n\geq 2)$ $[resp.\ n\geq 3]$. \\

$(4)$ $W_{1,3}$ is a plug of $E(n)_K\# \overline{\mathbf{CP}}^2$ $(n\geq 1)$. 
\end{theorem}

\vspace{.02in}

\begin{proof}
Start from Figure~$\ref{fig27}$. We can easily construct $W_1$ (resp.\ $W_{1,3}$) in $E(n)\# \overline{\mathbf{CP}}^2$ when $n\geq 2$ (resp.\ $n\geq 1$)  away from the regular neighborhood of a cusp fiber of $E(n)$, by an argument similar to the proof of Proposition~$\ref{prop:cork}$. Perform logarithmic transformations of multiplicities $p$ and $q$ (resp.\ knot surgery along $K$) in the cusp neighborhood. Then an argument similar to the proof of Theorem ~$\ref{ex:corks and plugs}$ shows the required claims. The case of $\overline{W}_1$ is almost the same. 
\end{proof}

\section{Rational blow-down and logarithmic transform}\label{sec:rbd}
In this section, we review the rational blow-down introduced by Fintushel-Stern~\cite{FS1} and logarithmic transform (see also Gompf-Stipsicz~\cite{GS}). We also give a new procedure to change handlebody diagrams by rational blow-downs and logarithmic transforms. 

\begin{figure}[ht!]
\begin{center}
\includegraphics[width=3.8in]{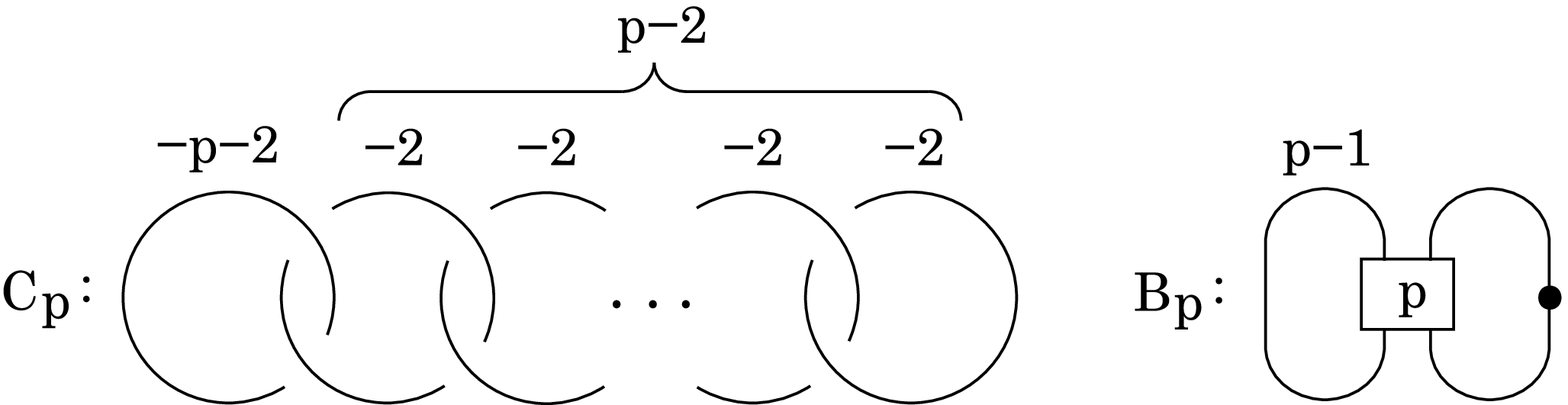}
\caption{}
\label{fig12}
\end{center}
\end{figure}

Let $C_p$ and $B_p$ ($p\geq 2$) be the smooth $4$-manifolds in Figure~$\ref{fig12}$. The boundary $\partial C_p$ of $C_p$ is diffeomorphic to the lens space $L(p^2,p-1)$ and to the boundary $\partial B_p$ of $B_p$. 

\begin{figure}[ht!]
\begin{center}
\includegraphics[width=1.75in]{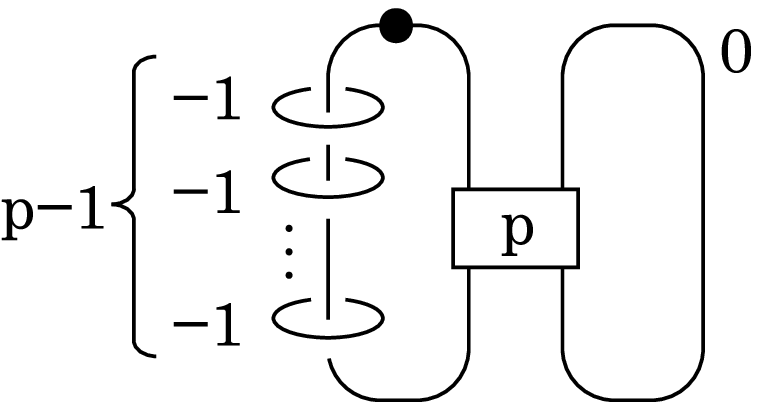}
\end{center}
\caption{$C_p$}
\label{fig13}
\end{figure}

\begin{figure}[ht!]
\begin{center}
\includegraphics[width=4.9in]{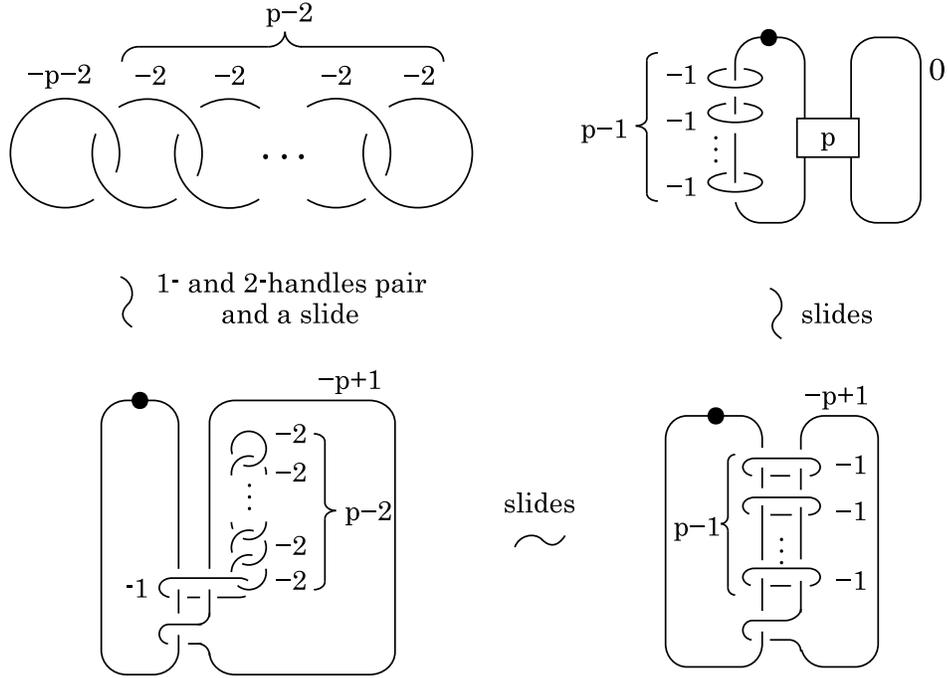}
\end{center}
\caption{diagrams of $C_p$}
\label{fig14}
\end{figure}

\begin{figure}[ht!]
\begin{center}
\includegraphics[width=4.2in]{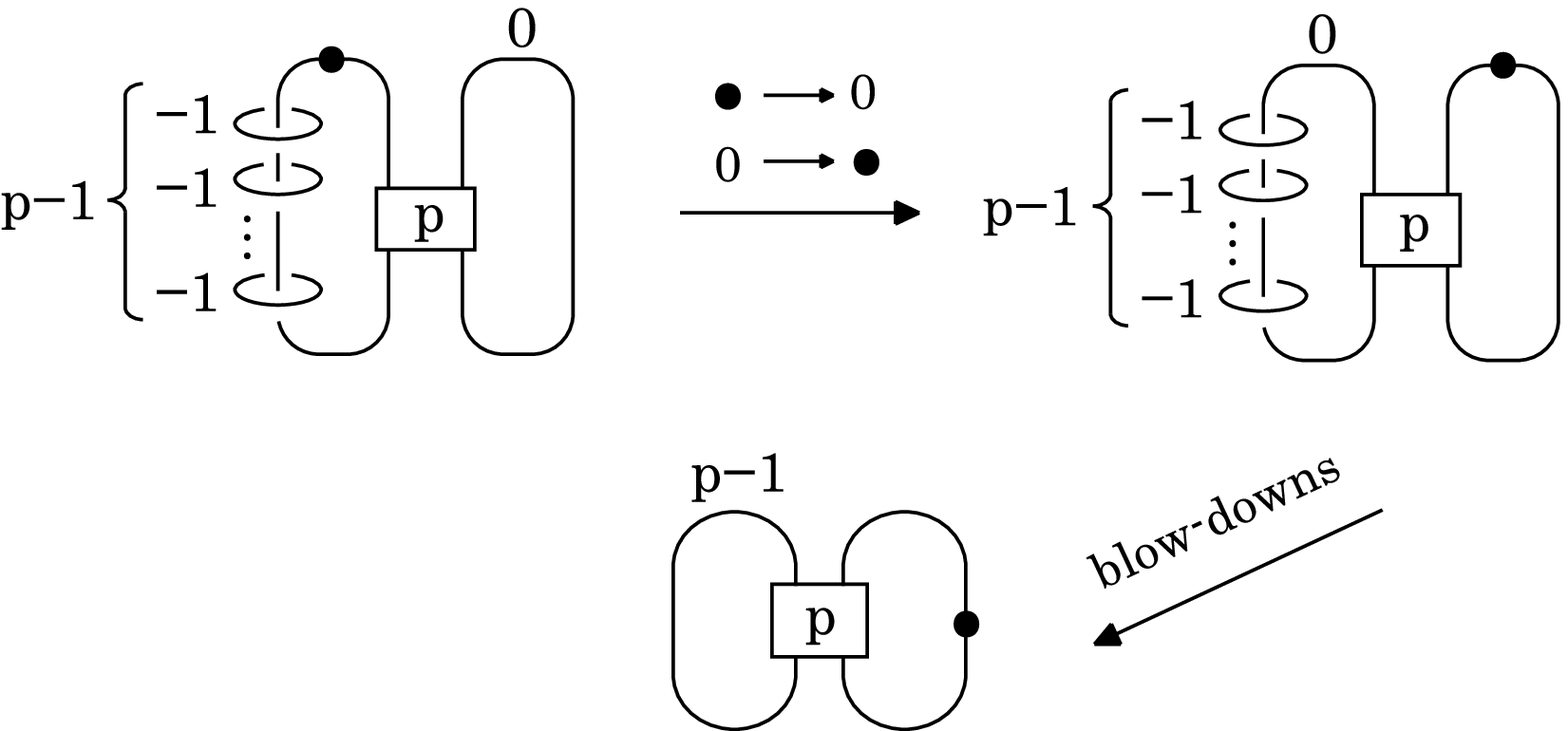}
\end{center}\medskip 
\caption{Rational blow-down along $C_p$}
\label{fig15}
\end{figure}

Suppose that $C_p$ embeds in a smooth $4$-manifold $X$. Let $X_{(p)}$ be a smooth $4$-manifold obtained from $X$ by removing $C_p$ and gluing $B_p$. The $4$-manifold $X_{(p)}$ is called the rational blow-down of $X$ along $C_p$. Note that $X_{(p)}$ is uniquely determined up to diffeomorphism by a fixed pair $(X,C_p)$. 
We briefly recall the procedure given in \cite[Section 8.5]{GS} to obtain a handlebody diagram of $X_{(p)}$ from a diagram of $X$. Change the diagram of $C_p$ in Figure~$\ref{fig12}$ into the diagram in Figure~$\ref{fig13}$ (introduce a 1-handle/2-handle pair and slide handles as shown in Figure~$\ref{fig14}$). Then surger the obvious $S^1\times B^3$ to $B^2\times S^2$ and surger the other $B^2\times S^2$ to $S^1\times B^3$ (i.e. change the dot and $0$ in Figure~$\ref{fig13}$). Finally blow down $-1$-framed unknots. See Figure~$\ref{fig15}$ for this procedure. 

\begin{figure}[ht!]
\begin{center}
\includegraphics[width=5.0in]{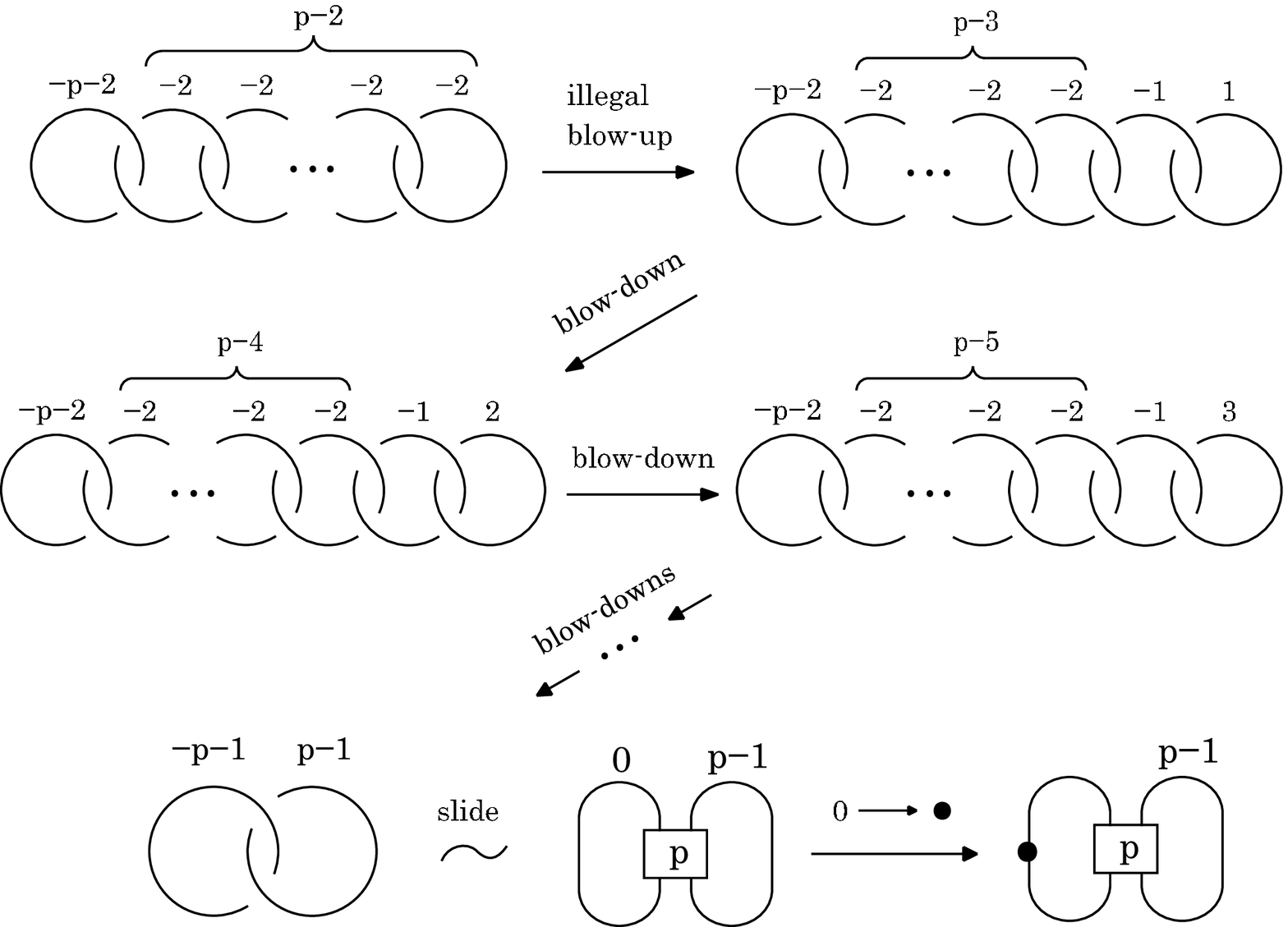}
\end{center}
\caption{Rational blow-down along $C_p$}
\label{fig16}
\end{figure}

We now give a new way to draw a handlebody diagram of $X_{(p)}$, as indicated in \linebreak Figure~$\ref{fig16}$: First we connected sum with $\mathbf{C}\mathbf{P}^2$ and get the second diagram. Now we have a $-1$-framed unknot, then blow down this $-1$-framed unknot. Now again we have a $-1$-framed unknot, by repeating blow-downs as before, we get the fifth diagram. A handle slide gives the sixth diagram. Finally we replace $0$ with a dot. Since this procedure changes $C_p$ into $B_p$, this is the rational blow-down operation. 

Next, we discuss the logarithmic transform. Let $\varphi_p$ $(p\geq 0)$ be the self-diffeomorphism of $T^3$ induced by the automorphism 
\begin{equation*}
\left(
\begin{array}{ccc}
1 &0 &0  \\
0 &0 &1  \\
0 &-1 &p 
\end{array}
\right)
\end{equation*}
of $H_1(S^1;\mathbf{Z})\oplus H_1(S^1;\mathbf{Z})\oplus H_1(S^1;\mathbf{Z})$ with the obvious basis. 
Suppose that a smooth $4$-manifold $X$ contains a torus $T$ with trivial normal bundle $\nu (T)\approx T\times B^2$. 

Now we remove $\text{int}\:\nu (T) $ from $X$ and glue $T^2\times B^2$ back via $\varphi_p: T^2\times S^1\to \partial \nu (T)$. We call this operation a logarithmic transform of $X$ with multiplicity $p$. Note that this definition is consistent with the one given in ~\cite{GS} when $T$ is in a cusp neighborhood (i.e. $B^4$ with a $2$-handle attached along a $0$-framed right trefoil knot). 

\vspace{.05in}

\begin{figure}[ht!]
\begin{center}
\includegraphics[width=4.3in]{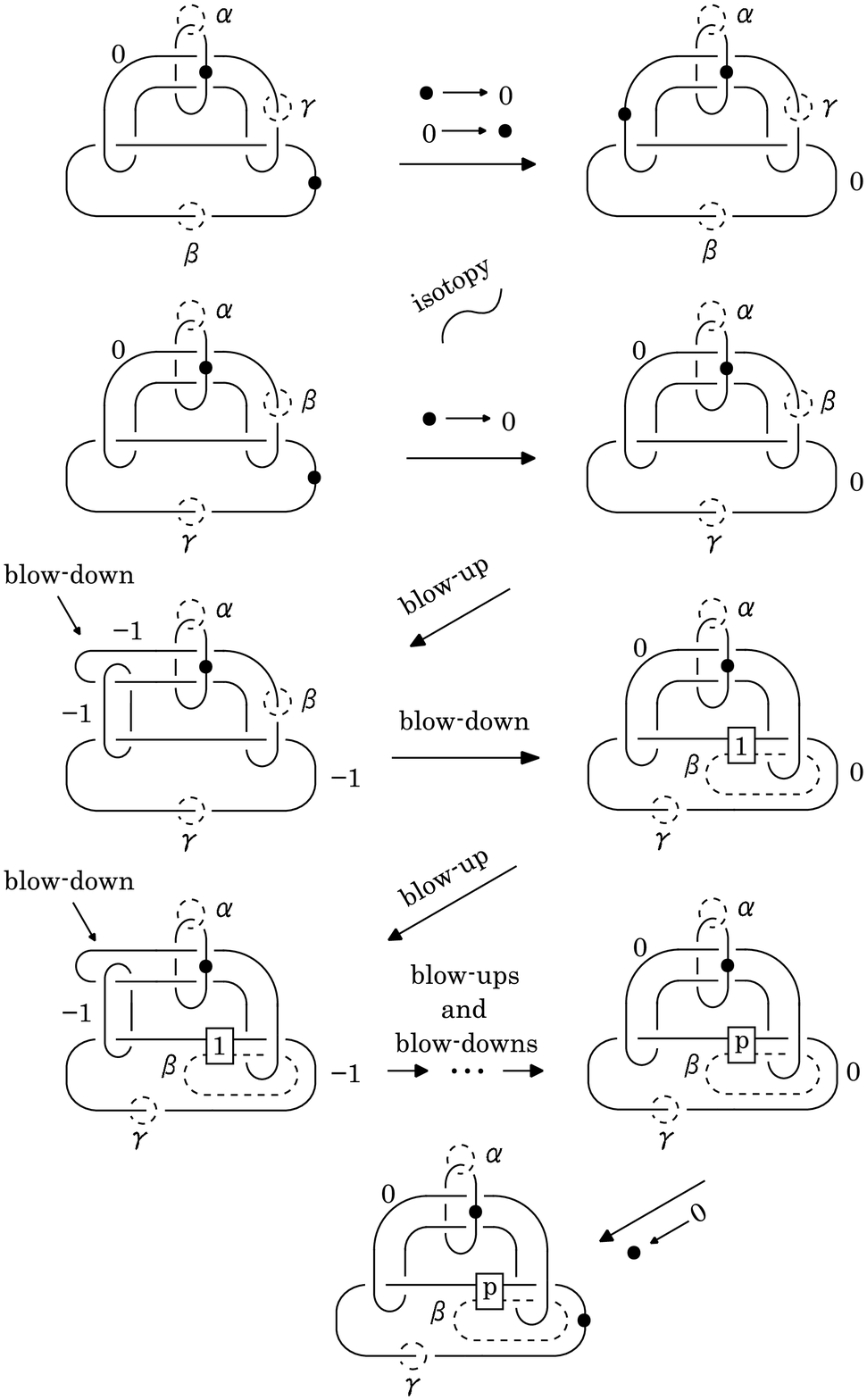}
\caption{Logarithmic transform with multiplicity $p$}
\label{fig17}
\end{center}
\end{figure}

Let $X_p$ be a logarithmic transform of $X$ with multiplicity $p$. We can draw a handlebody diagram of $X_p$ as indicated in Figure~$\ref{fig17}$. The first diagram denotes $\nu (T)\approx T\times B^2$ in $X$. We will keep track of circles $\alpha, \beta, \gamma$ on the boundary $\partial \nu (T)$ during the diffeomorphisms we describe. Replace the dot of the lower circle with $0$ of the middle circle. Then we get the third diagram by an isotopy. This corresponds to a logarithmic transform of $X$ with multiplicity $0$. We next change the dot of the lower circle into $0$ and blow up as shown in the figure. Then we blow down the upper $-1$-framed circle. We again blow up slide and blow down as shown in the Figure 14. By repeating blow-ups and blow-downs similarly, we get the eighth diagram. Finally change the lower $0$-framed unknot into a dotted circle. This gives the last diagram, which is the handlebody of $X_p$. By inspecting the positions of $\alpha, \beta, \gamma$ in the last diagram, we can easily verify that this gluing map corresponds to the map ${\varphi}^{-1}_p: \partial \nu (T) \to T^2\times S^1$.

\begin{remark}
Note that the $0$-log transform operation removes a Stein manifold $T^2\times B^2$ and reglues it via an involution on the boundary. It follows from Gompf~\cite{G} that $0$-log transforms can be plug operations. 
\end{remark}
The rational blow-down has the following relation with the logarithmic transform.

\begin{figure}[htb!]
\begin{center}
\includegraphics[width=1.2in]{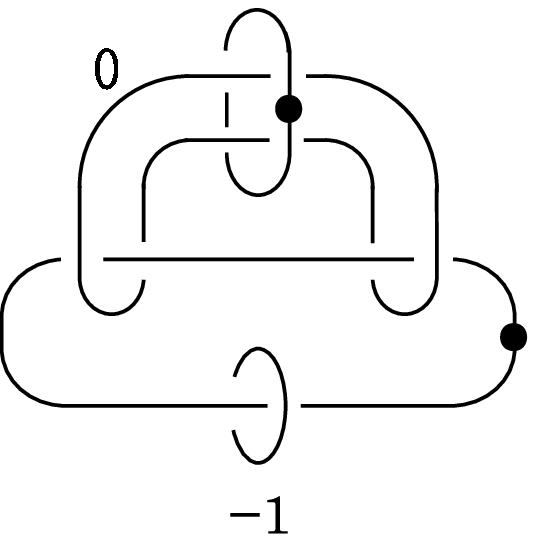}
\end{center}
\caption{fishtail neighborhood}
\label{fig18}
\end{figure}

\begin{figure}[ht!]
\begin{center}
\includegraphics[width=4.2in]{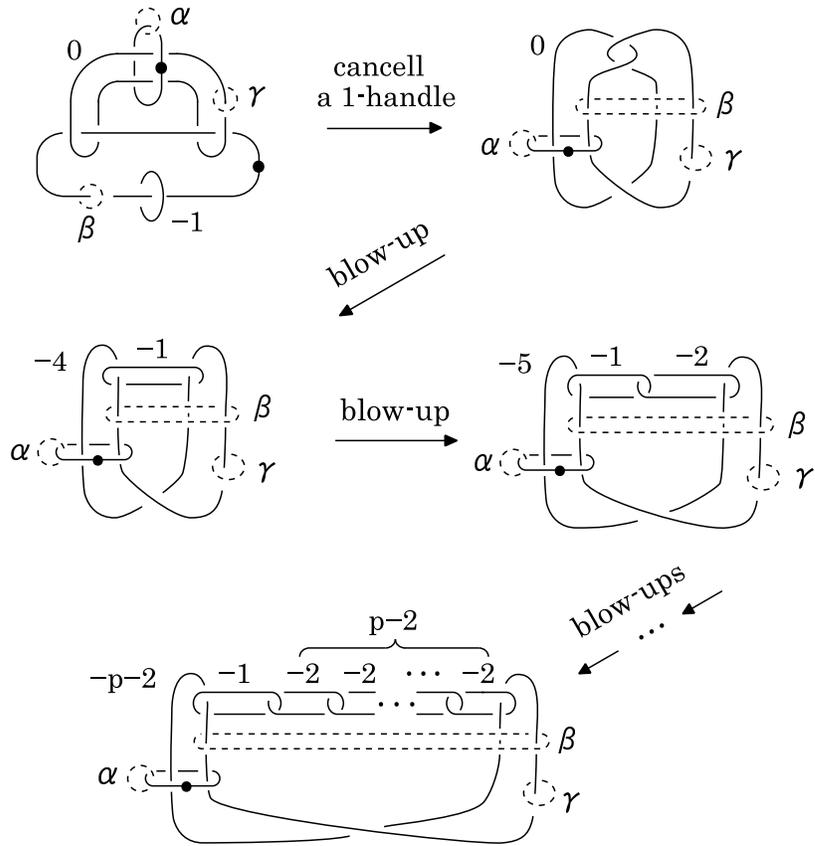}
\caption{Construction of $C_p$ in $X\# (p-1)\overline{\mathbf{CP}}^2$}
\label{fig19}
\end{center}
\end{figure}

\begin{theorem}[{Fintushel-Stern~\cite{FS1}, see also Gompf-Stipsicz~\cite{GS}}]\label{th:log}
Suppose that a smooth $4$-manifold $X$ contains a fishtail neighborhood, that is, the smooth $4$-manifold in Figure~$\ref{fig18}$. 
Let $X_{p}$ $(p\geq 2)$ be the smooth $4$-manifold obtained from $X$ by 
a logarithmic transform with multiplicity $p$ in the fishtail neighborhood. 
Then there exists a copy of $C_p$ in 
$X\# (p-1)\overline{\mathbf{CP}}^2$ such that the rational blow-down of 
$X\# (p-1)\overline{\mathbf{CP}}^2$ along the copy of $C_p$ is diffeomorphic to $X_p$.
\end{theorem}
\begin{proof}
Fintushel-Stern constructed $C_p$ in $X\# (p-1)\overline{\mathbf{CP}}^2$ as indicated in 
 Figure~\ref{fig19}. We keep track of circles $\alpha, \beta, \gamma$ on the boundary $\partial (T^2\times B^2)$. We cancel a 1-handle of the fishtail neighborhood as in the second diagram. By blowing up $(p-1)$ times, we get $C_p$ as in the last diagram. 

We draw a diagram of the rational blow down of $X\# (p-1)\overline{\mathbf{CP}}^2$ along this $C_p$ by the procedure in Figure~\ref{fig16}. Then we can easily show that this is the same operation as logarithmic transform of $X$ with multiplicity $p$, by checking positions of $\alpha, \beta, \gamma$. 
\end{proof}

\clearpage

\section{Corks, Plugs and rational blow-down}\label{section:relation}

In this section we relate the rational blow-down operation with regluing corks and plugs. 
Recall that the rational blow-down operation can change the first and second homology groups. Since regluing corks $W_n$ and plugs $W_{m,n}$ does not change these groups, we need suitable assumptions to relate these operations. 

\begin{figure}[ht!]
\begin{center}
\includegraphics[width=4.9in]{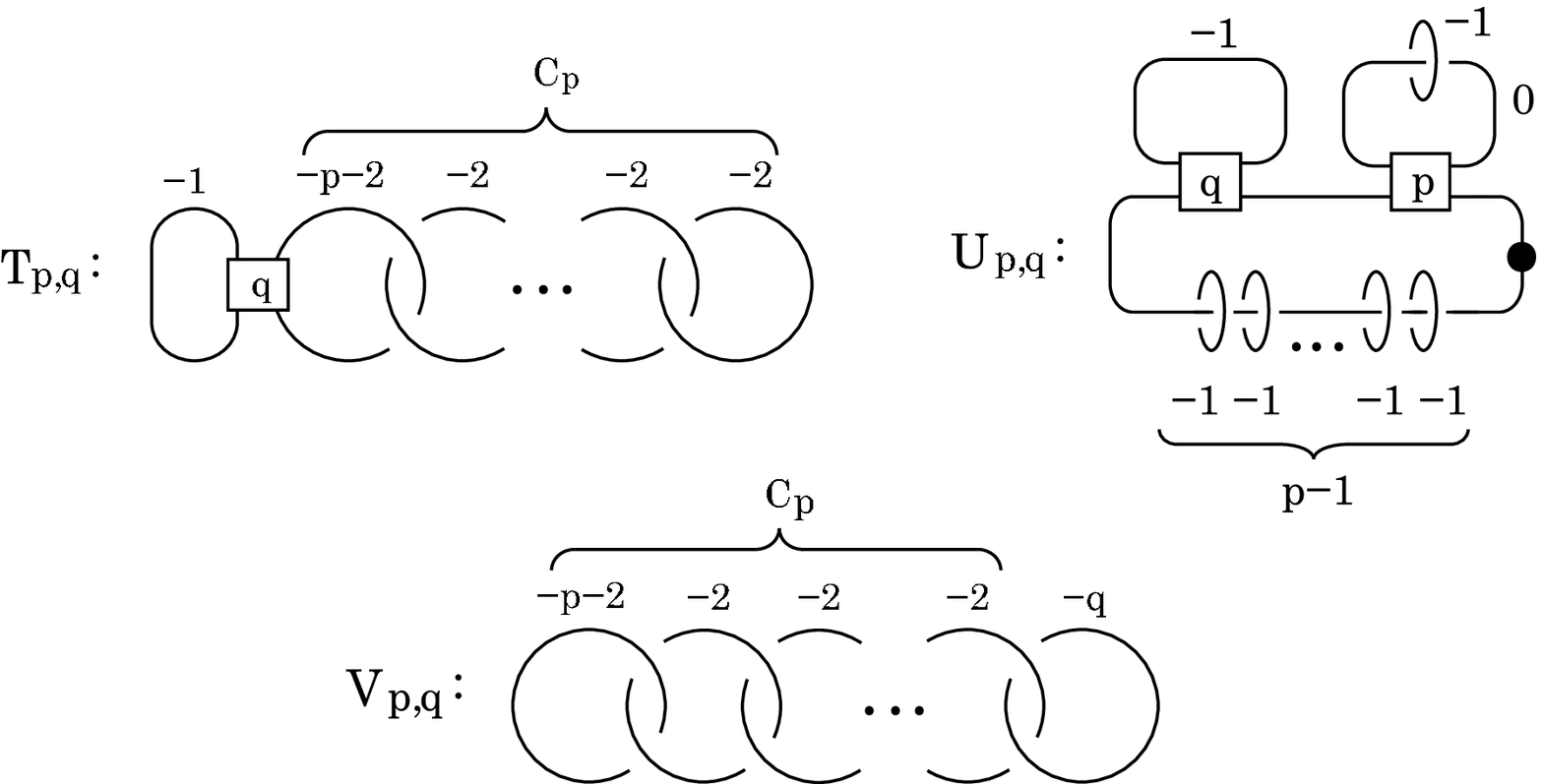}
\caption{}
\label{fig20}
\end{center}
\end{figure}

Let $T_{p,q}$, $U_{p,q}$ and $V_{p,q}$ be the smooth $4$-manifolds in Figure~$\ref{fig20}$. Note that $4$-manifolds $T_{p,q}$, $U_{p,q}$ and $V_{p,q}$ contain $C_p$.
\begin{theorem}\label{prop:cork and rbd}

$(1)$ Suppose that a smooth $4$-manifold $X$ contains $T_{p,p-1}$ $(\text{resp.}$ $T_{p,p+1})$. Let $X_{(p)}$ be the rational blow-down of $X$ along $C_p$. 
Then $X$ contains $W_{p-2}$ $(\text{resp. }$ $W_{p-1})$ such that $X_{(p)}\# (p-1)\overline{\mathbf{CP}}^2$ is obtained from $X$ by removing $W_{p-2}$ $(\text{resp.}$ $W_{p-1})$ and regluing it via $f_{p-2}$ $(\text{resp.}$ $f_{p-1})$. \\

$(2)$ Suppose that a smooth $4$-manifold $X$ contains $U_{p,p-1}$ $(\text{resp.}$ $U_{p,p+1})$. Let $X_{(p)}$ be the rational blow-down of $X$ along $C_p$. 
Then $X$ contains $W_{p-2}$ $(\text{resp.}$ $W_{p-1})$ such that $X_{(p)}\# (p-1)\overline{\mathbf{CP}}^2$ is obtained from $X$ by removing $W_{p-2}$ $(\text{resp.}$ $W_{p-1})$ and regluing it via $f_{p-2}$ $(\text{resp.}$ $f_{p-1})$. \\

$(3)$ Suppose that a smooth $4$-manifold $X$ contains $V_{p,q}$. Let $X_{(p)}$ be the rational blow-down of $X$ along $C_p$. 
Then $X$ contains $W_{q-1,p}$ such that $X_{(p)}\# (p-1)\overline{\mathbf{CP}}^2$ is obtained from $X$ by removing $W_{q-1,p}$ and regluing it via $f_{q-1,p}$. 
\end{theorem}

\begin{proof}
$(1)$ The $q=p+1$ case: We can easily get the handle decomposition of $T_{p,p+1}$ in Figure~$\ref{fig37}$ by handle slides of $C_p$ as in Section $\ref{sec:rbd}$ (See Figure~\ref{fig14}.). Slide handles as shown in Figure~$\ref{fig37}$. In the third diagram of Figure~$\ref{fig37}$, we can find a $0$-framed unknot which links the dotted circle geometrically once. Replace this dot and $0$ as shown in the first row of Figure~$\ref{fig38}$. This operation keeps the diffeomorphism type of $X$ because this corresponds to removing $B^4$ and regluing $B^4$. As a consequence, we can easily find $W_{p-1}$ in $X$. Note that Figure~$\ref{fig39}$ is isotopic to the standard diagram of $W_{p-1}$. By removing $W_{p-1}$ in $X$ and regluing it via $f_{p-1}$, we get the lower diagram of Figure~$\ref{fig38}$. 

We can easily check that $X_p\# (p-1)\overline{\mathbf{CP}}^2$ is obtained by replacing the dot and $0$ as shown in the left side of Figure~$\ref{fig38}$. Hence, we obtain $X_p\# (p-1)\overline{\mathbf{CP}}^2$ from $X$ by removing $W_{p-1}$ and regluing it via $f_{p-1}$. 

\vspace{.05in}

The $q=p-1$ case is similar.
$(2)$ is similar to the $(1)$ case. For $(3)$: handle slides in Figure~$\ref{fig20}$ give Figure~$\ref{fig40}$ (See Figure~$\ref{fig14}$). Then the required claim easily follows from this figure. 
\end{proof}

Applying the theorem above, we easily get the following examples.  Let $E_3$ be the smooth $4$-manifold constructed by the second author~\cite{Y1}. Note that $E_3$ is homeomorphic but not diffeomorphic to $\mathbf{C}\mathbf{P}^2\# 9\overline{\mathbf{CP}}^2$

\begin{example}\label{example:E_3}

$(1)$ $E_3\# 2\overline{\mathbf{CP}}^2$ is obtained from $\mathbf{C}\mathbf{P}^2\# 11\overline{\mathbf{CP}}^2$ by removing a cork $W_1$ and regluing it via $f_1$. Note that $E_3\# 2\overline{\mathbf{CP}}^2$ is homeomorphic but not diffeomorphic to $\mathbf{C}\mathbf{P}^2\# 11\overline{\mathbf{CP}}^2$. \\

$(2)$ $E_3\# 2\overline{\mathbf{CP}}^2$ is obtained from $\mathbf{C}\mathbf{P}^2\# 11\overline{\mathbf{CP}}^2$ by removing a plug $W_{1,3}$ and regluing it via $f_{1,3}$.\\

$(3)$ $E_3\# \mathbf{C}\mathbf{P}^2\#3\overline{\mathbf{CP}}^2$ is diffeomorphic to $2\mathbf{C}\mathbf{P}^2\# 12\overline{\mathbf{CP}}^2$.
\end{example}

\begin{proof}(1) See \cite[Definition~$3.4$ and Figure~$25$]{Y1}. 
The $4$-manifold $E_3$ is obtained by rationally blowing down $C_3$ in a submanifold $T_{3,2}$ of $\mathbf{C}\mathbf{P}^2\# 11\overline{\mathbf{CP}}^2$. Thus Theorem~\ref{prop:cork and rbd} implies the required claim.\\
(2)  See \cite[Proposition~3.2.(1) and Definition~$3.4$]{Y1}. The $4$-manifold $E_3$ is obtained by rationally blowing down $C_3$ in a submanifold $V_{3,2}$ of $\mathbf{C}\mathbf{P}^2\# 11\overline{\mathbf{CP}}^2$. Thus Theorem~\ref{prop:cork and rbd} implies the required claim.\\
(3) This clearly follows from (2) and Proposition~\ref{prop:stabilization}.
\end{proof}
\begin{remark}\label{remark:E_3}
Example~\ref{example:E_3}.(1) and (2) show that regluing a cork $W_3$ of $\mathbf{C}\mathbf{P}^2\# 11\overline{\mathbf{CP}}^2$ is the same operation as regluing a plug $W_{1,3}$, in this case. So, investigating plug structures of $4$-manifolds might be helpful to study cork structures of $4$-manifolds.
\end{remark}
\begin{example}For $n\geq 1$, the elliptic surface $E(2n)\# 2\overline{\mathbf{CP}}^2$ contains $C_{2n}$ such that $(E(2n)\# 2\overline{\mathbf{CP}}^2)_{(2n)}\# (2n-1)\overline{\mathbf{CP}}^2$ is obtained from $E(2n)\# (2n+1)\overline{\mathbf{CP}}^2$ by removing $W_{2n-1}$ and regluing it via $f_{2n-1}$. Here $(E(2n)\# 2\overline{\mathbf{CP}}^2)_{(2n)}$ denotes the rational blow-down of $E(2n)\# 2\overline{\mathbf{CP}}^2$ along the copy of $C_{2n}$. Note that $(E(2n)\# 2\overline{\mathbf{CP}}^2)_{(2n)}\# (2n-1)\overline{\mathbf{CP}}^2$ is homeomorphic but not diffeomorphic to $E(2n)\# (2n+1)\overline{\mathbf{CP}}^2$. 
\end{example}

\begin{proof}Start from Figure~\ref{fig28}. By blowing ups, we can find $U_{2n,2n+1}$ in $E(2n)\# 2\overline{\mathbf{CP}}^2$. Then the required claim follows from Theorem~\ref{prop:cork and rbd}.
\end{proof}

\begin{example}For $n,p\geq 2$, the elliptic surface $E(n)_p\# (p-1)\overline{\mathbf{CP}}^2$ is obtained from $E(n)\# (p-1)\overline{\mathbf{CP}}^2$ by removing a plug $W_{n,p}$ and regluing it via $f_{n,p}$. 
\end{example}

\begin{proof} Start from Figure~\ref{fig8}. Construct $C_p$ in $E(n)\# (p-1)\overline{\mathbf{CP}}^2$ by blowing ups the cusp neighborhood, following the procedure in Figure~\ref{fig19}. Note that a cusp neighborhood naturally contains a fishtail neighborhood. Then we can easily find $V_{p,n+1}$. In this case in the last diagram of Figure~\ref{fig19}, we have $\gamma $ with $-n$ framing.
In the last diagram of Figure~\ref{fig19}, we also have a $-1$-framed unknot.
Slide $\gamma$ over the $-1$-framed knot.
Then we have $V_{p,n+1}$. Now the required claim follows from Theorem~\ref{th:log} and \ref{prop:cork and rbd}.
\end{proof}

\section{Exotic $\mathbf{C}\mathbf{P}^2\# 9\overline{\mathbf{CP}}^2$}\label{exotic}

The second author~(\cite{Y0}, \cite{Y1}) constructed the minimal smooth $4$-manifold $E_3'$ which is homeomorphic but not diffeomorphic to $\mathbf{C}\mathbf{P}^2\# 9\overline{\mathbf{CP}}^2$, and it has neither $1$- nor $3$-handles. He also gave the way to draw the handlebody diagram of $E_3'$, but he did not draw it. In this section, we draw a whole diagram of $E_3'$, by slightly changing the construction in \cite{Y1} to obtain a simple diagram. From this we find a cork and a plug of $E_3'$. As far as the authors know, this is the first example of a cork and a plug in minimal $4$-manifolds. 

\vspace{.05in}

\begin{theorem}\label{th:figure}
$(1)$ Figure~$\ref{fig21}$ is a handlebody diagram of a smooth $4$-manifold which is homeomorphic but not diffeomorphic to $\mathbf{CP}^2\# 9\overline{\mathbf{CP}}^2$ and has the same Seiberg-Witten invariant as the elliptic surface $E(1)_{2,3}$. Consequently, this $4$-manifold is minimal, that is, it does not contain any $2$-sphere with self intersection number $-1$. It is clear from Figure~$\ref{fig21}$ that this minimal $4$-manifold has neither $1$- nor $3$-handles in that handle decomposition. \\

\noindent $(2)$ $W_1$ is a cork of the minimal $4$-manifold above. \\

\noindent $(3)$ $W_{1,3}$ is a plug of the minimal $4$-manifold above. 

\begin{figure}[ht]
\begin{center}
\includegraphics[width=4.8in]{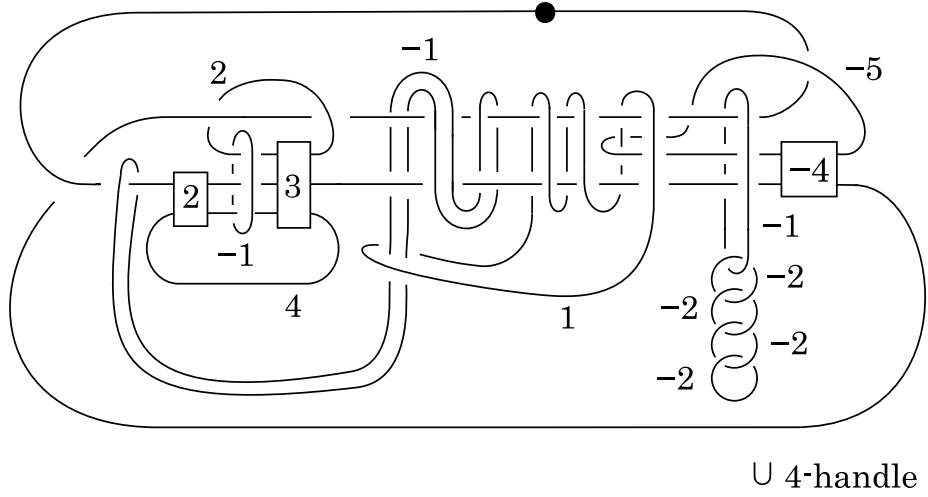}
\caption{A minimal exotic $\mathbf{CP}^2\# 9\overline{\mathbf{CP}}^2$}
\label{fig21}
\end{center}
\end{figure}
\end{theorem}

\vspace{.05in}

\begin{remark}
(1) In ~\cite{A5} the first author proved that $\overline{W}_1$ is a cork of the elliptic surface $E(1)_{2,3}$. \\
(2) In \cite{Y2}, the second author constructed the smooth $4$-manifold $X_{a,3}$ $(3\leq a\leq 7)$ which is homeomorphic but not diffeomorphic to $\mathbf{CP}^2\# (12-a)\overline{\mathbf{CP}}^2$. As in Theorem~\ref{th:figure}, we can  prove that $W_2$ is a cork of $X_{4,3}$. 
\end{remark}

\vspace{.1in}

Here we also show that corks and plugs can be knotted, by similar constructions of exotic rational surfaces in \cite{Y0} and \cite{Y1}. 

\begin{theorem}[Knotting corks]\label{th:knotted cork}
There exist two copies of $W_3$ in $\mathbf{CP}^2\# 14\overline{\mathbf{CP}}^2$ with the following two properties:\\
$(1)$ Two $4$-manifolds obtained from $\mathbf{CP}^2\# 14\overline{\mathbf{CP}}^2$ by removing each copy of $W_3$ and regluing it via $f_3$ are homeomorphic but not diffeomorphic to $\mathbf{CP}^2\# 14\overline{\mathbf{CP}}^2$. In particular, each copy of $W_3$ is a cork of $\mathbf{CP}^2\# 14\overline{\mathbf{CP}}^2$.\\
$(2)$ Two $4$-manifolds obtained from $\mathbf{CP}^2\# 14\overline{\mathbf{CP}}^2$ by removing each copy of $W_3$ and regluing it via $f_3$ are homeomorphic but not diffeomorphic to each other. In particular, the above two corks of $\mathbf{CP}^2\# 14\overline{\mathbf{CP}}^2$ are not isotopic to each other.  
\end{theorem}

\vspace{.02in}

\begin{theorem}[Knotting plugs]\label{th:knotted plug}
There exist two copies of $W_{4,5}$ in $\mathbf{CP}^2\# 16\overline{\mathbf{CP}}^2$ with the following two properties:\\
$(1)$ Two $4$-manifolds obtained from $\mathbf{CP}^2\# 16\overline{\mathbf{CP}}^2$ by removing each copy of $W_{4,5}$ and regluing it via $f_{4,5}$ are homeomorphic but not diffeomorphic to $\mathbf{CP}^2\# 16\overline{\mathbf{CP}}^2$. In particular, each copy of $W_{4,5}$ is a plug of $\mathbf{CP}^2\# 16\overline{\mathbf{CP}}^2$.\\
$(2)$ Two $4$-manifolds obtained from $\mathbf{CP}^2\# 16\overline{\mathbf{CP}}^2$ by removing each copy of $W_{4,5}$ and regluing it via $f_{4,5}$ are homeomorphic but not diffeomorphic to each other. In particular, the above two plugs of $\mathbf{CP}^2\# 16\overline{\mathbf{CP}}^2$ are not isotopic to each other.  
\end{theorem}

Let $h,e_1,e_2,\dots,e_n$ be a canonical orthogonal basis of $H_2(\mathbf{CP}^2\# n\overline{\mathbf{CP}}^2;\mathbf{Z})$ such that $h^2=1$ and $e_1^2=e_2^2=\dots=e_n^2=-1$. 

In handlebody diagrams, we often write the second homology classes given by $2$-handles, instead of usual framings. Note that the square of the homology class given by a 2-handle is equal to the usual framing. 

\begin{proof}[Proof of Theorem~$\ref{th:figure}$]
$(1)$ Recall the construction of $E_3'$ in \cite{Y1}: Figure~\ref{fig41} is a standard diagram of $\mathbf{C}\mathbf{P}^2\# 2\overline{\mathbf{CP}}^2$. By sliding handles as  in 
\cite[Figure~$8\sim17$]{Y1}, we obtain Figure~\ref{fig42}. Blowing ups gives Figure~\ref{fig43}. We have Figure~\ref{fig44} by an isotopy. We get Figure~\ref{fig46} by handle slides as indicated in Figure~\ref{fig44} and \ref{fig45}. By blowing up, we obtain Figure~\ref{fig47}. Note that a copy of $C_5$ (see Figure~\ref{fig48}) is inside Figure~\ref{fig47}. By rationally blowing down this copy of $C_5$, we get a smooth $4$-manifold $E'_3$. By using the rational blow-down procedure in Figure~\ref{fig15}, we obtain Figure~\ref{fig21} of $E'_3$.
Though this construction of $E'_3$ is slightly different from \cite{Y1}, the same argument as in \cite{Y1} shows that this $E'_3$ is homeomorphic to $\mathbf{CP}^2\# 9\overline{\mathbf{CP}}^2$ and has the same Seiberg-Witten invariant as the elliptic surface $E(1)_{2,3}$. Therefore this $E'_3$ is minimal. \\

$(2)$ We can get Figure~\ref{fig49} from Figure~\ref{fig41} by introducing a $2$-handle/$3$-handle pair and sliding the new $2$-handle. We obtain Figure~\ref{fig50} of $\mathbf{CP}^2\# 13\overline{\mathbf{CP}}^2$ from Figure~\ref{fig49}, similar to $(1)$. Note that there is a copy of $C_5$ is inside of Figure~\ref{fig50}. This copy of $C_5$ is isotopic to the one in Figure~\ref{fig47}, because we constructed the handle decomposition in Figure~\ref{fig50} without sliding any handles over the new $2$-handle which comes from a $2$-handle/$3$-handle pair. We thus get Figure~\ref{fig51} of $E'_3$ by rationally blowing down. The first diagram in Figure~\ref{fig52} is obtained from Figure~\ref{fig51} by ignoring handles. Handle slides give the last diagram in Figure~\ref{fig52}. Note that there is $W_1$ inside the last diagram in Figure~\ref{fig52}. By removing this $W_1$ from $E'_3$ and regluing it via $f_1$, we get a $2$-sphere with self intersection number $-1$. Therefore $W_1$ is a cork of $E'_3$. \\

(3) Similar to (2), we obtain the first diagram in Figure~\ref{fig53} from Figure~\ref{fig51} of $E'_3$ by ignoring handles. Handle slides give the last diagram in Figure~\ref{fig53}. Notice the $W_{1,3}$  inside of the last diagram in Figure~\ref{fig53}.  By removing $W_{1,3}$ from $E'_3$ and regluing it via $f_{1,3}$, we get a $2$-sphere with self intersection number $-1$. Therefore $W_{1,3}$ is a plug of $E'_3$. 
\end{proof}

\begin{remark}
$(1)$ In the proof above, we gave one construction of $E'_3$. We can give a lot of different constructions by using knotted bands in handle slides. We do not know if choices of bands in handle slides affect diffeomorphism types of $E'_3$, though they do not affect homeomorphism types and Seiberg-Witten invariants. See also \cite[Remark~6.1]{Y1} and \cite{Y3}. \\
$(2)$ In the proof above, we drew a digram of $E'_3$ by the rational blow-down procedure in Figure~~\ref{fig15}. Another rational blow-down procedure in Figure~\ref{fig16} gives a different diagram of $E'_3$, though they are diffeomorphic. \end{remark}

\subsection{Knotted corks and plugs in $4$-manifolds}

\begin{proposition}\label{prop:construction_knotted_cork}

$(1)$ $\mathbf{CP}^2\# 13\overline{\mathbf{CP}}^2$ admits a handle decomposition as in Figure~$\ref{fig22}$. Note that Figure~$\ref{fig22}$ contains $C_5$. \\

$(2)$ $\mathbf{CP}^2\# 13\overline{\mathbf{CP}}^2$ admits a handle decomposition as in Figure~$\ref{fig23}$. Note that Figure~$\ref{fig23}$ contains $C_5$.

\begin{figure}[h!]
\begin{center}
\includegraphics[width=3.2in]{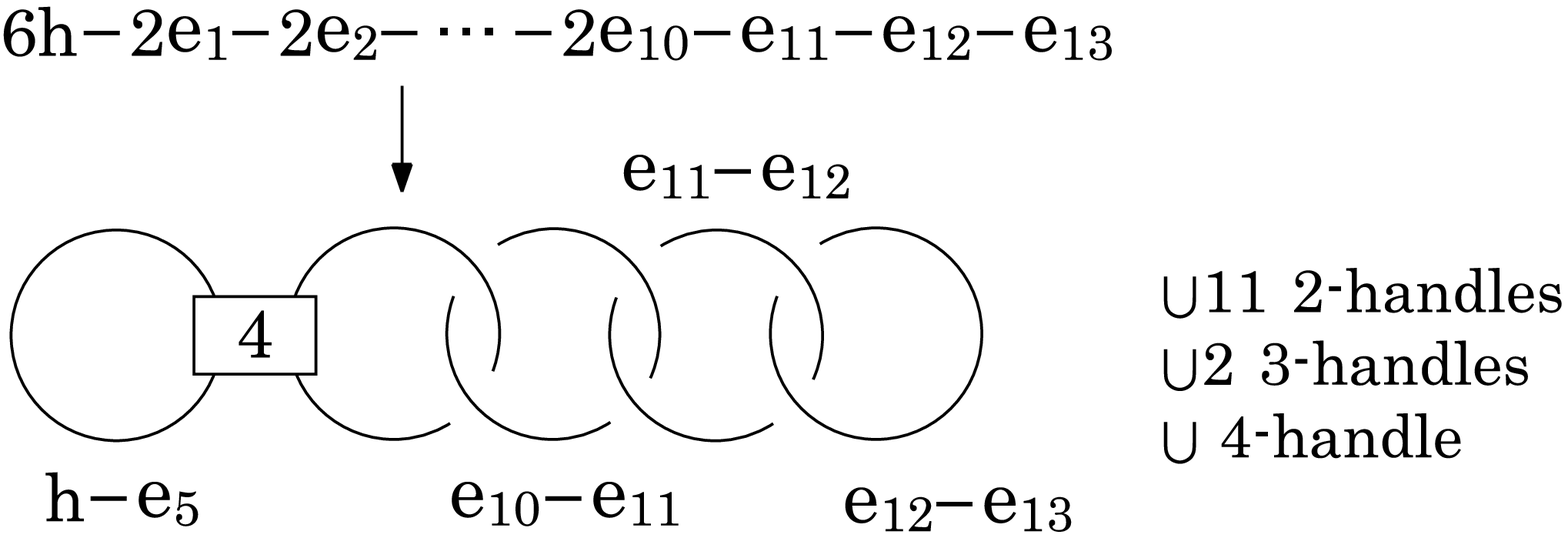}
\caption{$\mathbf{CP}^2\# 13\overline{\mathbf{CP}}^2$}
\label{fig22}
\end{center}
\end{figure}
\begin{figure}[h!]
\begin{center}
\includegraphics[width=3.2in]{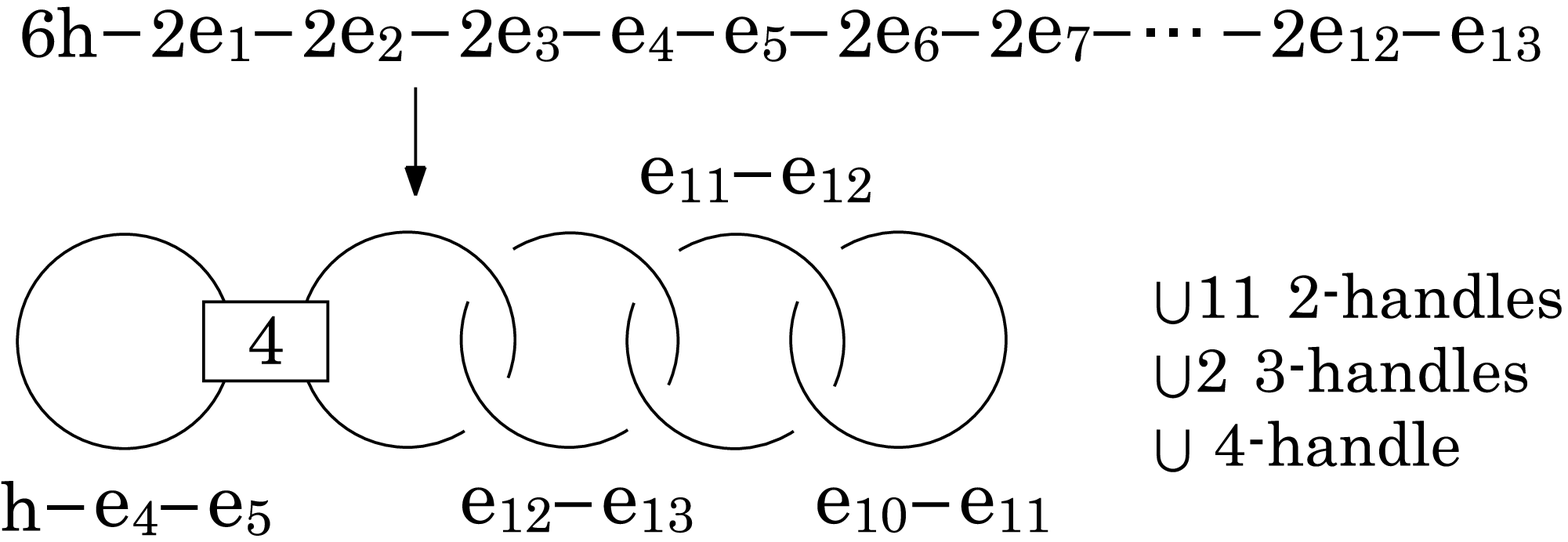}
\caption{$\mathbf{CP}^2\# 13\overline{\mathbf{CP}}^2$}
\label{fig23}
\end{center}
\end{figure}
\end{proposition}
\begin{proof}(1) Handle slides, introducing $2$-handle/$3$-handle pairs, and blowing ups as in \cite[Figure~$8\sim18$]{Y1}, i.e.  the moves similar to Figure~\ref{fig41} $\leadsto $ Figure~\ref{fig44},  give Figure \ref{fig54} of $\mathbf{CP}^2\# 3\overline{\mathbf{CP}}^2$. By an isotopy, we get Figure~\ref{fig55}. By the handle slide indicated in Figure~\ref{fig55}, we get Figure~\ref{fig56}. Blowing ups give Figure~\ref{fig57}. We blow up as in Figure~\ref{fig19}. Then we obtain Figure~\ref{fig22}. \\
(2) Blowing ups in Figure~\ref{fig55} of $\mathbf{CP}^2\# 3\overline{\mathbf{CP}}^2$ give Figure~\ref{fig58} of $\mathbf{CP}^2\# 5\overline{\mathbf{CP}}^2$. By further blow-ups we get Figure~\ref{fig59}. We then obtain Figure~\ref{fig23} of $\mathbf{CP}^2\# 13\overline{\mathbf{CP}}^2$, by a similar construction to \cite[Figure~$19\sim25$]{Y1}. 
\end{proof}

\vspace{.02in}

\begin{definition}
Let $X_5$ (resp.\ $X_3$) be the smooth $4$-manifold obtained from $\mathbf{CP}^2\# 13\overline{\mathbf{CP}}^2$ by rationally blowing down the copy of $C_5$ in Proposotion~$\ref{prop:construction_knotted_cork}$.(1) (resp.\ in Proposotion~$\ref{prop:construction_knotted_cork}$.(2)). 
\end{definition}

\vspace{.02in}

\begin{corollary}\label{cor:X and cork}
$(1)$ $X_5\# \overline{\mathbf{CP}}^2$ is obtained from $\mathbf{CP}^2\# 14\overline{\mathbf{CP}}^2$ by rationally blowing down the copy of $C_5$ in Figure~$\ref{fig24}$. \\

$(2)$ $X_3\# \overline{\mathbf{CP}}^2$ is obtained from $\mathbf{CP}^2\# 14\overline{\mathbf{CP}}^2$ by rationally blowing down the copy of $C_5$ in Figure~$\ref{fig25}$.
\begin{figure}[h!]
\begin{center}
\includegraphics[width=3.2in]{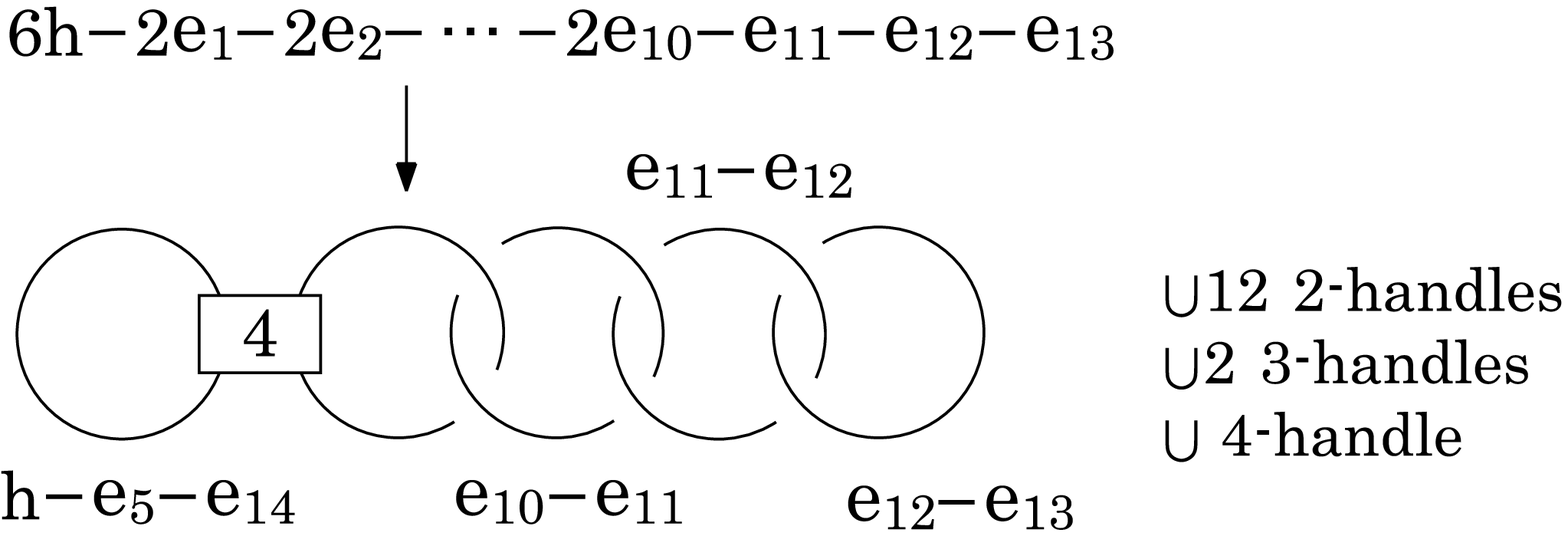}
\caption{$\mathbf{CP}^2\# 14\overline{\mathbf{CP}}^2$}
\label{fig24}
\end{center}
\end{figure}
\begin{figure}[h!]
\begin{center}
\includegraphics[width=3.2in]{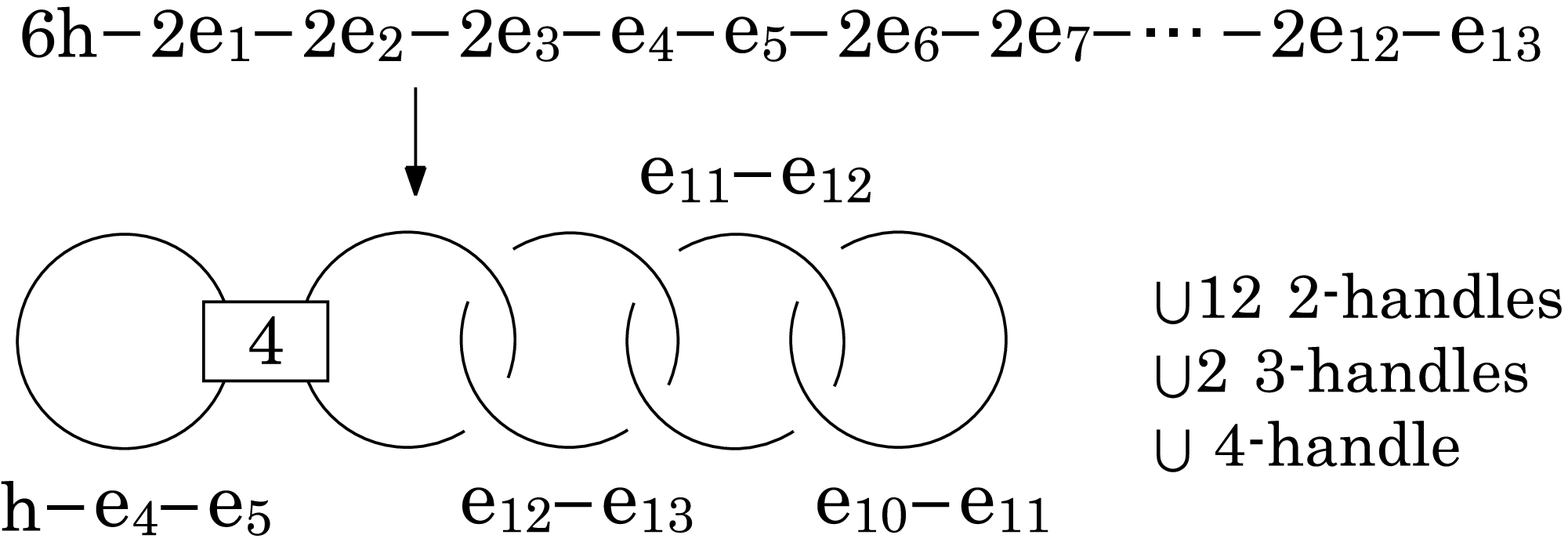}
\caption{$\mathbf{CP}^2\# 14\overline{\mathbf{CP}}^2$}
\label{fig25}
\end{center}
\end{figure}
\end{corollary}
\begin{corollary}\label{cor:X is simp. conn}
$(1)$ $X_5\# 5\overline{\mathbf{CP}}^2$ is obtained from $\mathbf{CP}^2\# 14\overline{\mathbf{CP}}^2$ by removing $W_3$ and \linebreak
 regluing it via $f_3$. \\

$(2)$ $X_3\# 5\overline{\mathbf{CP}}^2$ is obtained from $\mathbf{CP}^2\# 14\overline{\mathbf{CP}}^2$ by removing $W_3$ and regluing it via $f_3$. 
\end{corollary}
\begin{proof}
This ovbiously follows from Corollary~\ref{cor:X and cork} and Theorem~\ref{prop:cork and rbd}. 
\end{proof}
\begin{corollary}\label{cor:X_3 and X_5 are not diffeomorphic}
$(1)$ $X_3$ is homeomorphic to the elliptic surface $E(1)_{2,3}$ and has the same Seiberg-Witten invariant as $E(1)_{2,3}$. In particular, $X_3$ is homeomorphic but not diffeomorphic to $\mathbf{CP}^2\# 9\overline{\mathbf{CP}}^2$.\\

$(2)$ $X_5$ is homeomorphic to the elliptic surface $E(1)_{2,5}$ and has the same Seiberg-Witten invariant as $E(1)_{2,5}$. In particular, $X_5$ is homeomorphic but not diffeomorphic to $\mathbf{CP}^2\# 9\overline{\mathbf{CP}}^2$.\\

$(3)$ $X_3\#k\overline{\mathbf{CP}}^2$, $X_5\#k\overline{\mathbf{CP}}^2$ and $\mathbf{CP}^2\# (9+k)\overline{\mathbf{CP}}^2$ $(k\geq 0)$ are mutually homeomorphic but not diffeomorphic to each other.
\end{corollary}
\begin{proof}
Corollary~\ref{cor:X is simp. conn} implies that $X_3$ and $X_5$ are simply connected. Thus Freedman's \linebreak
theorem and Rochlin's theorem imply that $X_3$ and $X_5$ are homeomorphic to $\mathbf{CP}^2\# 9\overline{\mathbf{CP}}^2$. The same argument as in \cite{Y1} (see also \cite{Y3}) shows $X_3$ (resp.\ $X_5$) has the same Seiberg-Witten invariant as $E(1)_{2,3}$ (resp.\ $E(1)_{2,5}$). ($X_3$ and $X_5$ corresponds to $E_3'$ and $E_5$ in~\cite{Y1}, respectively.) We thus can prove that the Seiberg-Witten invariants of $X_3\#k\overline{\mathbf{CP}}^2$, $X_5\#k\overline{\mathbf{CP}}^2$ and $\mathbf{CP}^2\# (9+k)\overline{\mathbf{CP}}^2$ $(k\geq 0)$ are mutually different, similar to Szab\'{o} \linebreak
\cite[Lemmas~3.6 and 3.5]{Sz}. 
\end{proof}
\begin{proof}[Proof of Theorem~$\ref{th:knotted cork}$]
Theorem~\ref{th:knotted cork} now clearly follows from Corollary~\ref{cor:X is simp. conn} and \ref{cor:X_3 and X_5 are not diffeomorphic}. 
\end{proof}
\begin{proof}[Proof of Theorem~$\ref{th:knotted plug}$] It easily follows from \cite[Proposition~3.3]{Y1} that $E'_3\# 3\overline{\mathbf{CP}}^2$ is \linebreak
obtained from $\mathbf{CP}^2\# 16\overline{\mathbf{CP}}^2$ by rationally blowing down $C_5$ inside $V_{5,5}$. (In this proof, $E_3'$ denotes the manifold defined in~\cite{Y1}, not the manifold defined in Theorem~\ref{th:figure} in this paper.) Theorem~\ref{prop:cork and rbd} thus shows that $E'_3\# 7\overline{\mathbf{CP}}^2$ in \cite{Y1} is obtained from $\mathbf{CP}^2\# 16\overline{\mathbf{CP}}^2$ by removing $W_{4,5}$ and regluing it via $f_{4,5}$. Similarly, \cite[Proposition~3.2.(2)]{Y1} implies that $E_5\# 7\overline{\mathbf{CP}}^2$ is obtained from $\mathbf{CP}^2\# 16\overline{\mathbf{CP}}^2$ by removing $W_{4,5}$ and regluing it via $f_{4,5}$. Since $E_5$ has the same Seiberg-Witten invariant as the elliptic surface $E(1)_{2,5}$, we can prove that $E'_3\# 7\overline{\mathbf{CP}}^2$, $E_5\# 7\overline{\mathbf{CP}}^2$ and $\mathbf{CP}^2\# 16\overline{\mathbf{CP}}^2$ are mutually homeomorphic but not diffeomorphic, similar to Corollary~\ref{cor:X_3 and X_5 are not diffeomorphic}. 
\end{proof}
We finally point out an interesting property of $E_3$ and $E'_3$. 
\begin{proposition}
The minimal $4$-manifold in Figure~$\ref{fig21}$ is obtained from $\mathbf{CP}^2\# 14\overline{\mathbf{CP}}^2$ by rationally blowing down $C_6$. 
\end{proposition}
\begin{proof}
By blowing up in Figure~\ref{fig47}, we get Figure~\ref{fig60}. We can easily find a copy of $C_6$ in Figure~\ref{fig60}. Draw a diagram of the rational blow-down of $\mathbf{CP}^2\# 14\overline{\mathbf{CP}}^2$ along this $C_6$. Then we get Figure~\ref{fig61}. By handle slides, we can change the diagram into Figure~\ref{fig21}. 
\end{proof}

\begin{remark}
We can change the construction of $E'_3$ so that, for any $5\leq p\leq 11$, $E'_3$ can be obtained from $\mathbf{CP}^2\# (8+p)\overline{\mathbf{CP}}^2$ by rationally blowing down $C_p$. The smooth $4$-manifold $E_3$ can also be obtained from $\mathbf{CP}^2\# (8+p)\overline{\mathbf{CP}}^2$ $(3\leq p\leq 11)$ by rationally blowing down $C_p$. Note that $E_3$ and $E'_3$ are homeomorphic to $E(1)_{2,3}$ and have the same Seiberg-Witten invariants as $E(1)_{2,3}$. However, it is not known whether or not $E(1)_{2,3}$ itself can be obtained from $\mathbf{CP}^2\# (8+p)\overline{\mathbf{CP}}^2$ by rationally blowing down $C_p$, for $p\geq 4$. Moreover, there is no known other minimal $4$-manifolds for which more than three different rational blow-down constructions exist. 

\vspace{.05in}

For various handle diagrams of $E(1)_{2,3}$, see Harer-Kas-Kirby~\cite{HKK}, Gompf~\cite{G}, the first author~\cite{A5}, and the second author~\cite{Y4}.  The handlebody of $E(1)_{2,3}$ in ~\cite{A5} has no $1$- and $3$- handles.  For more diagrams of exotic rational surfaces, see ~\cite{A5}, ~\cite{Y2}, \cite{Y4}.
\end{remark}
\section{Further remarks}We conclude this paper by making some remarks.
\subsection{General remarks}
\begin{remark}(1) Theorem~\ref{th:cork} and~\ref{ex:corks and plugs} show that $W_n$ and $\overline{W}_n$ are corks and that $W_{m,n}$ are plugs. By a similar technique, we can easily give more examples of corks (which consist of one 1-handle and one 2-handle) and plugs. \\
(2) In Figure~\ref{fig1}, we defined the smooth $4$-manifold $\overline{W}_n$. Note that this symbol does not denote the reverse orientaion of $W_n$. 
\end{remark}

\begin{remark}Theorem~\ref{prop:cork and rbd}.\ (1) and (2) describe relations between rational blow-down operations and cork operations, under some conditions. We can similarly prove the following theorem. Notice the slight difference from Theorem~\ref{prop:cork and rbd}. \end{remark}
\begin{theorem}
$(1)$ Suppose that a smooth $4$-manifold $X$ contains $T_{p,p-1}$ $(\text{resp.}$ $T_{p,p+1})$. Let $X_{(p)}$ be the rational blow-down of $X$ along $C_p$. 
Then $X$ contains $W_{p-1}$ $(\text{resp. }$ $W_{p})$ such that $X_{(p)}\# (p-1)\overline{\mathbf{CP}}^2$ is obtained from $X$ by removing $W_{p-1}$ $(\text{resp.}$ $W_{p})$ and regluing it via $f_{p-1}$ $(\text{resp.}$ $f_{p})$. \medskip \\
$(2)$ Suppose that a smooth $4$-manifold $X$ contains $U_{p,p-1}$ $(\text{resp.}$ $U_{p,p+1})$. Let $X_{(p)}$ be the rational blow-down of $X$ along $C_p$. 
Then $X$ contains $W_{p-1}$ $(\text{resp.}$ $W_{p})$ such that $X_{(p)}\# (p-1)\overline{\mathbf{CP}}^2$ is obtained from $X$ by removing $W_{p-1}$ $(\text{resp.}$ $W_{p})$ and regluing it via $f_{p-1}$ $(\text{resp.}$ $f_{p})$. 
\end{theorem}

\begin{remark}
Let $Z$ be a simply connected closed smooth $4$-manifold, and $\{ Z_n \}_{n\geq 1}$  mutually different smooth structures on $Z$. Then, by Theorem~\ref{th:1.1}, there exist corks $(C_i,\tau_i)$ of $Z$ such that $Z_i$ is obtained from $Z$ by removing $C_i$ and regluing it via $\tau_i$. It seems natural to ask whether or not $(C_i,\tau_i)$ and $(C_j,\tau_j)$ $(i\neq j)$ can be the same. Theorem~\ref{th:knotted cork} shows that $(C_1,\tau_1)$ and $(C_2,\tau_2)$ can be the same, namely, the embeddings of corks are different. In a forthcoming paper, we further discuss knotting corks and plugs. 
\end{remark}
\subsection{$1$- and $3$-handles of $4$-manifolds} 
%
It is not known if the minimal numbers of $1$- and $3$-handles of smooth $4$-manifolds (which are invariants of smooth structures) are non-trivial for simply connected closed $4$-manifolds. However, the following example shows that this invariant is non-trivial for simply connected $4$-manifolds with boundary. 

\begin{example}
There exist two homeomorphic simply connected smooth $4$-manifolds with boundary with the following property: one admits a handle decomposition without \linebreak
 $1$- and $3$-handles; the other has at least a $1$- or $3$-handle in each of its handle decompositions. In particular, these smooth structures can be detected by minimal numbers of its $1$- and $3$-handles. 
\end{example}
\begin{proof}
The first author~\cite{A4} constructed smooth 4-manifolds which are homeomorphic but not diffeomorphic to the cusp neighborhood, by knot surgery in the cusp neighborhood. 
Ozsv\'{a}th-Szab\'{o}~\cite{OS} proved that if Dehn surgery along a knot $K$ in $S^3$ is orientation-preserving diffeomorphic to $0$-framed right trefoil surgery, then $K$ is isotopic to right trefoil. 
These two theorems imply the claim.
\end{proof}

As far as the authors know, no other such examples are known. 
\begin{question}
Does $W^2_n$ $(n: \text{large})$ in Figure~$\ref{fig4}$ admit a handle decomposition without $1$- and $3$-handles? Note that $W^1_n$ $(n\geq 1)$ is homeomorphic but not diffeomorphic to $W^2_n$ and has a handle decomposition without 1- and 3-handles. 
\end{question}
\begin{remark}
$W^2_1$ admits a handle decomposition without $1$- and $3$-handles. 
\end{remark}

\begin{figure}[p]
\begin{center}
\includegraphics[width=4.9in]{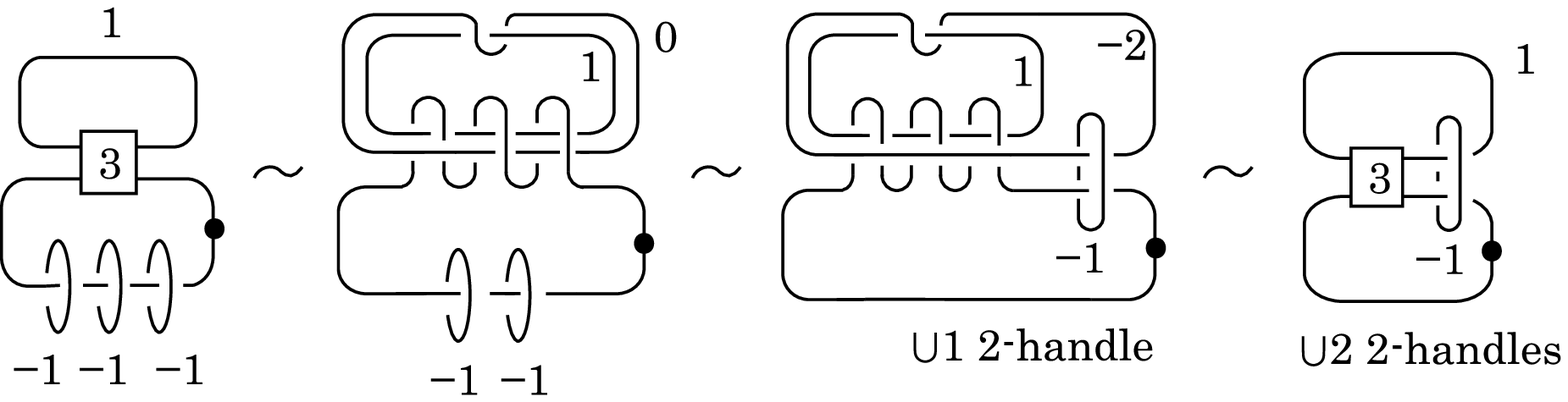}
\caption{Handle slides}
\label{fig26}
\end{center}
\end{figure}
\begin{figure}[ht!]
\begin{center}
\includegraphics[width=4.0in]{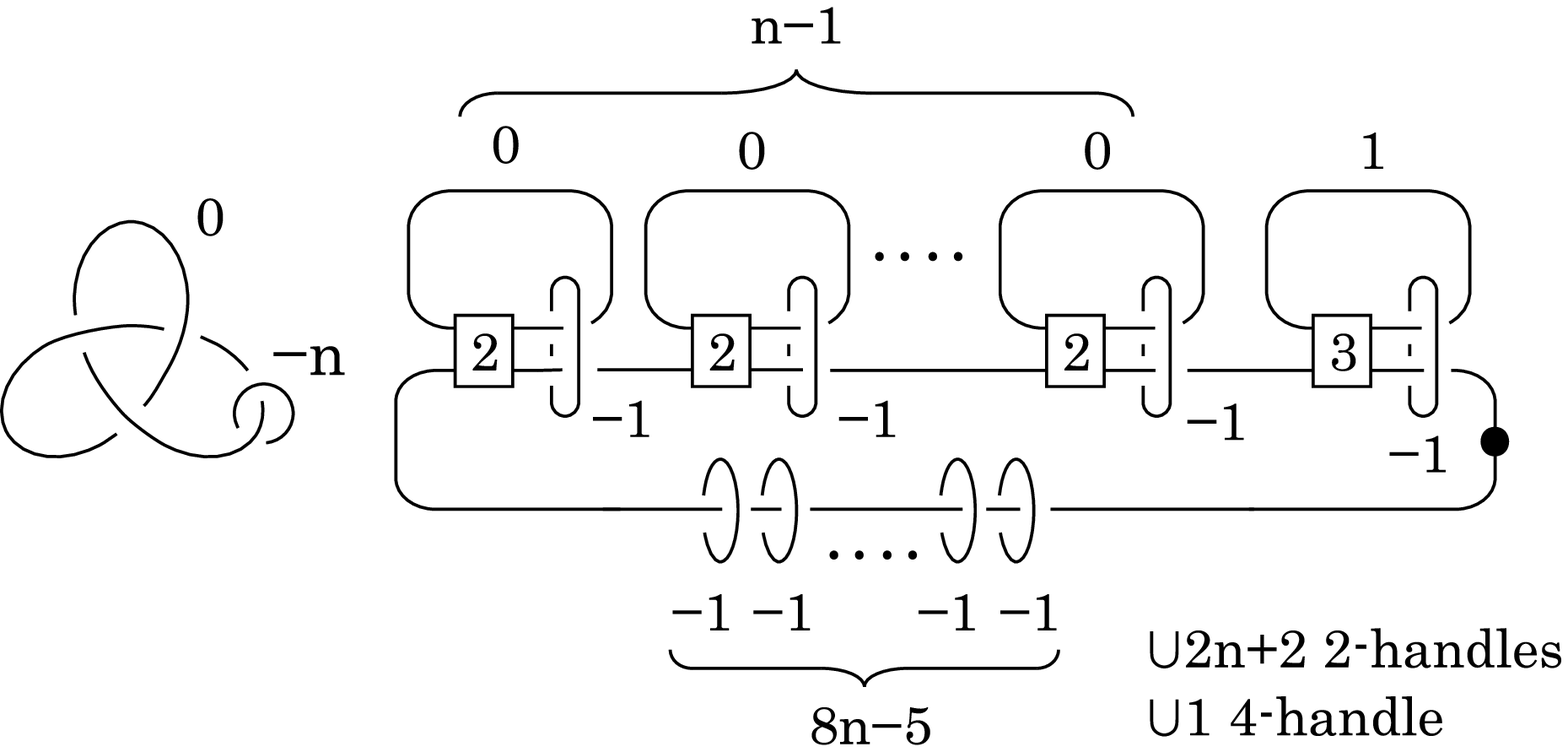}
\caption{$E(n)$}
\label{fig27}
\end{center}
\end{figure}
\begin{figure}[ht!]
\begin{center}
\includegraphics[width=2.2in]{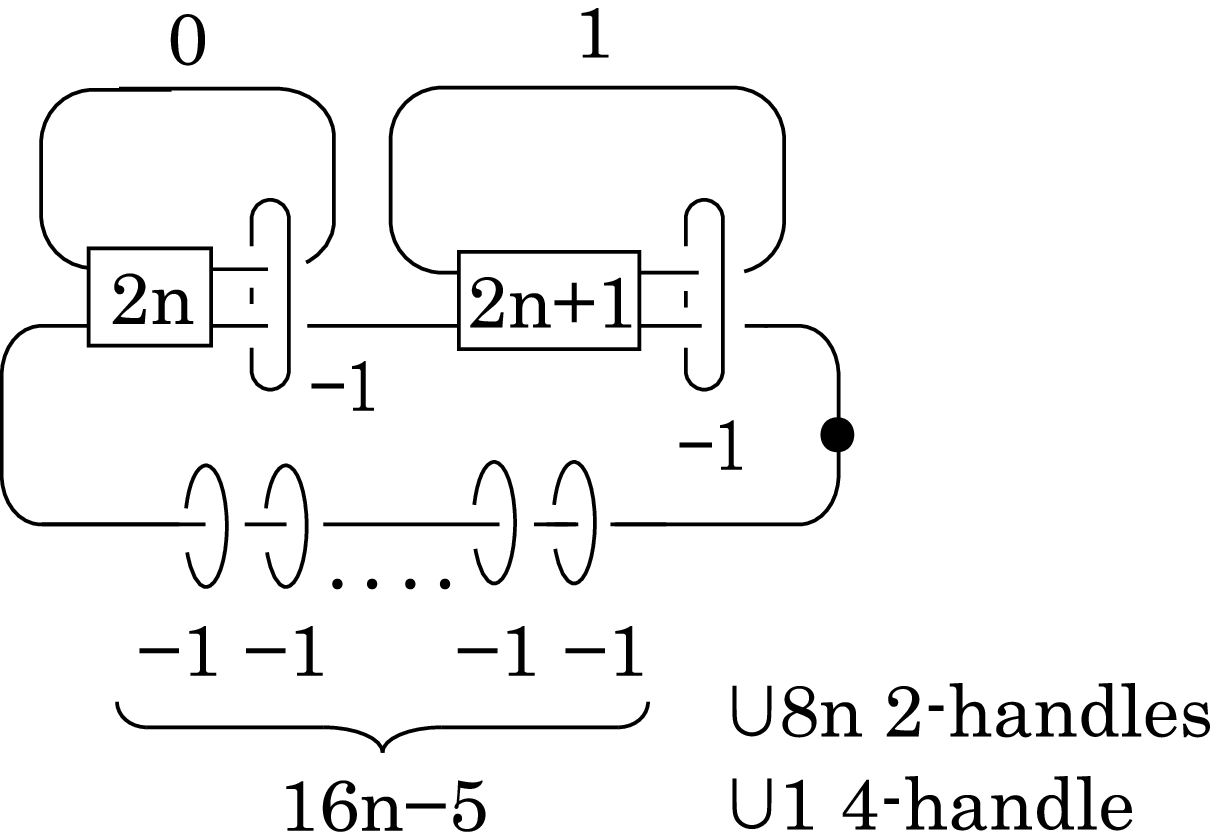}
\caption{$E(2n)$}
\label{fig28}
\end{center}
\end{figure}
\begin{figure}[ht!]
\begin{center}
\includegraphics[width=2.6in]{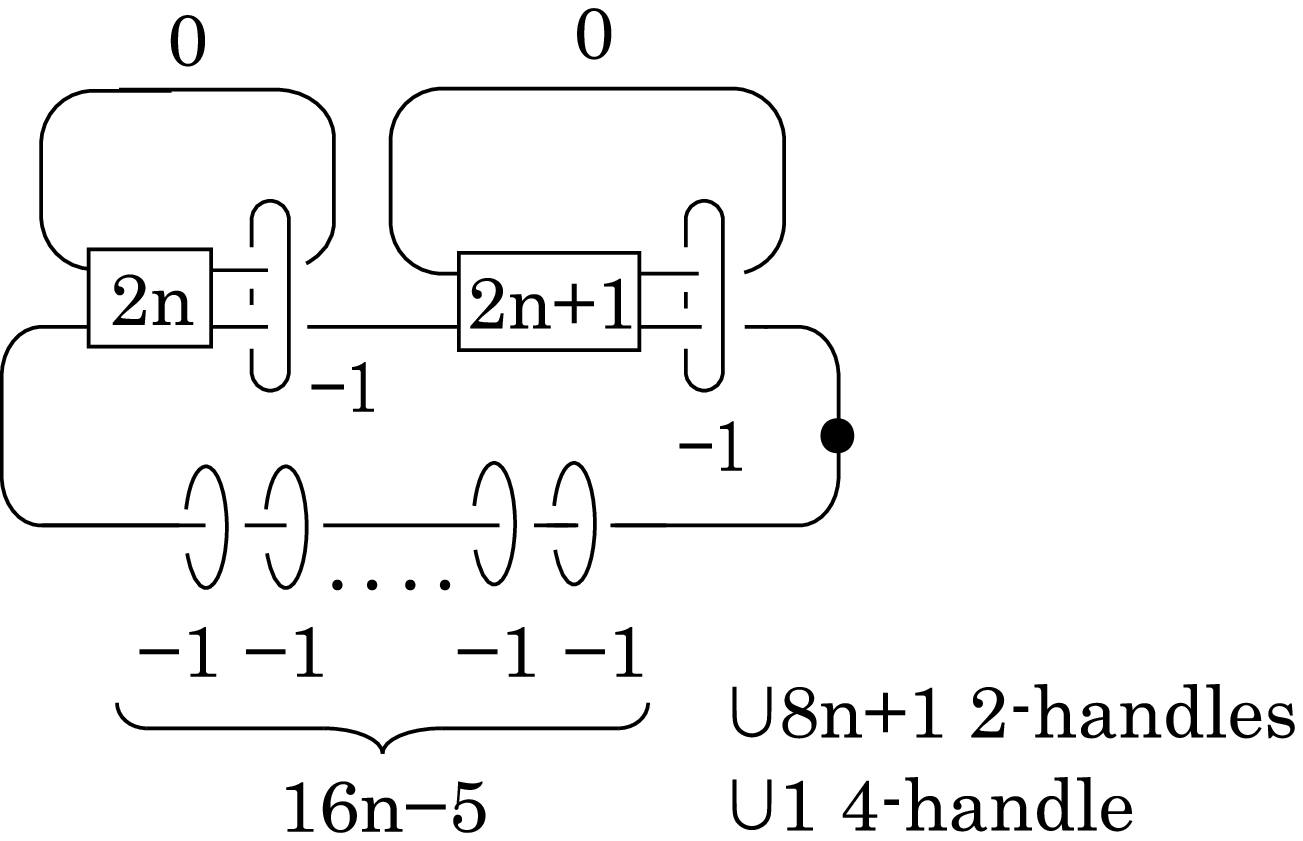}
\caption{$E(2n)\# \overline{\mathbf{CP}}^2$}
\label{fig29}
\end{center}
\end{figure}
\begin{figure}[ht!]
\begin{center}
\includegraphics[width=2.4in]{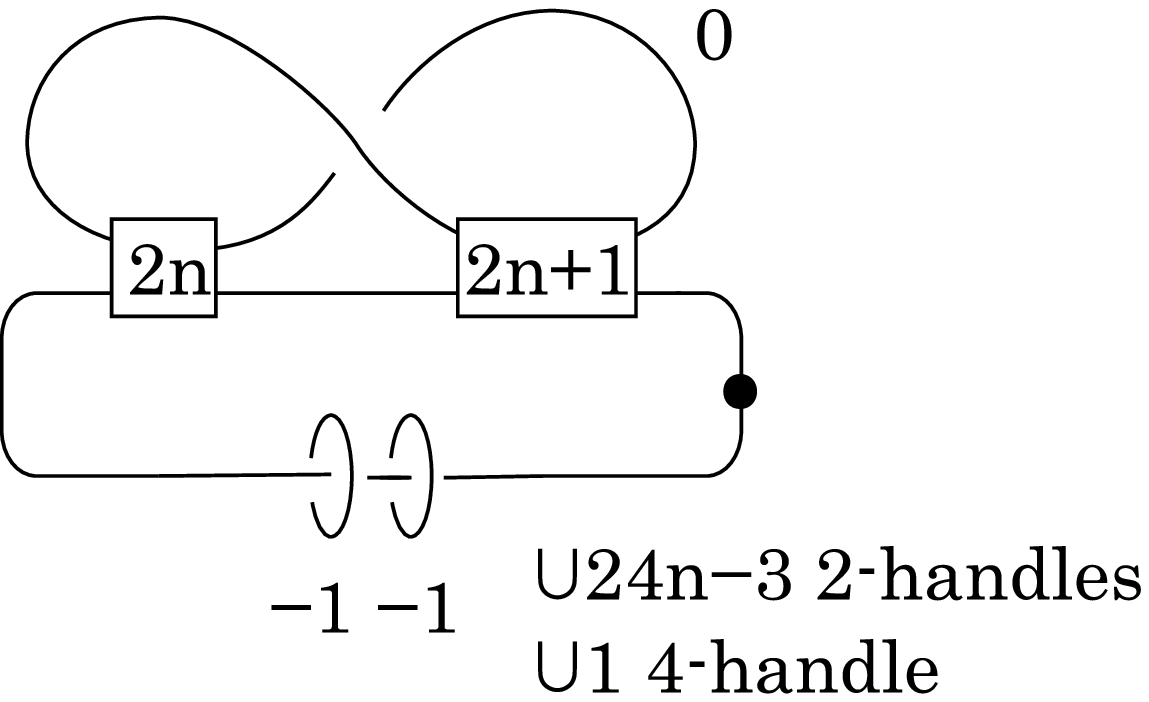}
\caption{$E(2n)\# \overline{\mathbf{CP}}^2$}
\label{fig30}
\end{center}
\end{figure}
\begin{figure}[ht!]
\begin{center}
\includegraphics[width=3.8in]{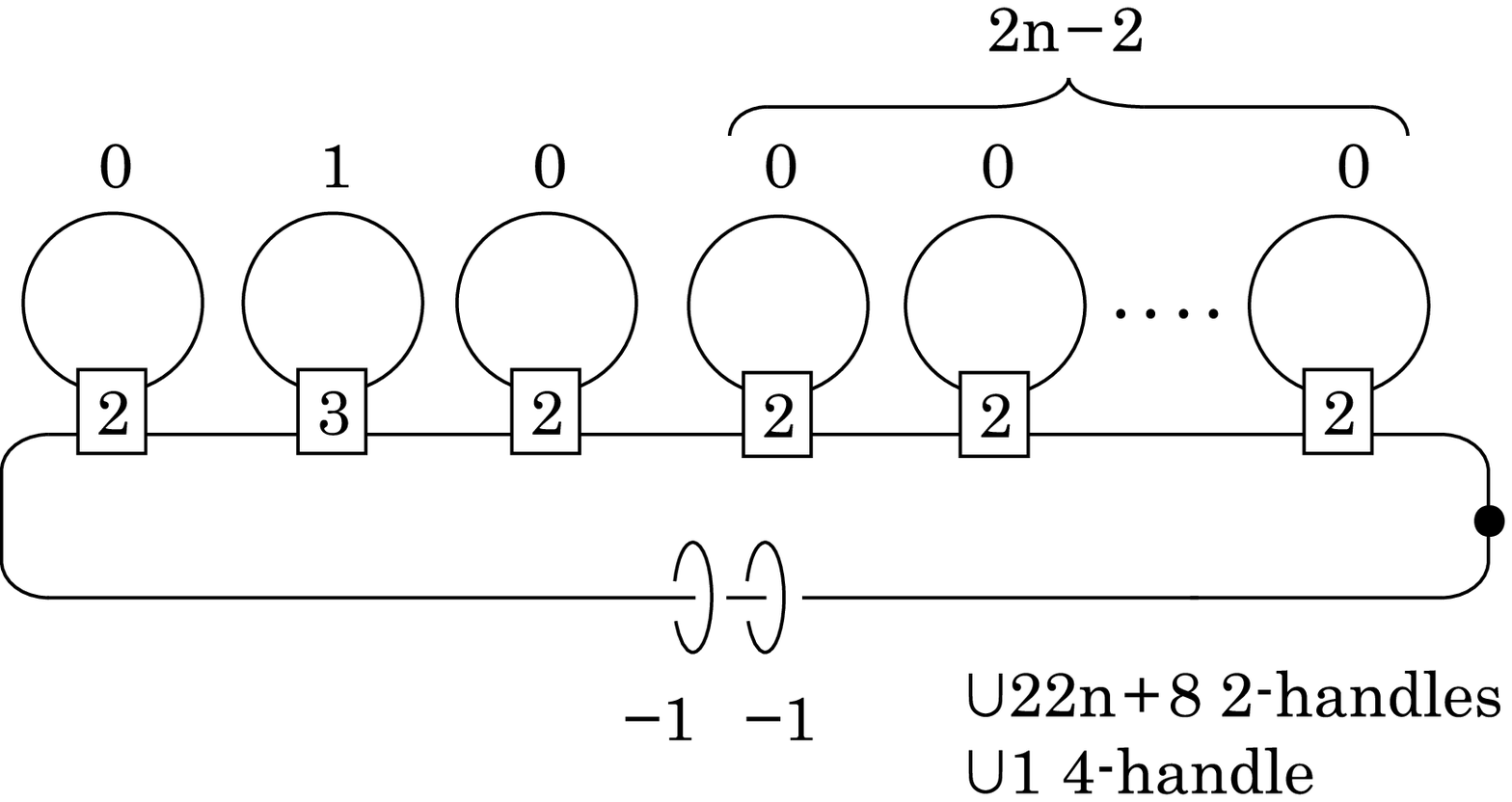}
\caption{$E(2n+1)$ $(n\geq 1)$}
\label{fig31}
\end{center}
\end{figure}
\begin{figure}[ht!]
\begin{center}
\includegraphics[width=4.1in]{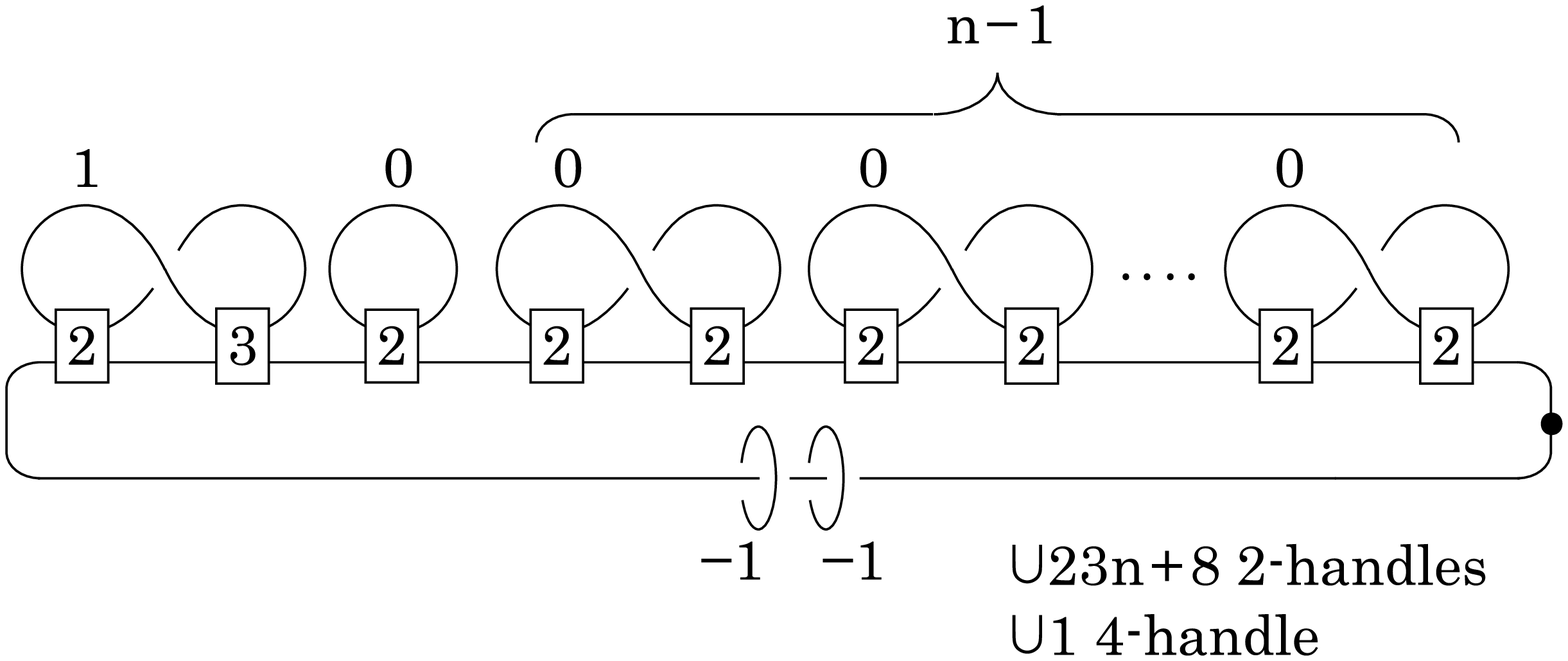}
\caption{$E(2n+1)$ $(n\geq 1)$}
\label{fig32}
\end{center}
\end{figure}
\begin{figure}[ht!]
\begin{center}
\includegraphics[width=4.1in]{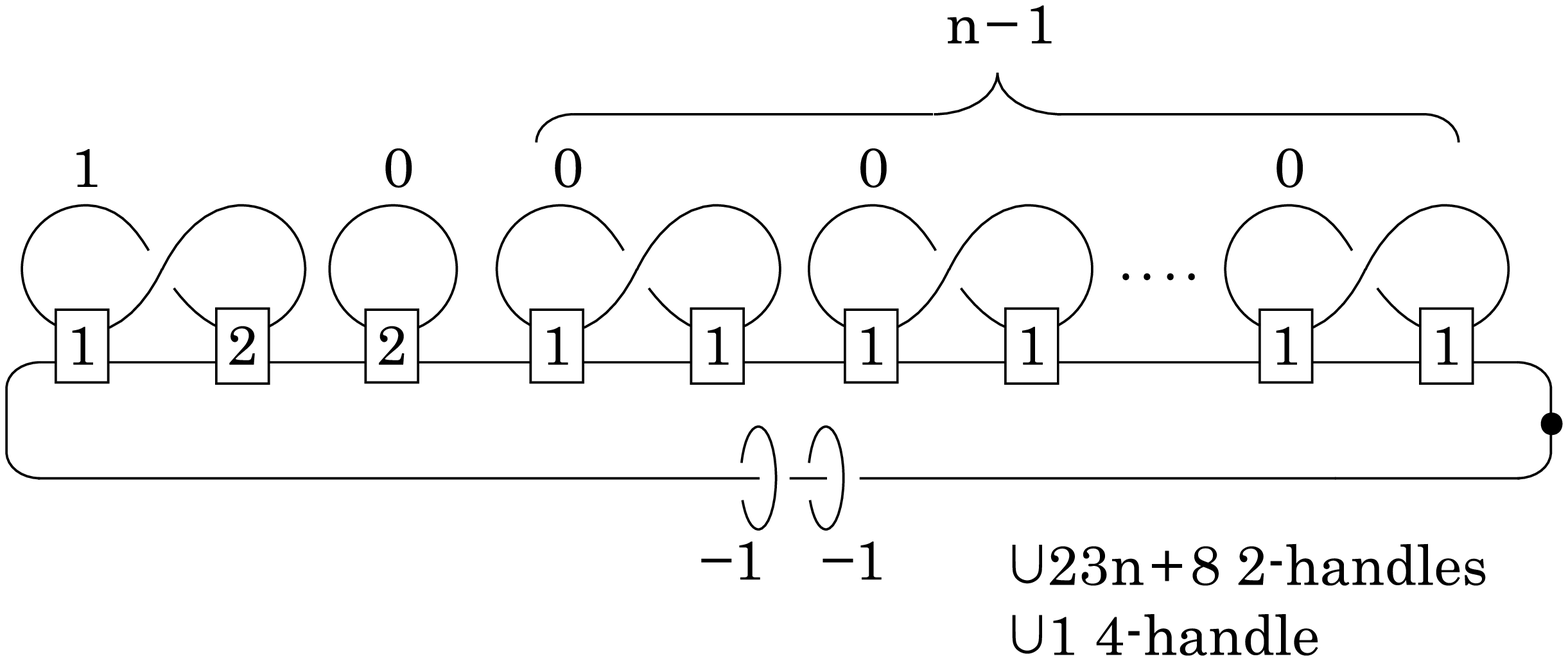}
\caption{$E(2n+1)$ $(n\geq 1)$}
\label{fig33}
\end{center}
\end{figure}
\begin{figure}[ht!]
\begin{center}
\includegraphics[width=3.9in]{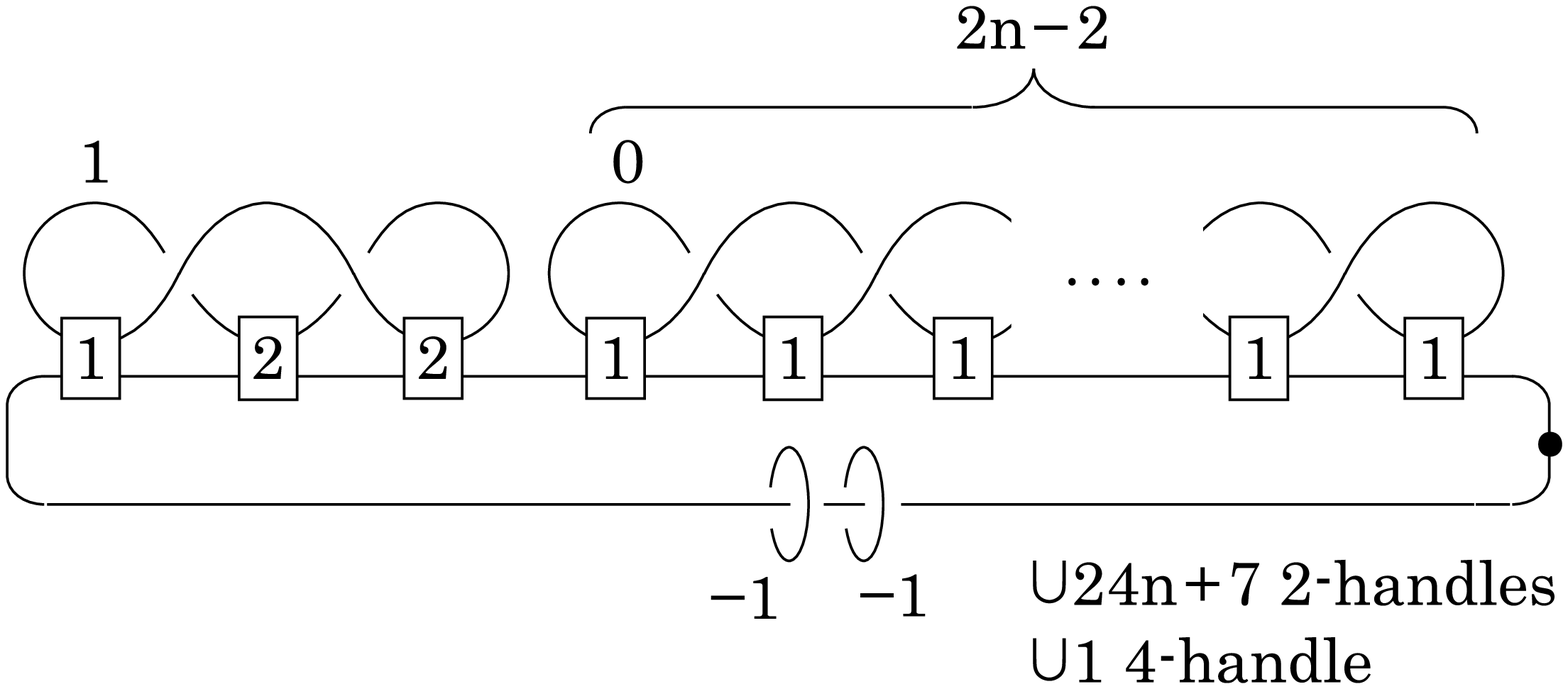}
\caption{$E(2n+1)$ $(n\geq 1)$}
\label{fig34}
\end{center}
\end{figure}
\begin{figure}[htb!]
\begin{center}
\includegraphics[width=4.1in]{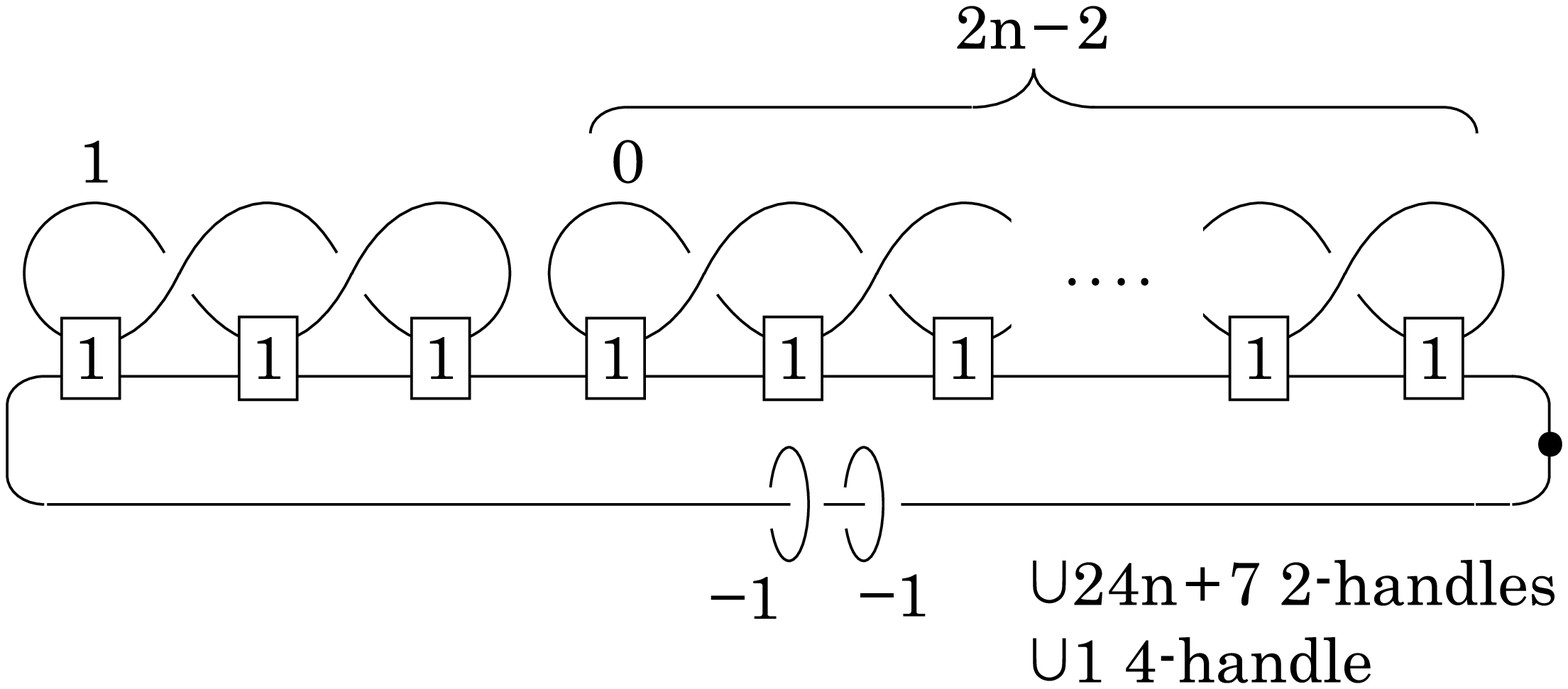}
\caption{$E(2n+1)$ $(n\geq 1)$}
\label{fig35}
\end{center}

\end{figure}
\begin{figure}[ht!]
\begin{center}
\includegraphics[width=2.6in]{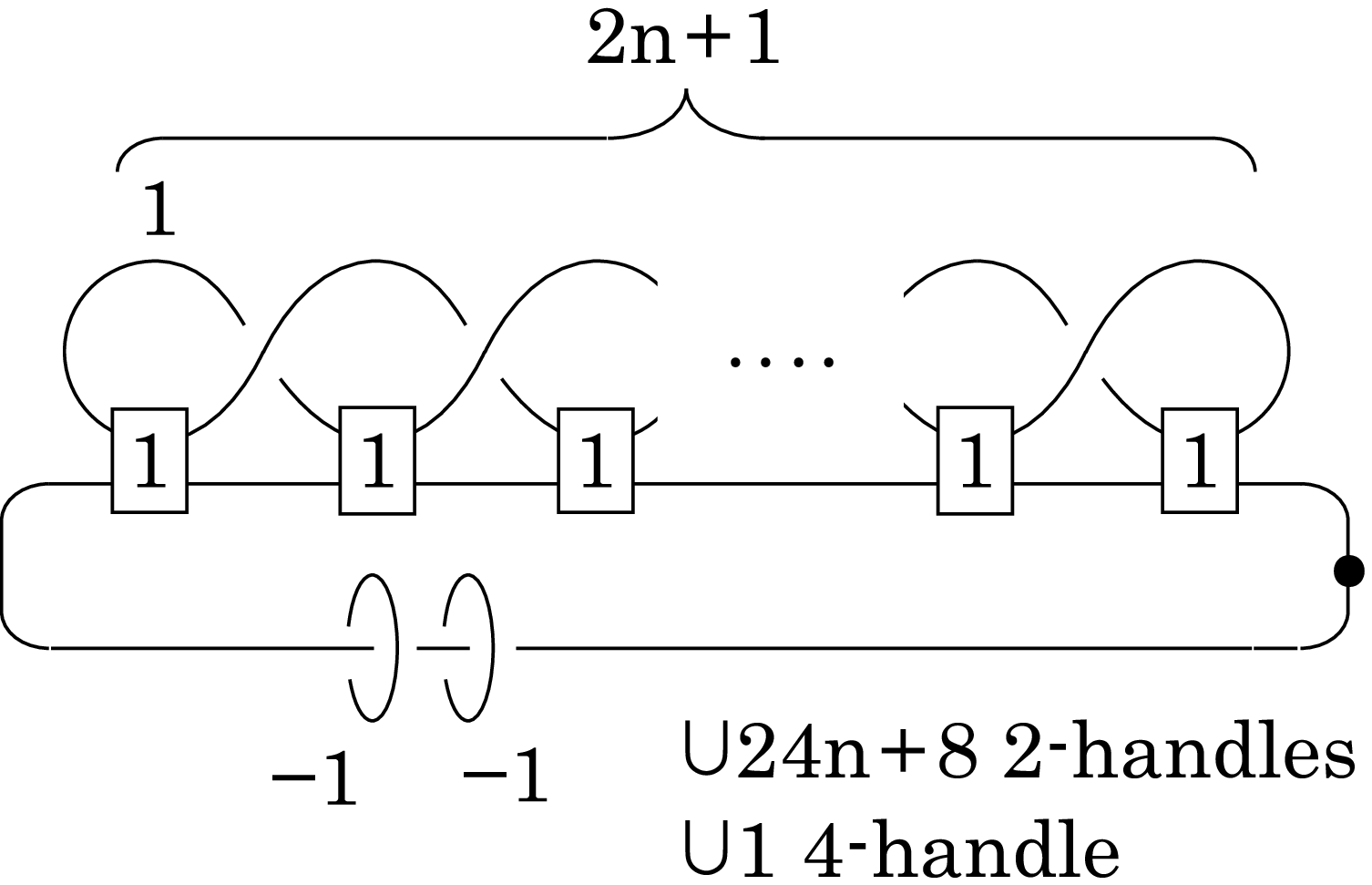}
\caption{$E(2n+1)$ $(n\geq 1)$}
\label{fig36}
\end{center}
\end{figure}
\begin{figure}[htb!]
\begin{center}
\includegraphics[width=4.8in]{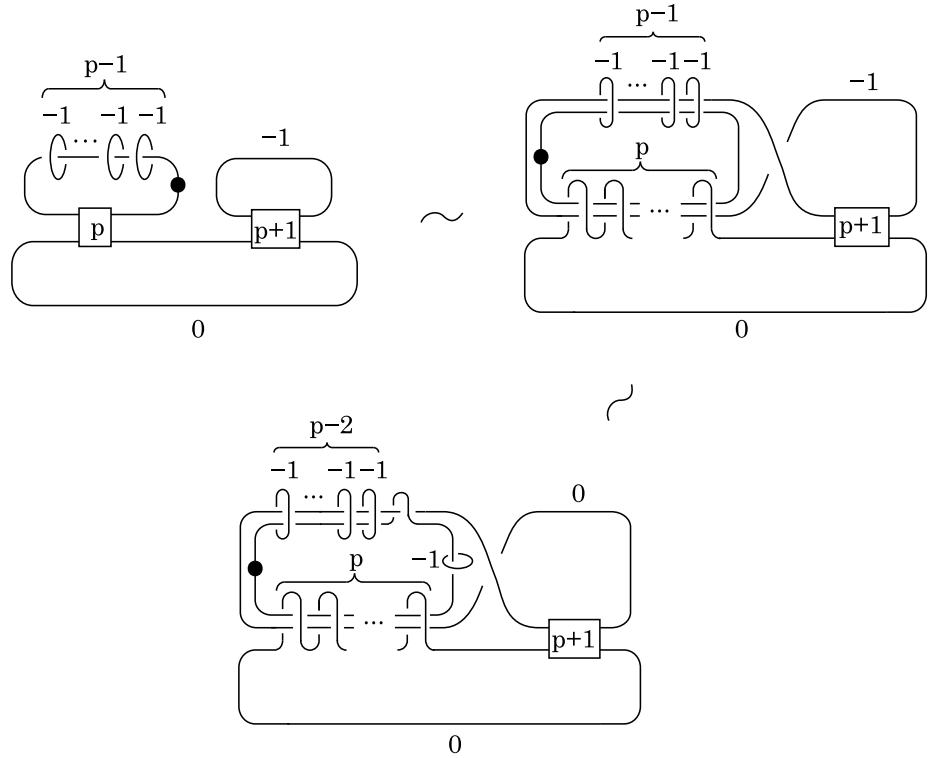}
\caption{Handle slides of $T_{p,p+1}$}
\label{fig37}
\end{center}
\end{figure}
\begin{figure}[htb!]
\begin{center}
\includegraphics[width=4.9in]{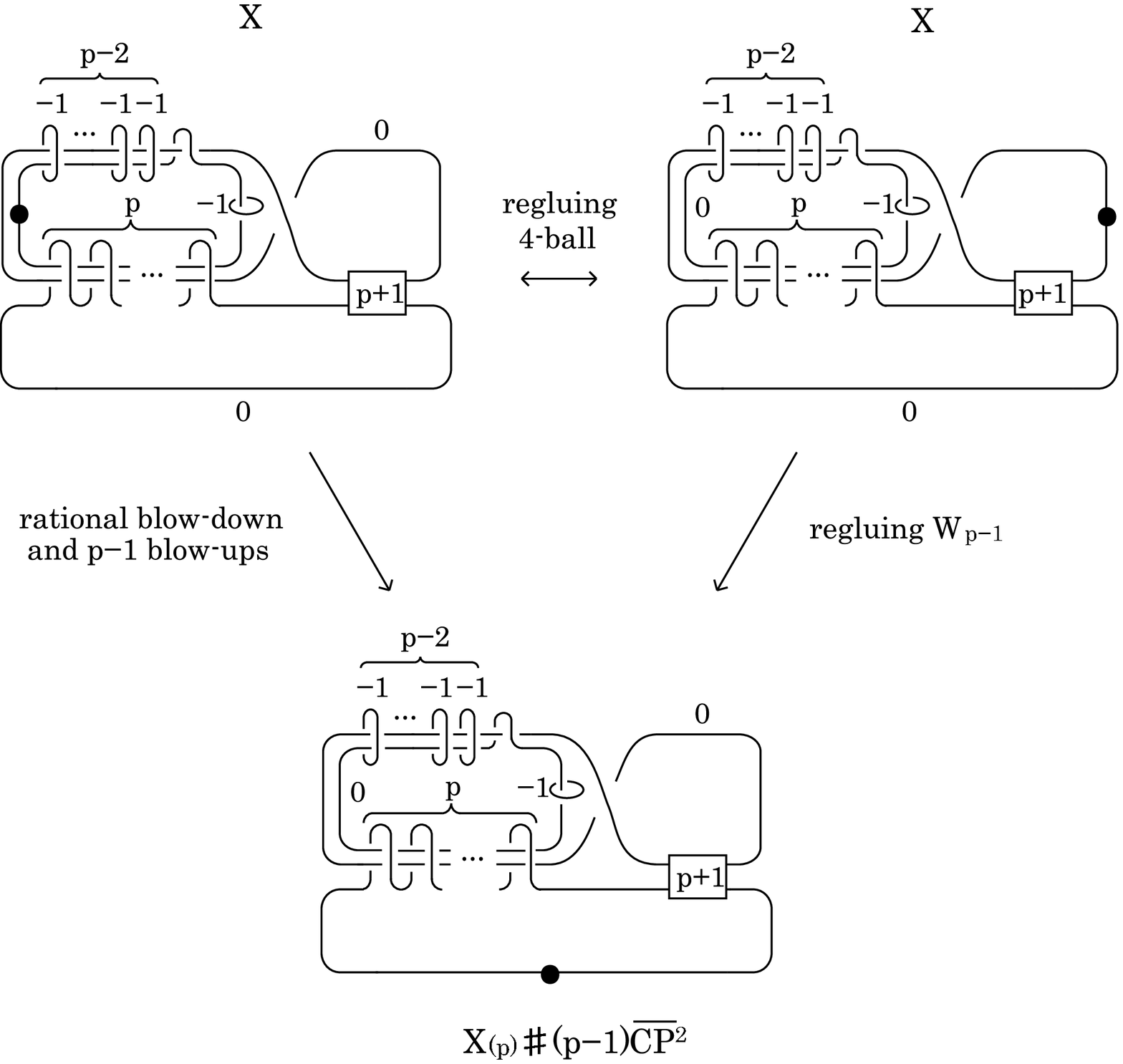}
\caption{}
\label{fig38}
\end{center}
\end{figure}
\begin{figure}[ht!]
\begin{center}
\includegraphics[width=1.5in]{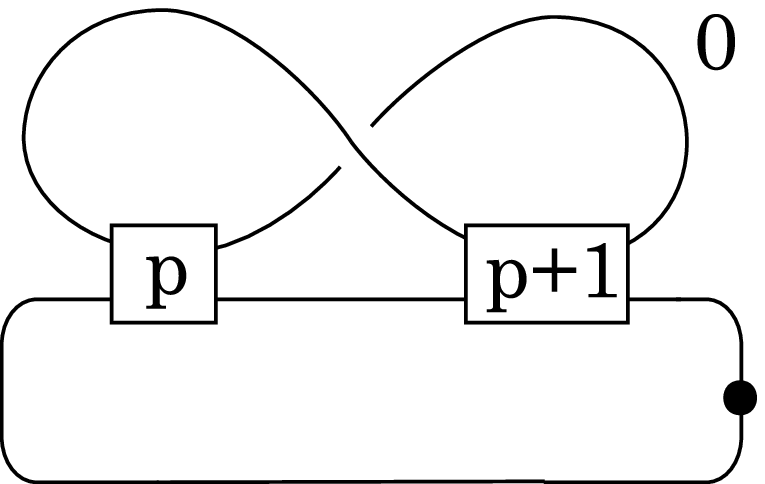}
\caption{}
\label{fig39}
\end{center}
\end{figure}
\begin{figure}[htb]
\begin{center}
\includegraphics[width=2.0in]{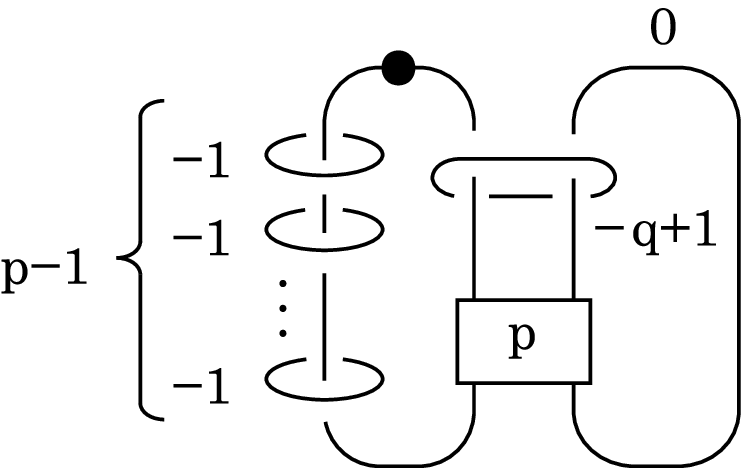}
\caption{}
\label{fig40}
\end{center}
\end{figure}
\begin{figure}
\begin{center}
\includegraphics[width=2.6in]{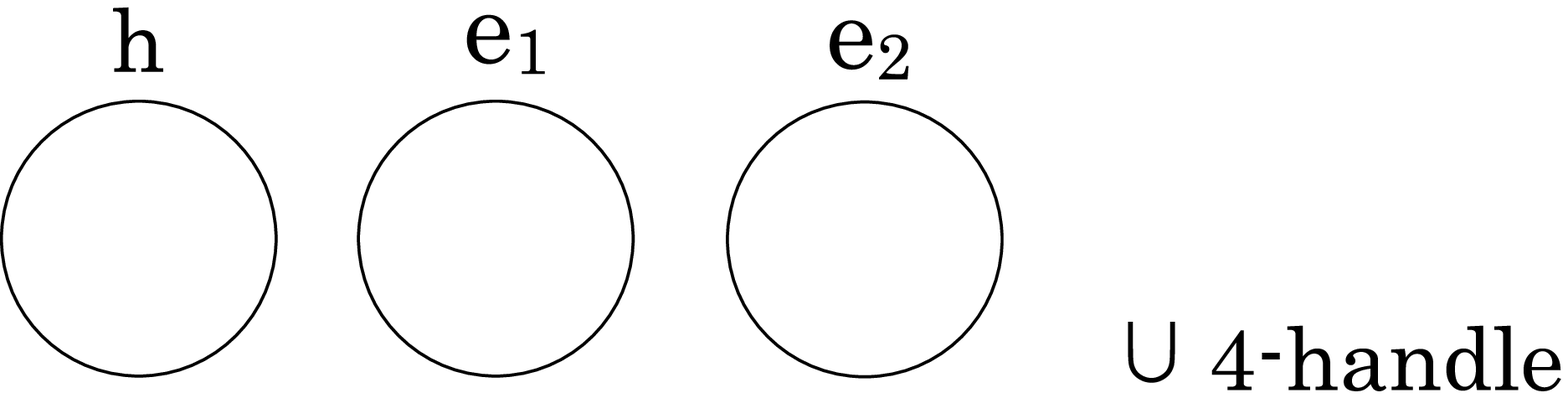}
\caption{$\mathbf{CP}^2\# 2\overline{\mathbf{CP}}^2$}
\label{fig41}
\end{center}
\end{figure}
\begin{figure}[ht!]
\begin{center}
\includegraphics[width=2.3in]{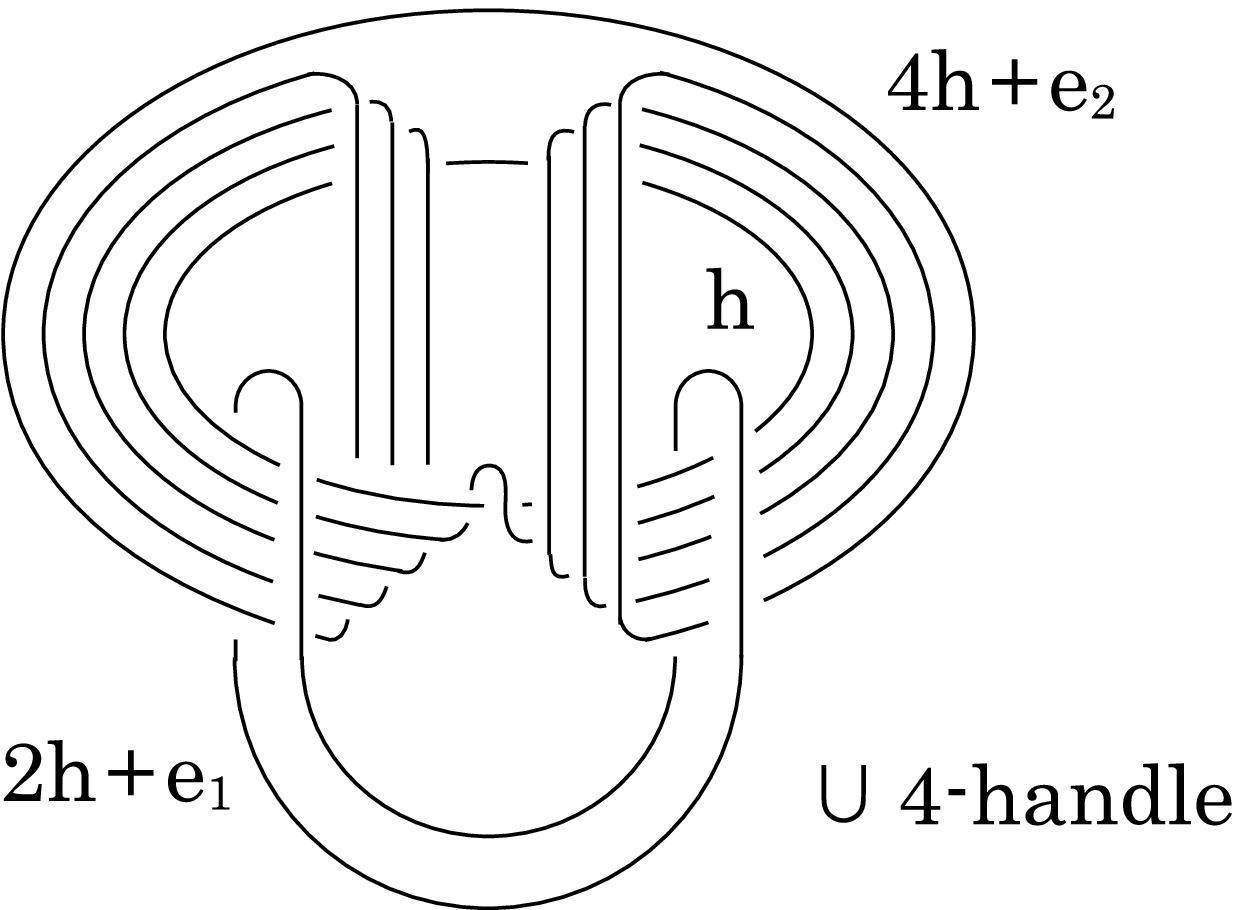}
\caption{$\mathbf{CP}^2\# 2\overline{\mathbf{CP}}^2$}
\label{fig42}
\end{center}
\end{figure}
\begin{figure}[ht!]
\begin{center}
\includegraphics[width=2.7in]{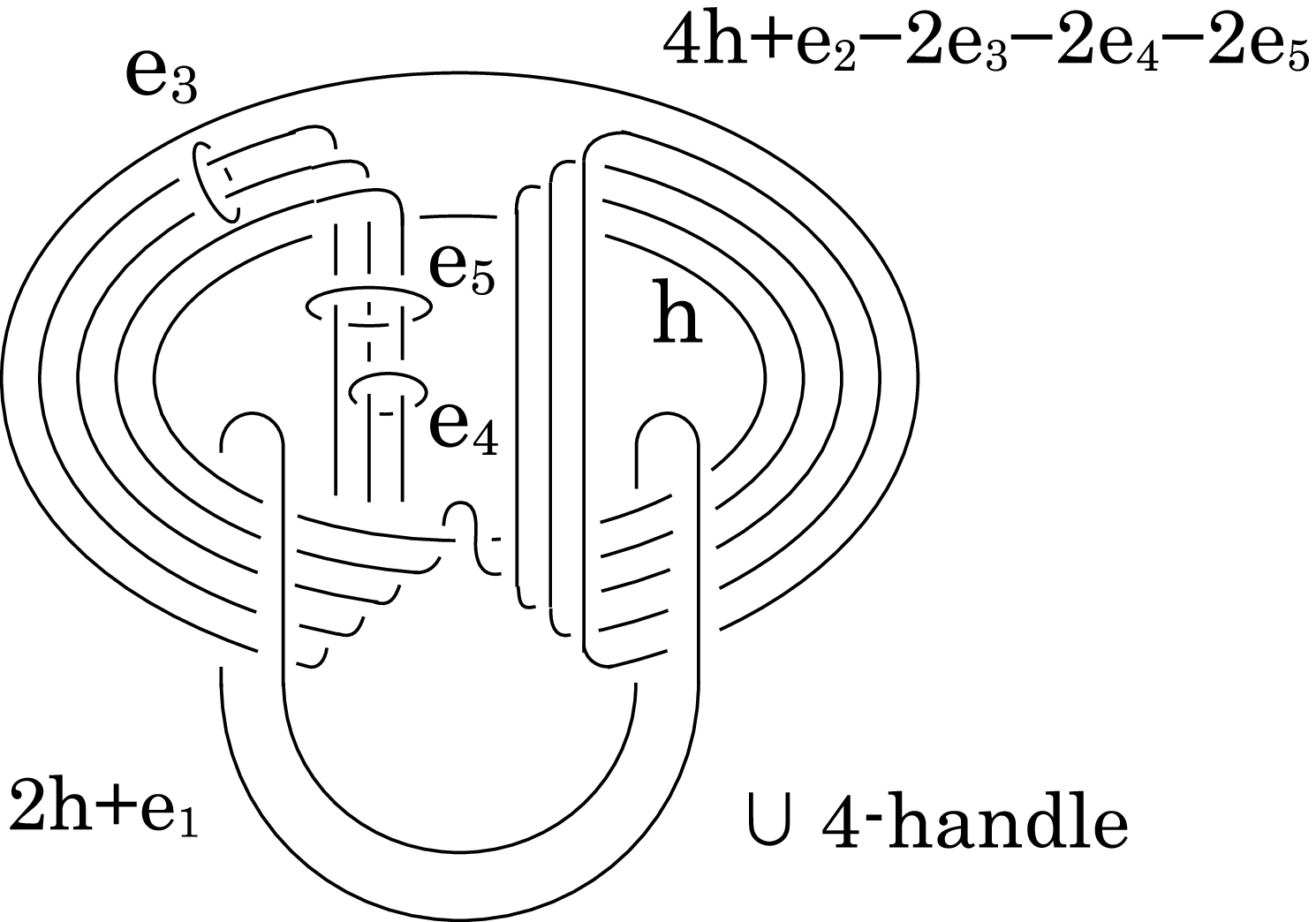}
\caption{$\mathbf{CP}^2\# 5\overline{\mathbf{CP}}^2$}
\label{fig43}
\end{center}
\end{figure}
\begin{figure}[htb!]
\begin{center}
\includegraphics[width=2.8in]{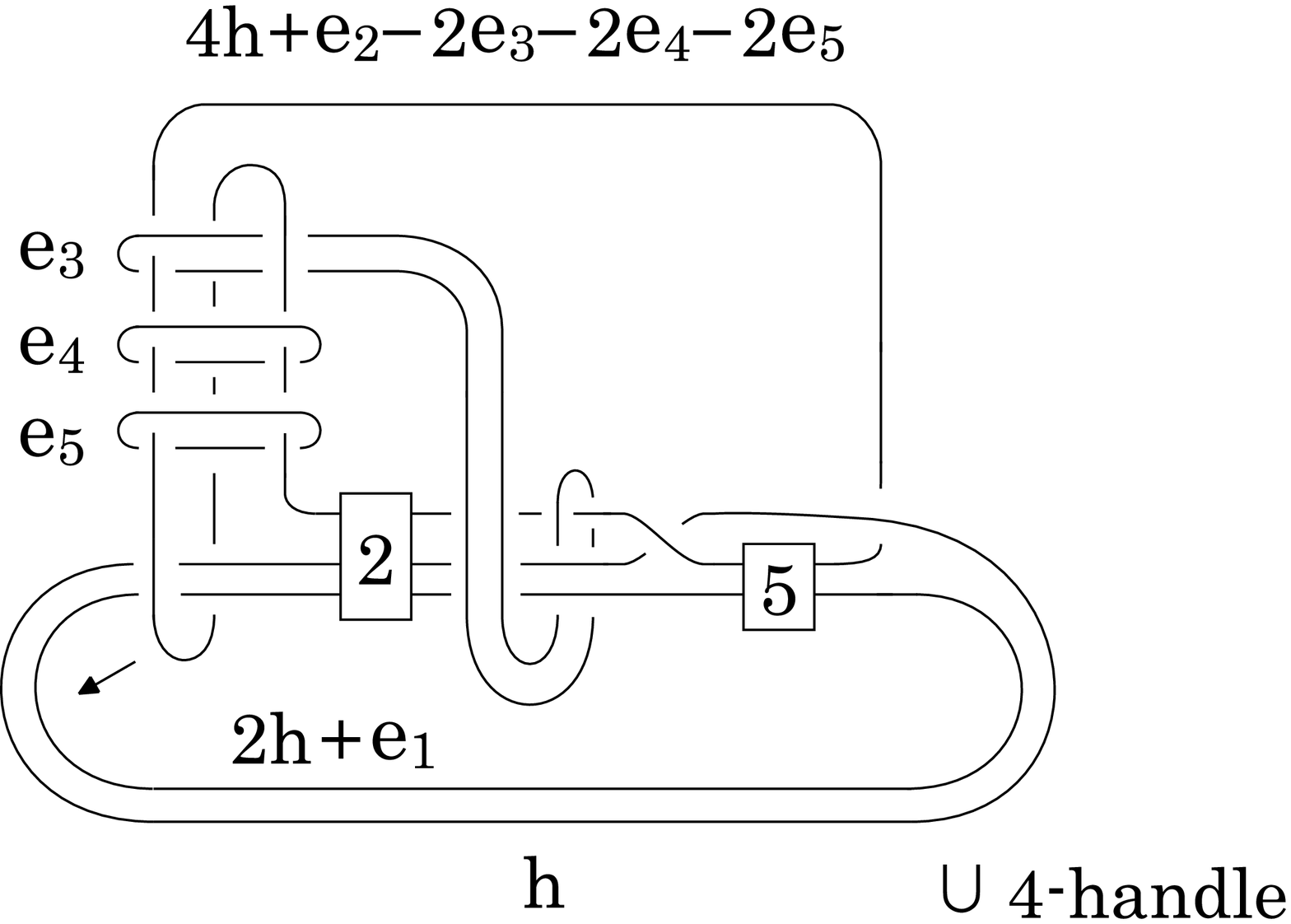}
\caption{$\mathbf{CP}^2\# 5\overline{\mathbf{CP}}^2$}
\label{fig44}
\end{center}
\end{figure}
\begin{figure}
\begin{center}
\includegraphics[width=2.9in]{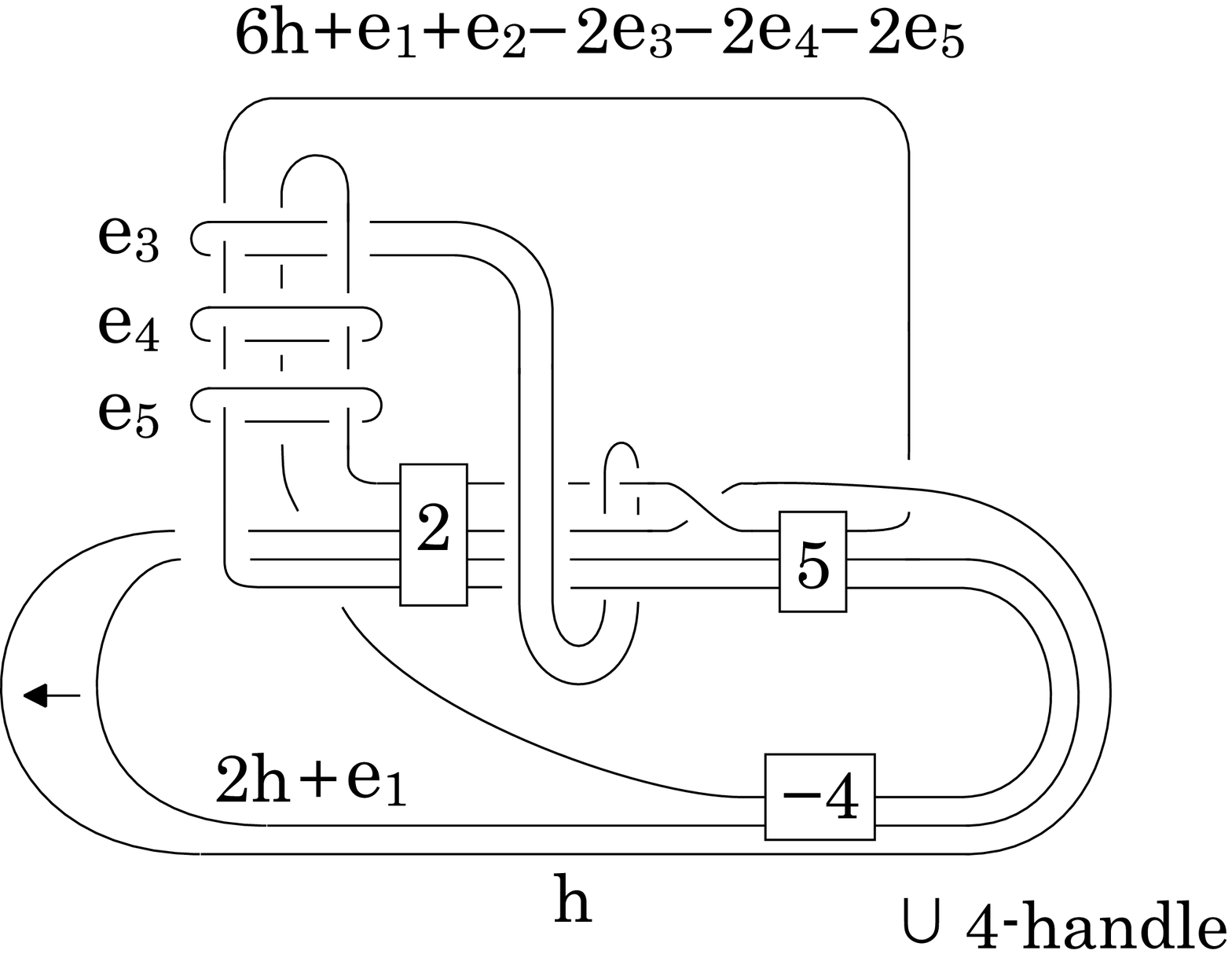}
\caption{$\mathbf{CP}^2\# 5\overline{\mathbf{CP}}^2$}
\label{fig45}
\end{center}
\end{figure}
\begin{figure}
\begin{center}
\includegraphics[width=3.8in]{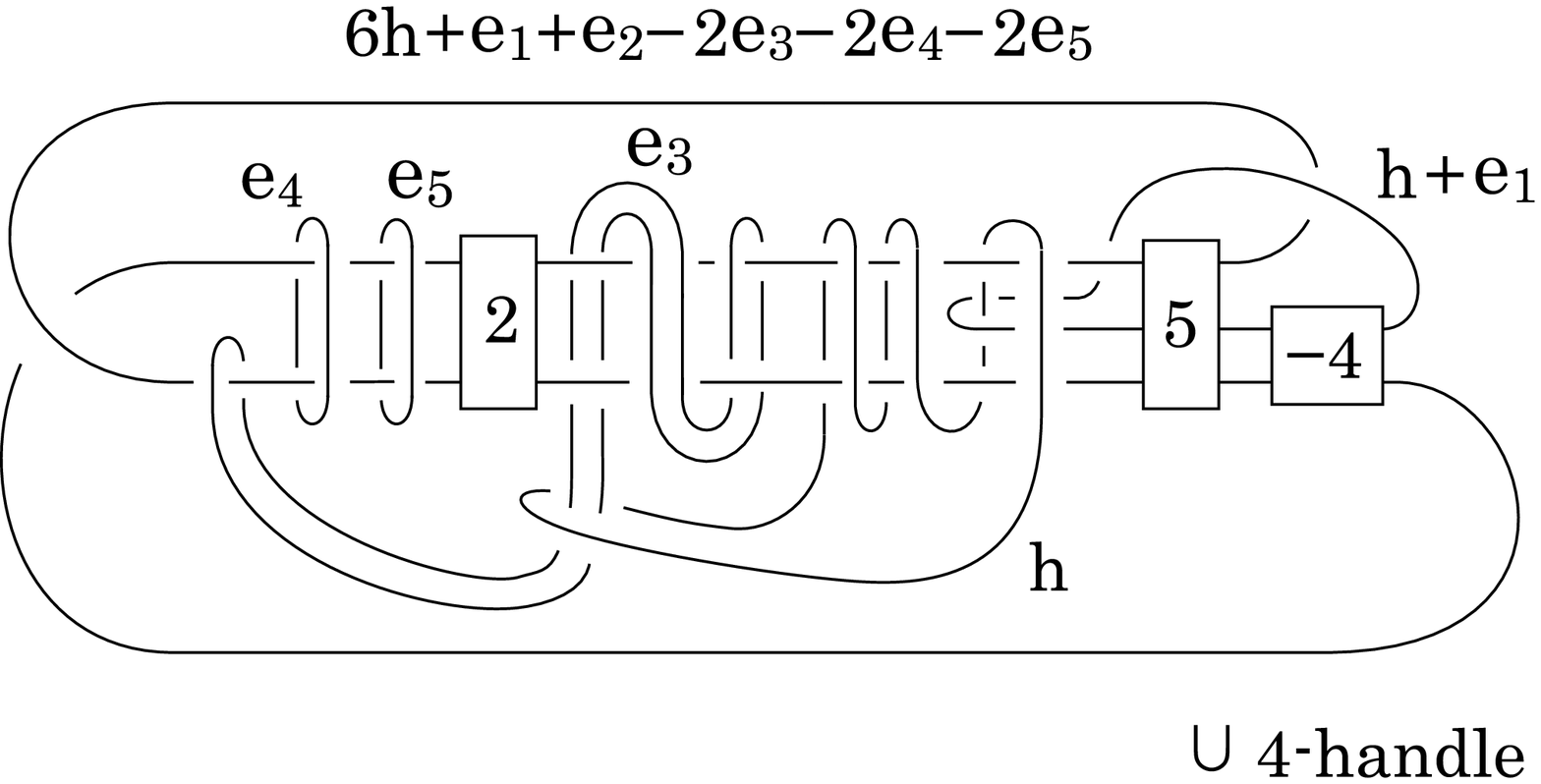}
\caption{$\mathbf{CP}^2\# 5\overline{\mathbf{CP}}^2$}
\label{fig46}
\end{center}
\end{figure}
\begin{figure}[htb!]
\begin{center}
\includegraphics[width=4.5in]{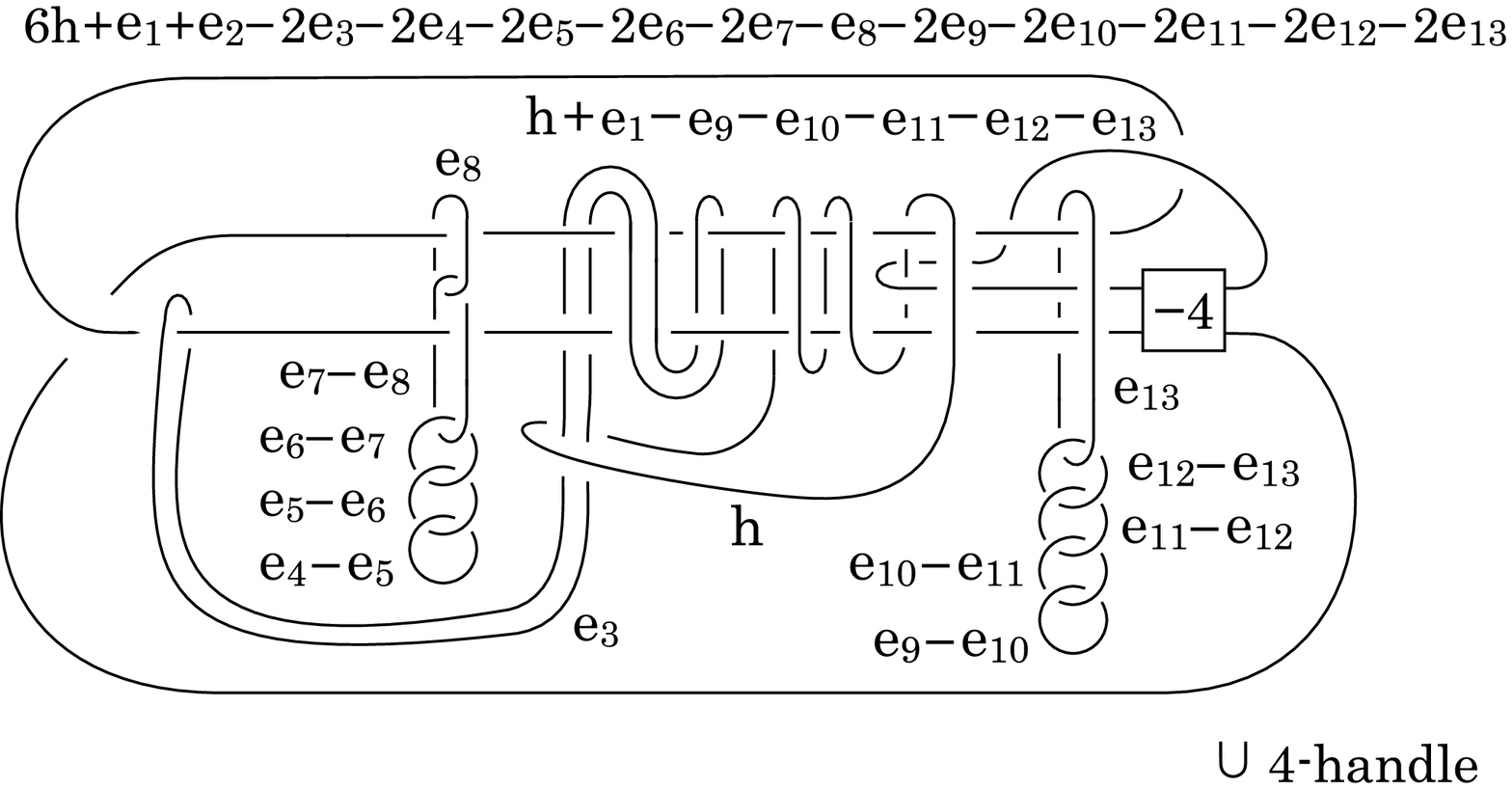}
\caption{$\mathbf{CP}^2\# 13\overline{\mathbf{CP}}^2$}
\label{fig47}
\end{center}
\end{figure}
\begin{figure}[htb!]
\begin{center}
\includegraphics[width=3.4in]{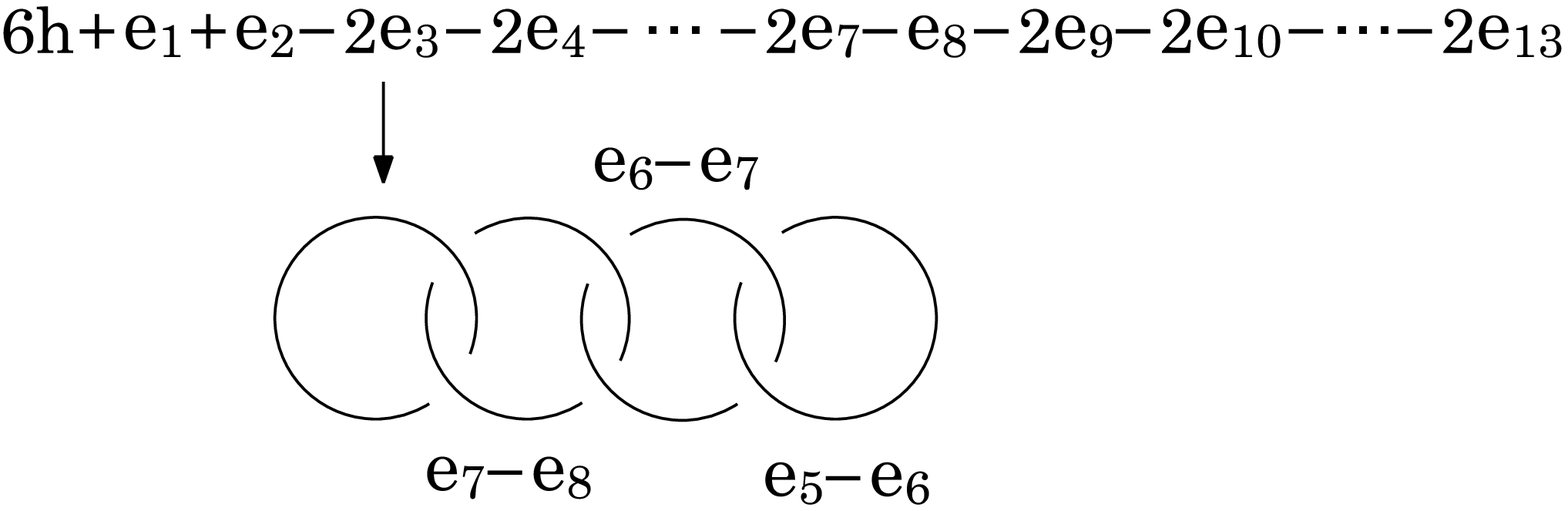}
\caption{$C_5$ in $\mathbf{CP}^2\# 13\overline{\mathbf{CP}}^2$}
\label{fig48}
\end{center}
\end{figure}
\begin{figure}[htb!]
\begin{center}
\includegraphics[width=2.3in]{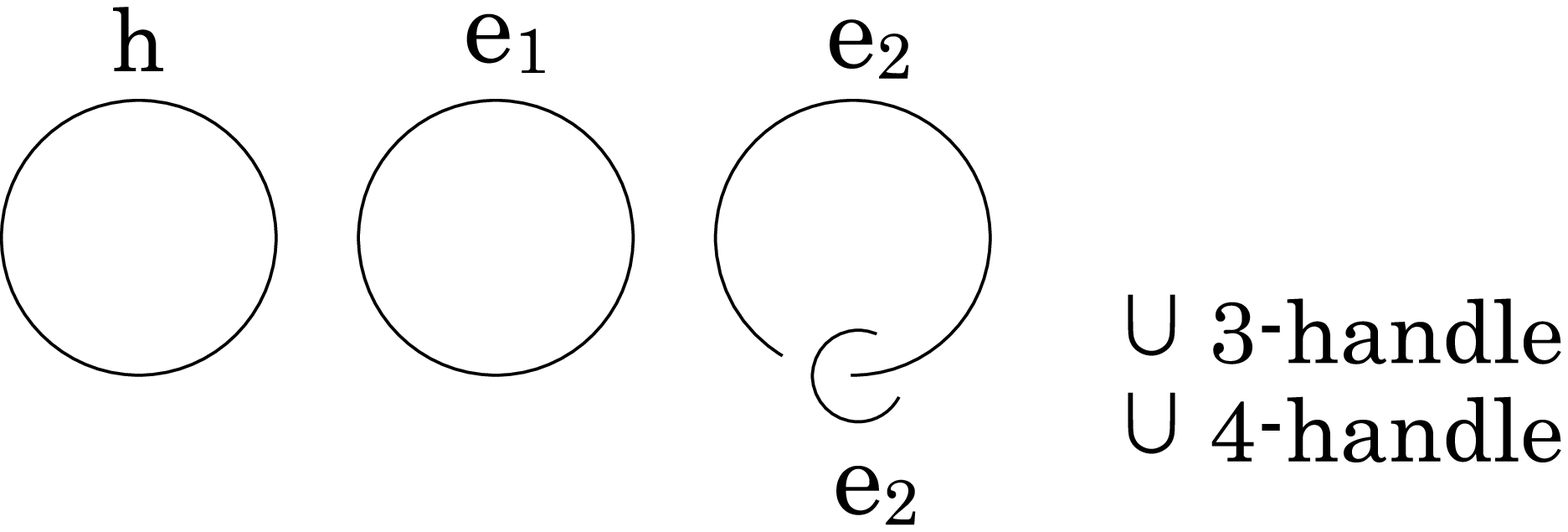}
\caption{$\mathbf{CP}^2\# 2\overline{\mathbf{CP}}^2$}
\label{fig49}
\end{center}
\end{figure}
\begin{figure}
\begin{center}
\includegraphics[width=4.5in]{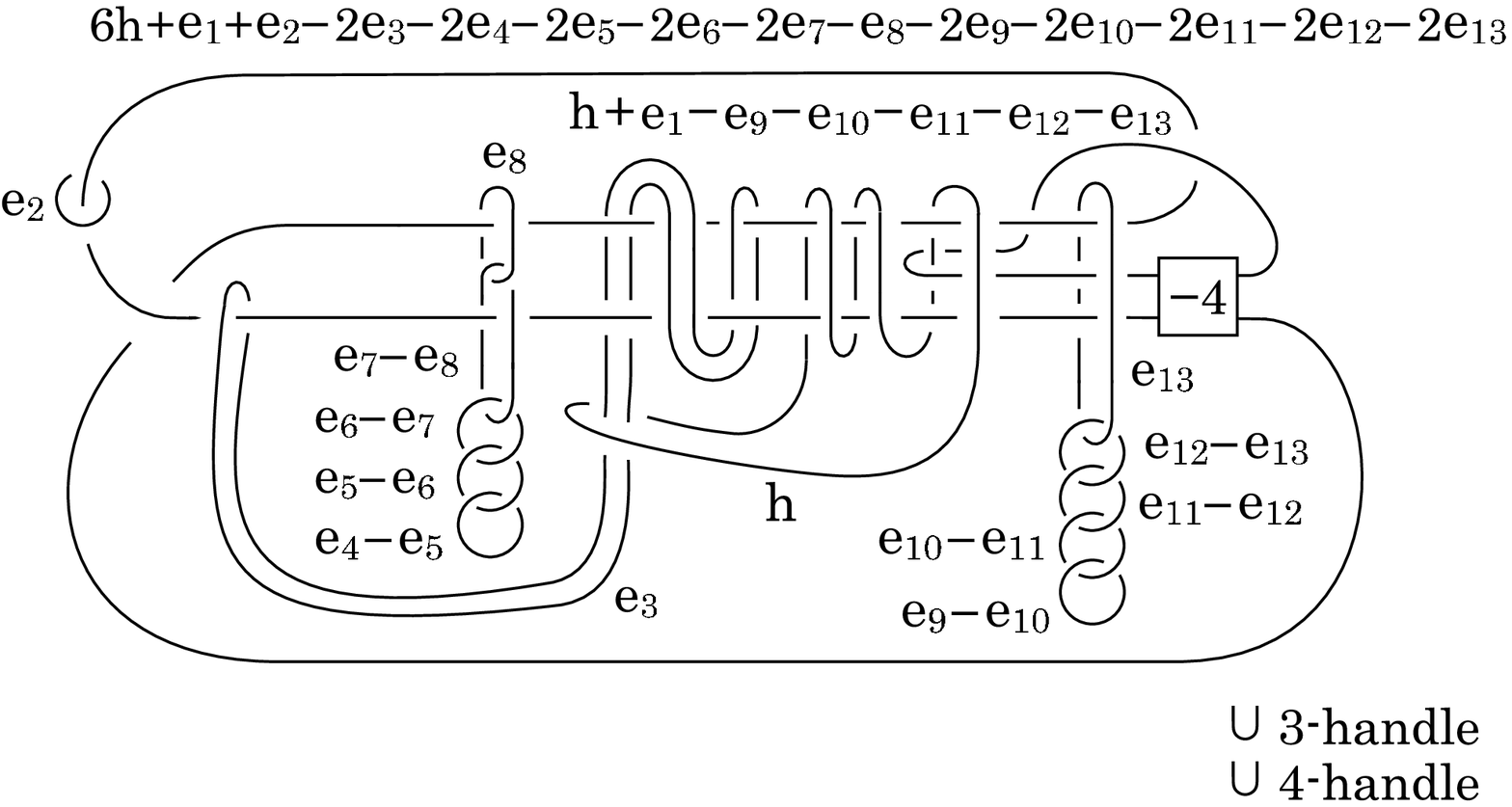}
\caption{$\mathbf{CP}^2\# 13\overline{\mathbf{CP}}^2$}
\label{fig50}
\end{center}
\end{figure}
\begin{figure}[ht!]
\begin{center}
\includegraphics[width=4.4in]{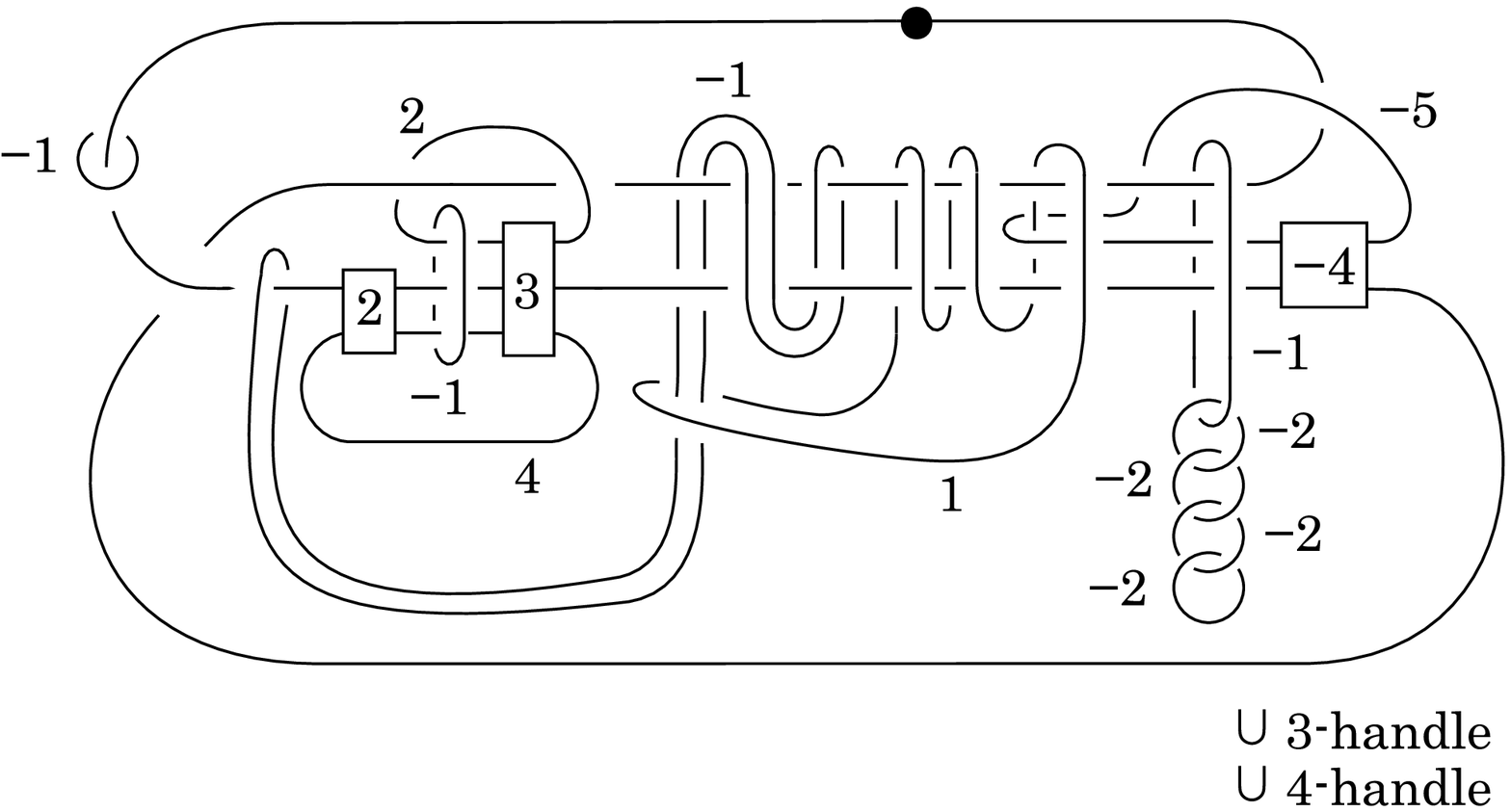}
\caption{A different diagram of the $4$-manifold in Figure~\ref{fig21}}
\label{fig51}
\end{center}\medskip \medskip 
\end{figure}
\begin{figure}[ht!]
\begin{center}
\includegraphics[width=4.8in]{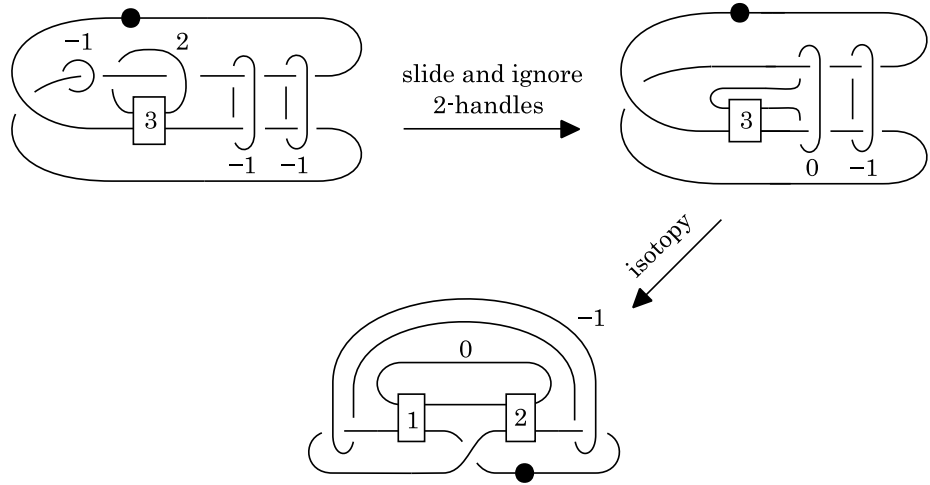}
\caption{Construction of a cork in $E'_3$}
\label{fig52}
\end{center}
\end{figure}
\begin{figure}[ht!]
\begin{center}
\includegraphics[width=4.8in]{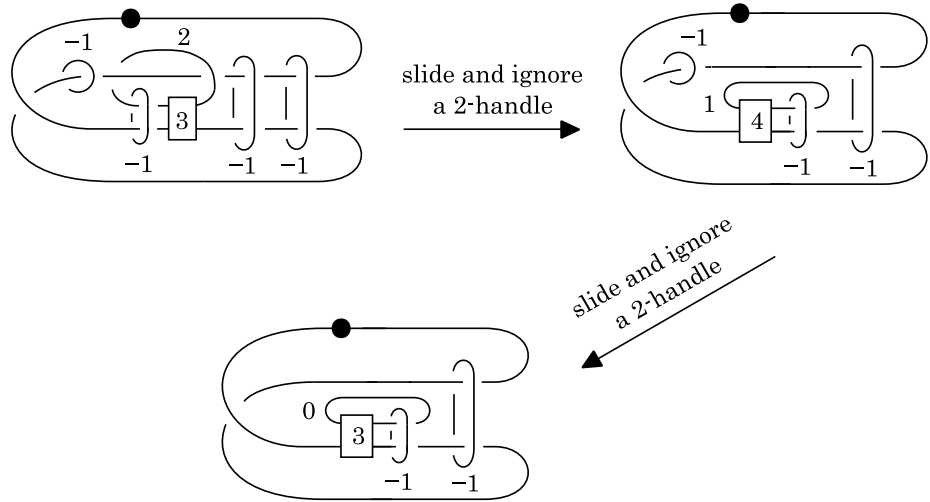}
\caption{Construction of a plug in $E'_3$}
\label{fig53}
\end{center}
\end{figure}
\begin{figure}[htb!]
\begin{center}
\includegraphics[width=2.5in]{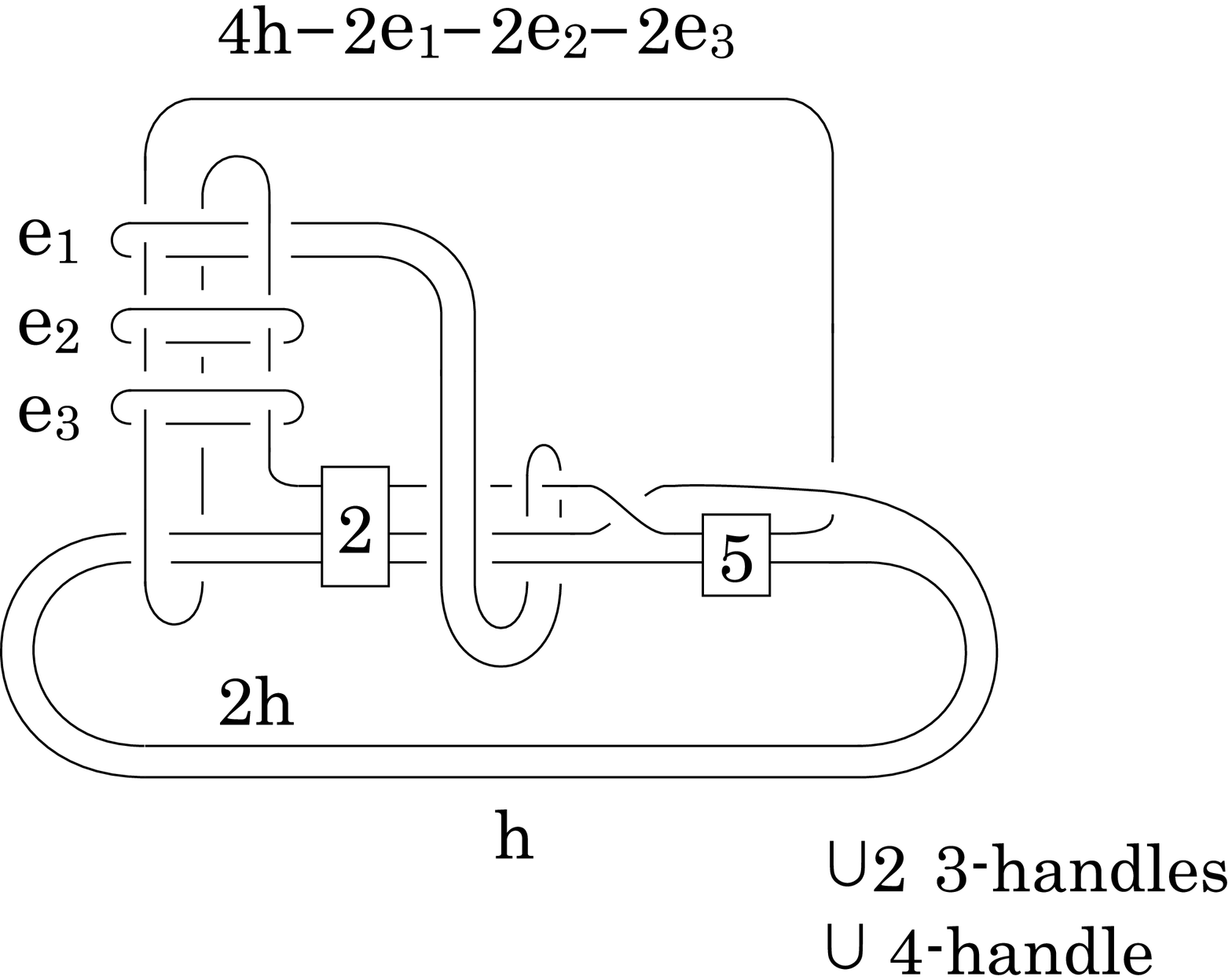}
\caption{$\mathbf{CP}^2\# 3\overline{\mathbf{CP}}^2$}
\label{fig54}
\end{center}
\end{figure}
\begin{figure}
\begin{center}
\includegraphics[width=2.3in]{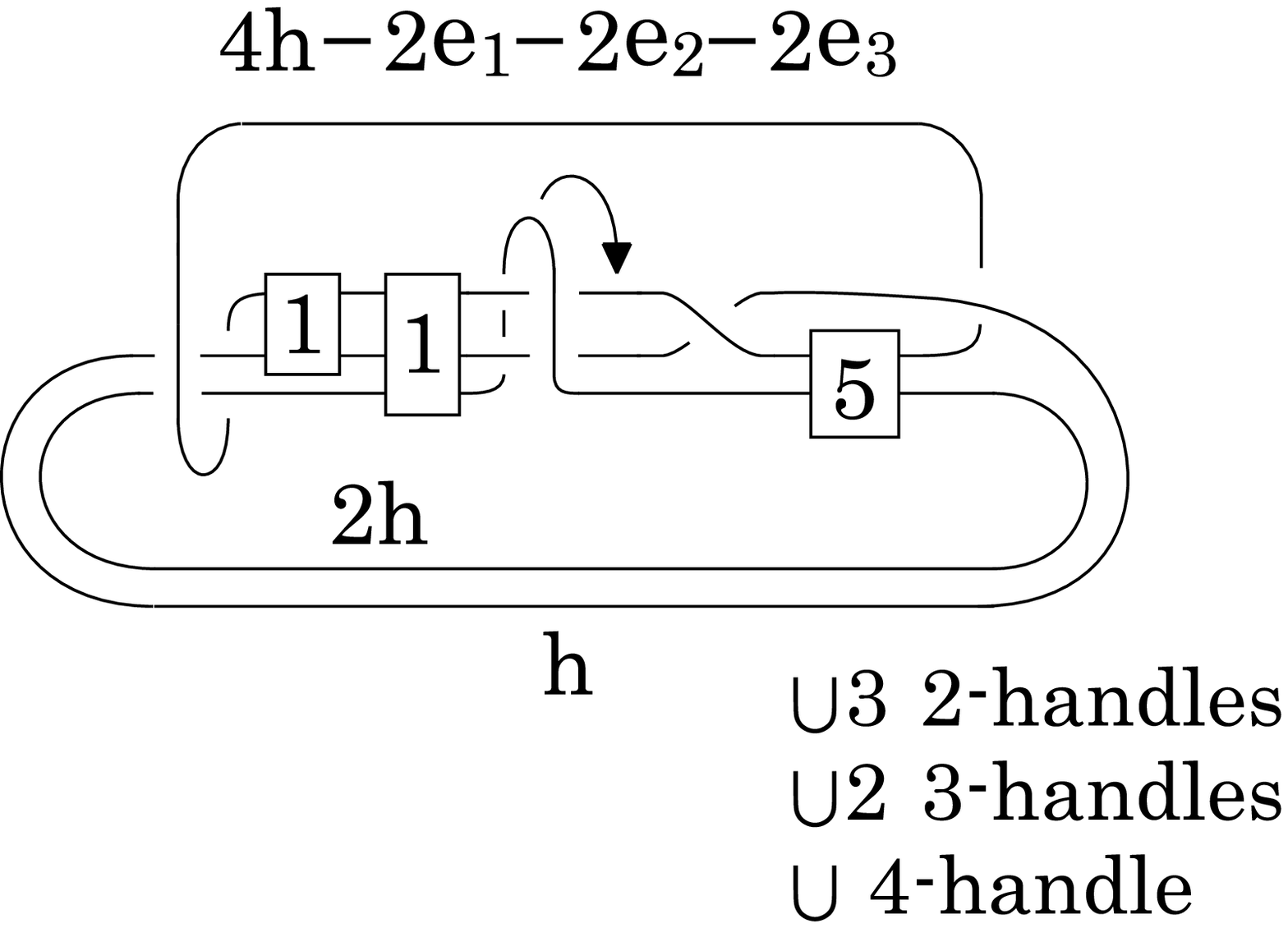}
\caption{$\mathbf{CP}^2\# 3\overline{\mathbf{CP}}^2$}
\label{fig55}
\end{center}
\end{figure}
\begin{figure}[ht!]
\begin{center}
\includegraphics[width=2.2in]{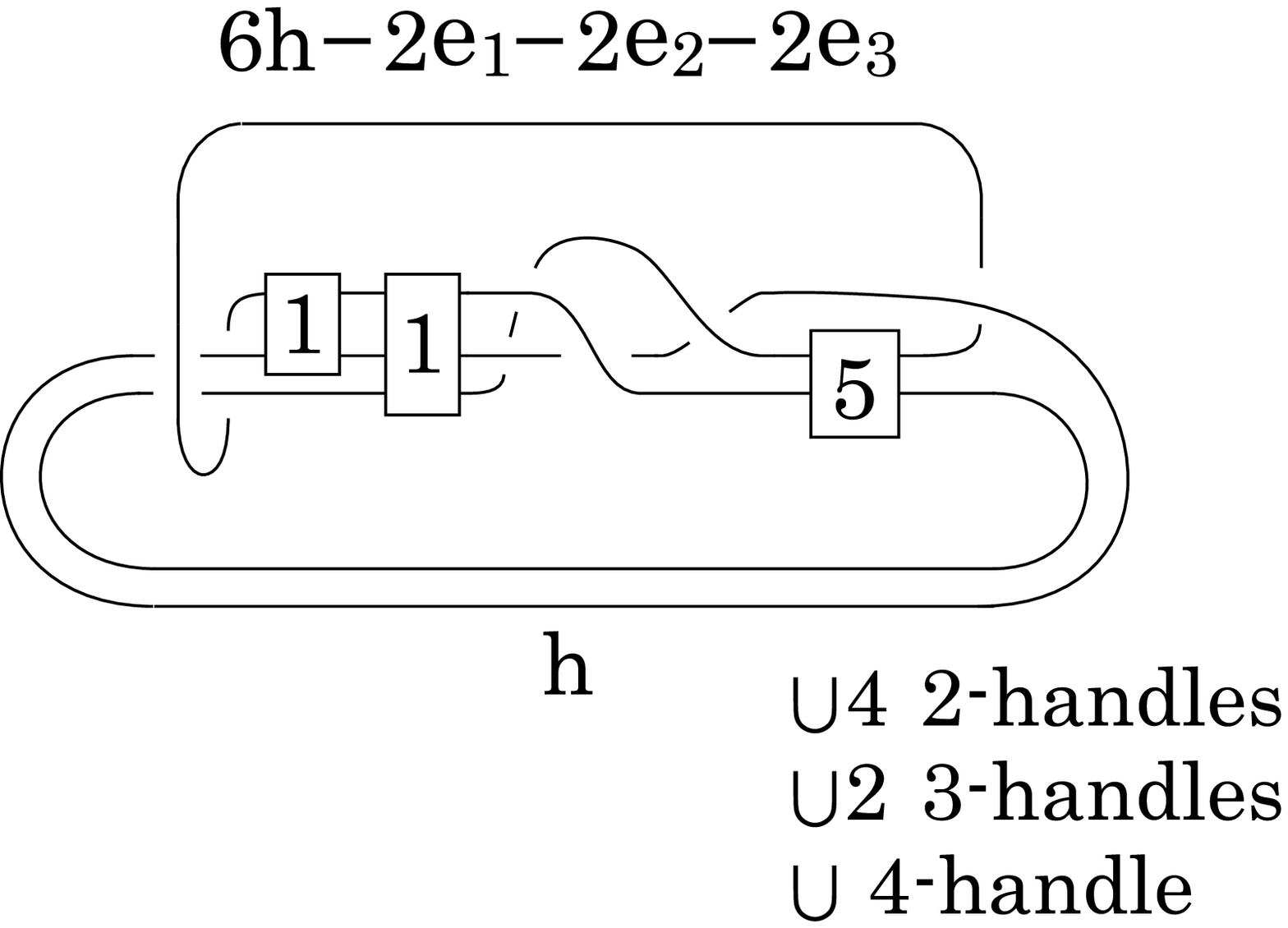}
\caption{$\mathbf{CP}^2\# 3\overline{\mathbf{CP}}^2$}
\label{fig56}
\end{center}
\end{figure}
\begin{figure}[ht!]
\begin{center}
\includegraphics[width=2.6in]{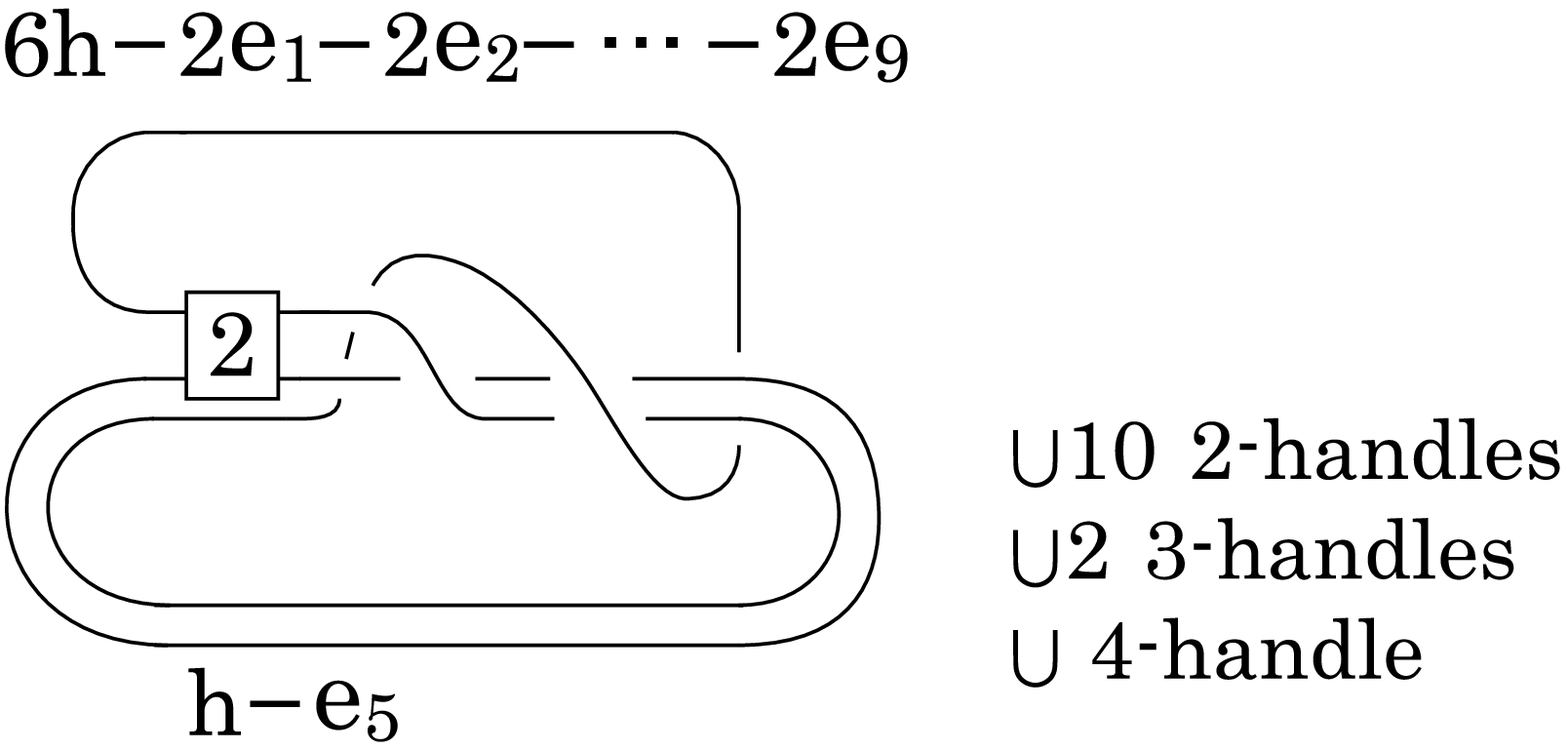}
\caption{$\mathbf{CP}^2\# 9\overline{\mathbf{CP}}^2$}
\label{fig57}
\end{center}
\end{figure}
\begin{figure}[ht!]
\begin{center}
\includegraphics[width=2.5in]{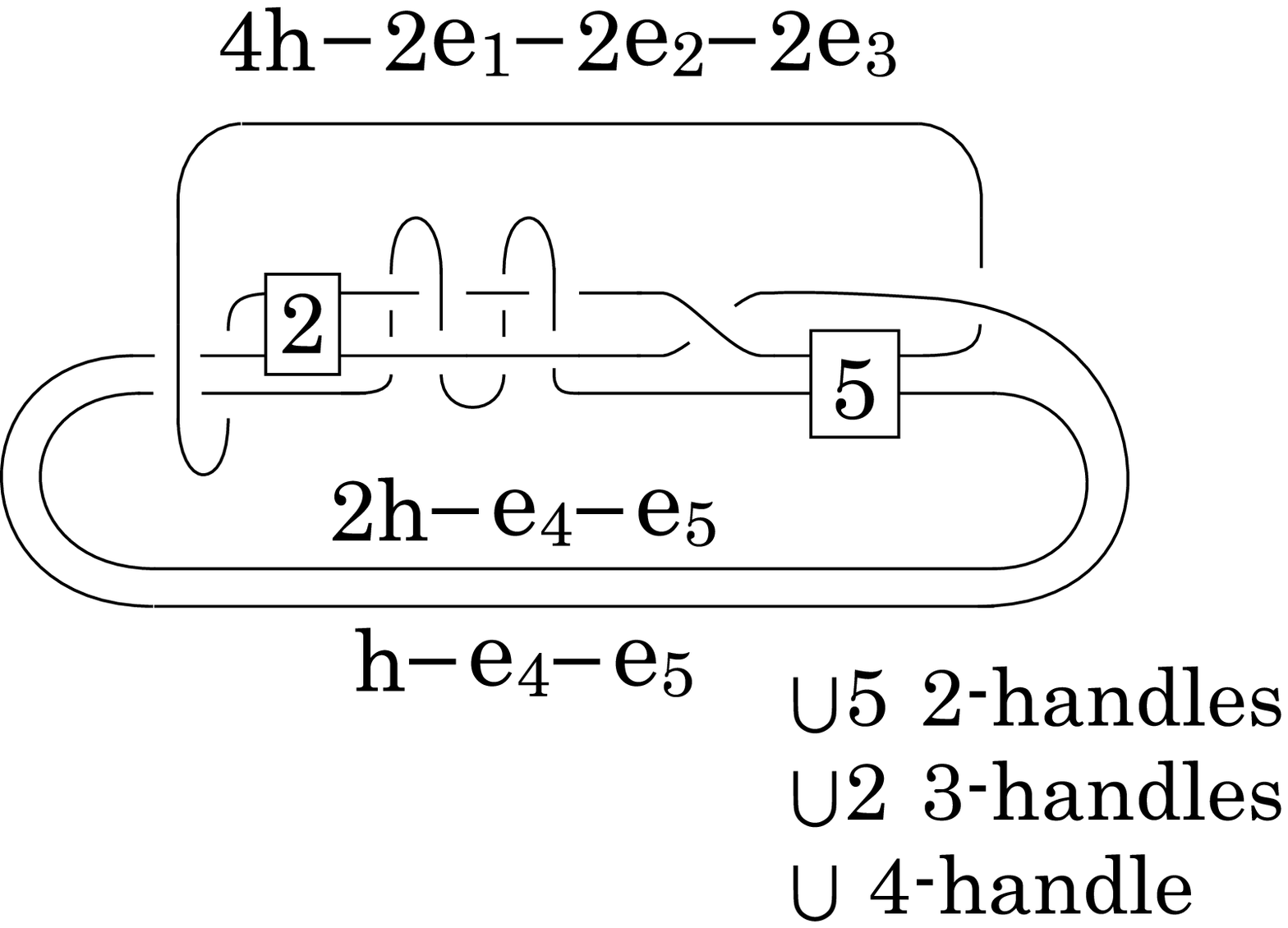}
\caption{$\mathbf{CP}^2\# 5\overline{\mathbf{CP}}^2$}
\label{fig58}
\end{center}
\end{figure}
\begin{figure}[ht!]
\begin{center}
\includegraphics[width=3.0in]{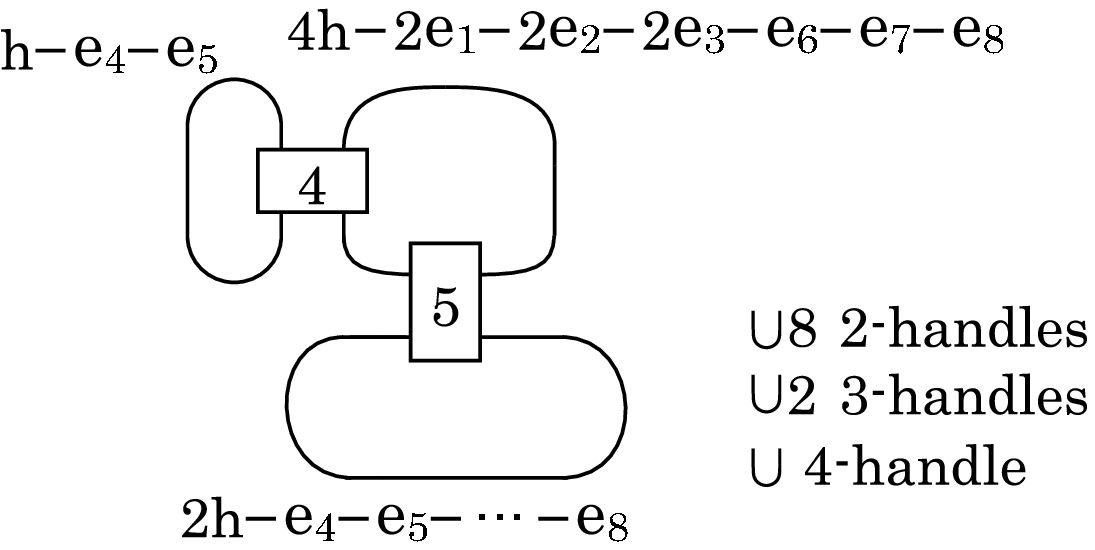}
\caption{$\mathbf{CP}^2\# 8\overline{\mathbf{CP}}^2$}
\label{fig59}
\end{center}
\end{figure}
\begin{figure}[htb!]
\begin{center}
\includegraphics[width=4.9in]{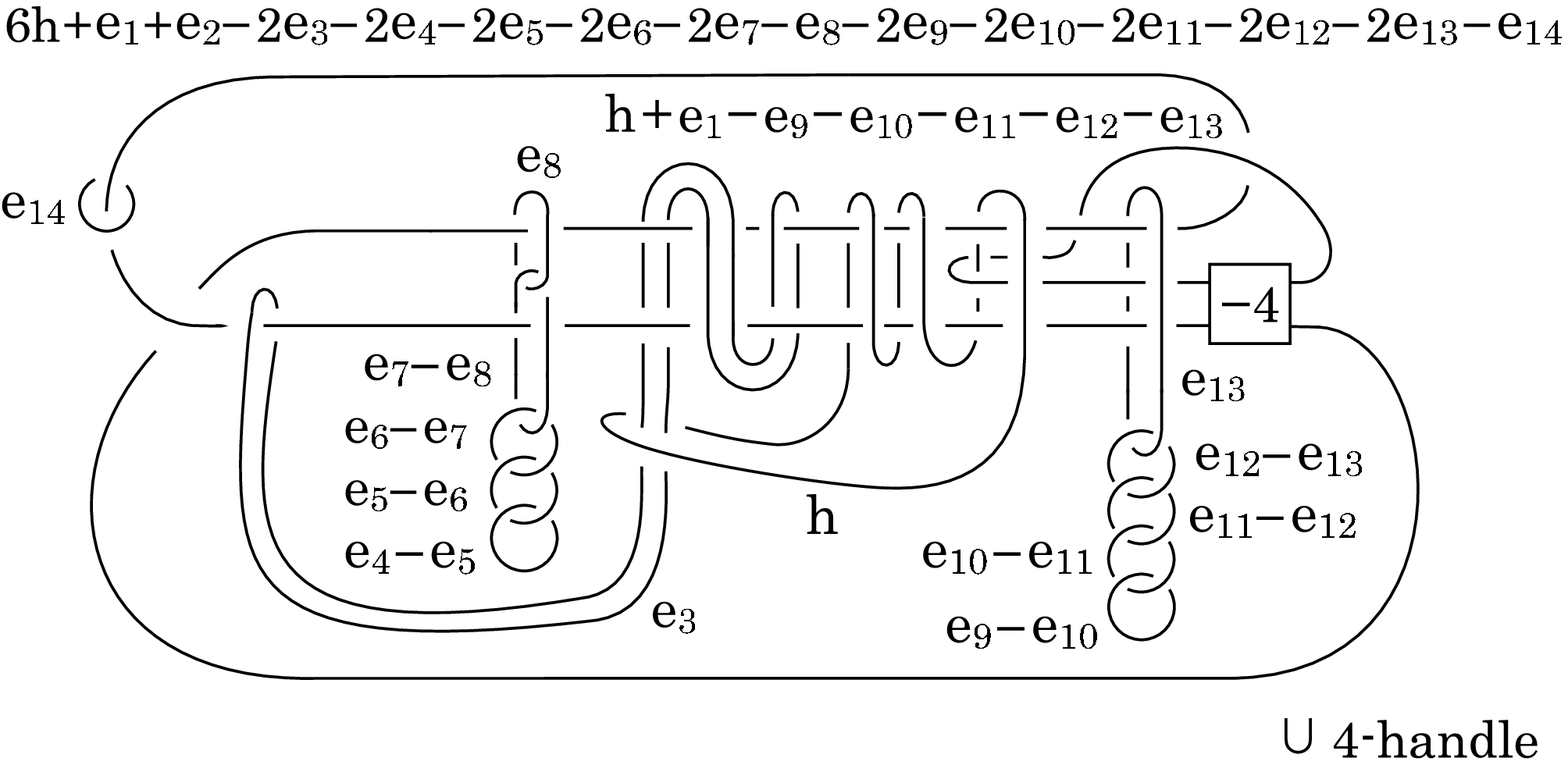}
\caption{$\mathbf{CP}^2\# 14\overline{\mathbf{CP}}^2$}
\label{fig60}
\end{center}
\end{figure}
\begin{figure}
\begin{center}
\includegraphics[width=4.9in]{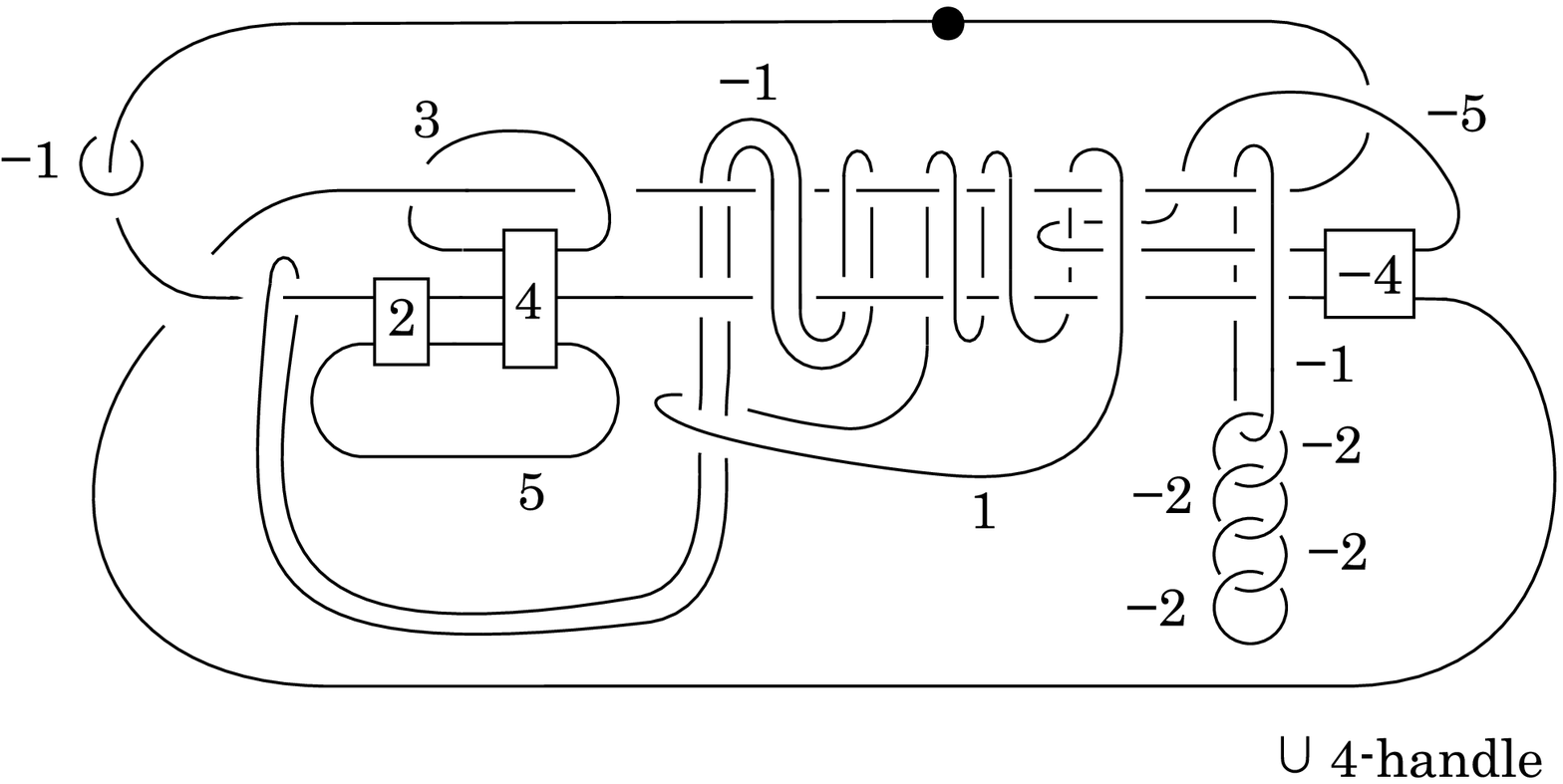}
\caption{A different diagram of the $4$-manifold in Figure~\ref{fig21}}
\label{fig61}
\end{center}
\end{figure}
\end{document}

%% file: goksty.tex

\newcommand{\ben}{\begin{enumerate}}
\newcommand{\een}{\end{enumerate}}
\newcommand{\be}{\begin{equation}}
\newcommand{\ee}{\end{equation}}
\newcommand{\bea}{\begin{eqnarray}}
\newcommand{\eea}{\end{eqnarray}}
\newcommand{\bc}{\begin{center}}
\newcommand{\ec}{\end{center}}

\newtheorem{thm}{Theorem}[section]
\newtheorem{cor}[thm]{Corollary}
\newtheorem{lem}[thm]{Lemma}
\newtheorem{prop}[thm]{Proposition}
\newtheorem{conj}[thm]{Conjecture}

\theoremstyle{definition}
\newtheorem{defn}[thm]{Definition}

\theoremstyle{remark}
\newtheorem{rem}[thm]{\rm\bfseries{Remark}}
\newtheorem*{notation}{Notation}

\newtheorem{ques}[thm]{\rm\bfseries{Question}}
\newtheorem{cons}[thm]{\rm\bfseries{Construction}}
\newtheorem{exm}[thm]{\rm\bfseries{Example}}
